\documentclass[11pt, letterpaper]{article}

\usepackage[]{multibib}
\newcites{A}{References}

\usepackage[final]{microtype}
\usepackage{color}
\usepackage[ruled]{algorithm2e}
\usepackage{subcaption}
\usepackage{varwidth}
\usepackage{capt-of}
\usepackage{boldline}
\usepackage{epsfig}
\usepackage{latexsym}
\usepackage{url}
\usepackage{float}
\usepackage[hypertexnames=false,hyperfootnotes=false]{hyperref}
\usepackage{texnansi}
\usepackage{tikz}
\usepackage{afterpage}
\usepackage{enumerate}
\usepackage[normalem]{ulem}
\usepackage{lmodern}
\usepackage{rotating}
\usepackage[ruled]{algorithm2e}
\usepackage[noend]{algpseudocode}
\usepackage{graphicx}
\usepackage{caption}
\usepackage{subcaption}

\usepackage{booktabs}

\usepackage{sectsty}
\usepackage{ifthen}

\usepackage[leqno]{amsmath}
\usepackage{amsthm}
\usepackage{amssymb}
\usepackage{amsfonts}
\usepackage[margin=10pt,font=small,labelfont=bf]{caption}

\usepackage{nicematrix}

\usepackage{enumitem}
\setlist[enumerate]{topsep=1pt,itemsep=1pt,partopsep=1pt,parsep=1pt}

\usepackage{natbib}
 \bibpunct[, ]{(}{)}{,}{a}{}{,}%
 \def\bibfont{\small}%

\makeatletter
\g@addto@macro{\UrlBreaks}{\UrlOrds}
\makeatother

\newtheoremstyle{thm-sf}{}{}{\itshape}{}{\sffamily\bfseries}{.}{ }{}
\theoremstyle{thm-sf}

\newtheorem{theorem}{Theorem}

\newtheorem{corollary}{Corollary}
\newtheorem{lemma}{Lemma}

\newtheorem{proposition}{Proposition}
\newtheorem{observation}{Observation}

\newtheoremstyle{thm-df}{}{}{}{}{\sffamily\bfseries}{.}{ }{}
\theoremstyle{thm-df}
\newtheorem{definition}{Definition}
\newtheorem{model}{Model}

\renewcommand{\proofname}[1]{{\normalfont\sffamily\bfseries #1}}

\newenvironment{myproof}[1][\proofname]{%
  \proof[\normalfont\sffamily\bfseries #1]%
}{\endproof}

\usetikzlibrary{arrows,patterns,plotmarks}
\tikzstyle{every picture} += [>=stealth]
\allsectionsfont{\sffamily}

\makeatletter
\def\@seccntformat#1{\csname the#1\endcsname.\quad}
\makeatother

\floatstyle{ruled}
\provideboolean{fastcompile}
\newcommand{\hidefastcompile}[1]{\ifthenelse{\boolean{fastcompile}}{}{#1}}

\usepackage[margin=1in]{geometry}
\usepackage[nodisplayskipstretch]{setspace}
\newtheorem*{le:ell}{Lemma~\ref{le:ell}}
\newtheorem*{le:feas}{Lemma~\ref{le:feas}}
\newtheorem*{le:sample_complexity}{Lemma~\ref{le:sample_complexity}}

\usepackage{setspace}

\title{\textsf{\textbf{Near-Optimal Real-Time Personalization \\ with Simple Transformers}}}
\author{Lin An \\ Tepper School of Business \\ Carnegie Mellon University \\ email: \url{linan@andrew.cmu.edu}
\and Andrew A. Li \\ Tepper School of Business \\ Carnegie Mellon University \\ email: \url{aali@andrew.cmu.edu}
\and Vaisnavi Nemala \\ Information Networking Institute \\ Carnegie Mellon University \\ email: \url{vnemala@andrew.cmu.edu}\and Gabriel Visotsky \\ Department of Mathematical Sciences \\ Carnegie Mellon University \\ email: \url{gvisotsk@andrew.cmu.edu}}

\begin{document}
\maketitle






\begin{abstract}
\noindent Real-time personalization has advanced significantly in recent years, with platforms utilizing machine learning models to predict user preferences based on rich behavioral data on each individual user. Traditional approaches usually rely on embedding-based machine learning models to capture user preferences, and then reduce the final real-time optimization task to one of {\em nearest-neighbors}, which can be performed extremely fast both theoretically and practically. However, these models struggle to capture some complex user behaviors, which are essential for making accurate recommendations. Transformer-based models, on the other hand, are known for their practical ability to model sequential behaviors, and hence have been intensively used in personalization recently to overcome these limitations. However, optimizing recommendations under transformer-based models is challenging due to their complicated architectures. In this paper, we address this challenge by considering a specific class of transformers, showing its ability to represent complex user preferences, and developing efficient algorithms for real-time personalization.

We focus on a particular set of transformers, called {\em simple transformers}, which contain a single self-attention layer. We show that simple transformers are capable of capturing complex user preferences, such as variety effects, complementarity and substitution effects, and irrational choice behaviors, which traditional embedding-based models cannot capture. We then develop an efficient algorithm for real-time personalization under simple transformer models. Our algorithm achieves near-optimal performance with sub-linear runtime with respect to the size of the item pool. Finally, we demonstrate the effectiveness of our approach through an empirical study on large datasets from Spotify and Trivago. Our experiment results show that (1) simple transformers predict user preferences substantially more accurately than non-transformer baselines and nearly as accurately as deeper transformer models, and (2) our algorithm completes personalized recommendation tasks both quickly and effectively. Specifically, under a fixed candidate budget, our method achieves objective values that are, on average, 20.86\% higher than those obtained using $k$-Nearest Neighbor and 20.56\% higher than those from Beam Search.

\end{abstract}



\setstretch{1.3}

\section{Introduction}
Personalization today is already immensely sophisticated. Media platforms, online retailers, and subscription services (just to name a few) capture rich data on their users in the form of their behavior and interactions with individual items/products. There are then two key ingredients: (1) this data is used {\em offline} to build machine-learning (ML) based models of users' preferences, and then (2) these models are used  in {\em real-time} to make personalized recommendations.

Zooming out from personalization for just a moment, the most jarring improvements in ML models over the last few years have been in generative models for language, and specifically {\em transformer}-based models (e.g.~the ``T'' in {\em ChatGPT}) that have proven to be extremely accurate in modeling {\em sequential} data. Perhaps unsurprisingly, these same models are well-equipped for personalization. To fix a concrete example, suppose an {\em Instacart} user is in the process of shopping online for groceries. This user’s behavior consists of interactions with grocery items: browsing through items, viewing a subset of these in more detail, and adding a subset of these to their shopping cart. The task of learning this user’s preferences essentially amounts to predicting their future interactions. The important observation here is that the user’s behavior is naturally sequential, and so this prediction task is similar to completing a sentence, where the ``words'' are the items themselves. This connection to language suggests that the same transformer-based models may succeed in learning preferences.

This is already being done in practice, often with substantial empirical success (e.g.~by
  {\em Alibaba} \citep{chen2019behavior}, {\em Amazon} \citep{lake2019large}, {\em Spotify} \citep{moor2023exploiting}, and {\em Wayfair} \citep{mei2022lightweight}). 
  However, as examples of such successes become increasingly common, there is little principled guidance on how transformers ``should'' be used for real-time personalization. This is the problem we seek to address.

\paragraph{Real-Time Personalization, Before Transformers:} To make the nature of this problem more precise, it is worth reviewing how real-time personalization is performed {\em without} transformers. Referring to the two key ingredients mentioned at the outset: first, the ML-based models of user preferences are, by and large, pure {\em embedding}-based models. Using past data, each item $i$ is mapped to some element $v_i \in \mathbb{R}^d$ in such a way that (a) ``similar'' items are ``close'' together, and perhaps more formally, (b) each user can be represented as some $u \in \mathbb{R}^d$ so that the inner products $v_i^\top u$ fully represent the preference/affinity of the user for each item $i$. 

These pure embedding models are not necessarily the most {\em accurate} models that can be estimated from past data, but they enable {\em fast} execution of the second ingredient, which is to optimize a set (or sequence) of items for each user in real time:
\begin{eqnarray} \label{eqn:pre-transformer}
         &\max&  f(S,u)   \\ 
        &\text{s.t.}& S\subset [n], |S|\leq k. \notag
\end{eqnarray}
The (extremely general) formulation above simply highlights that real-time personalization consists of solving cardinality-constrained {\em set-optimization} (or {\em sequence-optimization}, whose equivalence to set-optimization we will discuss later on) problems where (a) the objective function is user-dependent, (b) the number of items $n$ is potentially quite large -- often in the hundreds of millions -- and (c) a solution must be found in real-time, often just milliseconds. This is of course hopeless in general, but feasible when the objective function $f(S,u)$ results from a pure embedding models. In particular, $f(S,u)$ typically takes one of two forms:
\begin{enumerate}
    \item An additive function: \[f(S,u) = \sum_{i \in S} g(v_i^\top u),\]
    for some non-decreasing function $g(\cdot)$. In this case, the optimal set $S^*$ can be characterized as follows: if the number of items $i$ such that $g(v_i^\top u) > 0 $ is no greater than $k$, then $S^*$ consists of all such items. Otherwise, $S^*$ is the set of $k$ items with the largest values of $g(v_i^\top u)$.
    \item A monotone submodular function (in the argument $S$, for any $u$), such that each item is fully encoded by $v_i^\top u$, so that a constant-factor approximation can be found using  greedy-style algorithms along with (multiple queries to) a black box which computes the item with largest inner product to a given point in $\mathbb{R}^d$ \citep{farias2020optimizing}.
\end{enumerate}
In both of the above cases, the pure embedding model essentially reduces the final real-time optimization task to one of {\em nearest neighbor} algorithms, where the goal is to find the items whose embeddings are most similar (under inner product similarity) to the user embedding. This task can be performed extremely fast, both theoretically (using \textit{approximate nearest neighbor} algorithms with runtime sub-linear in $n$) and practically (given the nonstop engineering and improvement of commercial vector databases).

\paragraph{An Attempt to Introduce Transformers:} Returning to transformers now, the natural opportunity is to improve the accuracy of the pure embedding models in representing user preferences: the embedding models essentially fail to capture the effects that items may have on each other when present in the same recommended set. Such effects are well-known to exist in multiple fields of study, as we will discuss shortly, and often representable via transformers. Unfortunately, the just-described synergy between the ``upstream'' user preference model estimated from data, and the ``downstream'' optimization of user-dependent sets of items completely breaks down here. As we will see in the next section, this is true in a formal sense: Problem (\ref{eqn:pre-transformer}) is NP-hard, and likely hard to even approximate in linear time, if the objective $f(\cdot,u)$ is transformer-based. As of now, any practical implementation using transformers either applies a pure nearest-neighbor or greedy-style algorithm (essentially ignoring this hardness), or a random search heuristic (such as {\em Beam Search}).

In the spirit of developing a principled approach to transformer-based real-time personalization, this raises two major questions:
\begin{enumerate}
    \item {\bf Fast Optimization:} The formal hardness results just alluded to imply that achieving fast (ideally sub-linear in $n$) optimization of Problem (\ref{eqn:pre-transformer}) with any meaningful optimality guarantee is impossible if the objective $f(\cdot,u)$ is to allow for {\em all} transformers. Naturally then, is there a non-trivial sub-class of transformers for which this is possible?
    \item {\bf Modeling Power:} Assuming a positive answer to the first question, i.e.~assuming the existence of a subset of transformers which enable fast optimization, does restricting to this subset come at a substantial cost in terms of modeling user preferences? Put another way, can this smaller sub-class of transformers achieve the same predictive accuracy as the family of all transformers?
\end{enumerate}






\subsection{Our Contributions}
In short, our contributions provide concrete, theoretically-backed answers to both of the above questions:

\paragraph{Modeling User Preferences with Simple Transformers.} 
We focus our study on a sub-class of transformers that we refer to as {\em simple transformers}. These are transformers which contain a single {\em self-attention} layer (to be defined in the next section), whereas transformers as a whole may contain multiple attention layers. Addressing the pair of questions above in reverse,  we first formally show that simple transformers are able to represent two known, popular parametric models of user preference (Proposition \ref{prop: variety and subsitution}):

    \begin{itemize}
        \item {\em Sequential variety effects} in the context of marketing;
        \item {\em Pairwise complementarity and substitution effects} in the context of economics.
    \end{itemize}
It should also be emphasized that none of these models are representable via pure embedding models.

\paragraph{Real-Time Personalization with Simple Transformers.} Our main result (Theorem \ref{thm: main theorem}) is an algorithm (Algorithm \ref{alg: complete}) which approximately solves Problem (\ref{eqn:pre-transformer}) in sub-linear time, when the objective function is given by a simple transformer: 
    
\begin{theorem}[Informal]
Under additional (rank) assumptions on the simple transformer, given any $n,k\in\mathbb{N}$ and $\epsilon>0$, there exists an algorithm that achieves $\textup{ALG}\geq(1-\epsilon)\textup{OPT}$  with expected amortized runtime \[\tilde{O}\left(n^{1-c(\epsilon,k)}\cdot k^{\mu(\epsilon)}\right)\] for functions $c,\mu$ satisfying $c(\epsilon,k),\mu(\epsilon)>0$. Here $\tilde{O}$ hides factors of order $n^{o(1)}$.
\end{theorem}

We will present and discuss the rank assumptions in Section \ref{Section Hardness}. We also show that these assumptions are necessary for approximately solving Problem (\ref{eqn:pre-transformer}) in sublinear time (Proposition~\ref{prop: lower bound} and Proposition~\ref{prop: MDK lower bound}), and that our algorithm's expected amortized runtime dependence on $k$ and $\epsilon$ is optimal (Proposition~\ref{prop: MDK lower bound}).

Our algorithm operates under the same two-phase {\em retrieval} and {\em ranking} paradigm that is used in many competition-winning personalization algorithms, though in our case both phases are adapted specifically to simple transformers, and enjoy provable guarantees (the combination of which generates our main result). Our algorithm has the added practical benefit of subsuming (given a particular, sub-optimal selection of tuning parameters) the {\em Beam Search} algorithm commonly used in practice.

\paragraph{Empirical Study.}

We empirically validated the theoretical results of the previous contributions on two large datasets from {\em Spotify}  \citep{chen2018recsys} (which includes $1,000,000$ playlists with $2,262,292$ unique songs) and the travel website {\em Trivago} \citep{Knees_etal:RSC:2019} (which includes user sessions of searching for hotel bookings, with around $730,000$ unique users and around $340,000$ unique hotels recorded in around $900,000$ different sessions). 

In support of the first contribution, our first set of experiments demonstrates that, given data on past user behavior, simple transformers can effectively model and predict user preferences. They achieve (a) substantially higher accuracy than non-attention-based models (such as pure embedding models), and (b) performance that is nearly comparable to more complex transformer architectures. Specifically, simple transformers achieved, on average, 14.1\% higher accuracy than non-attention models (e.g., logistic regression, random forest, support vector machine), and only 2.5\% lower accuracy than general transformers with multiple self-attention layers.

In support of the second contribution, our second set of experiments demonstrates that our algorithm performs simple-transformer-based personalized recommendation tasks both efficiently and effectively. We solved instances of Problem (\ref{eqn:pre-transformer}) using the simple transformers trained in the first set of experiments and compared our algorithm to two widely used benchmark methods: $k$-Nearest Neighbor (for retrieval) and Beam Search (for ranking). Under a fixed candidate solution budget (to partially standardize runtime), our algorithm achieved objective values that were, on average, 20.86\% higher than those obtained using $k$-Nearest Neighbor and 20.56\% higher than those obtained using Beam Search. Therefore our algorithm achieved strong empirical performance in both the retrieval and ranking phases.

\subsection{Literature Review}

\paragraph{Transformers in Recommender Systems.} Recommender systems have significantly evolved with the emergence of transformer-based architectures, first introduced by \cite{vaswani2017attention}. Transformers' ability to model long-range dependencies and efficiently process sequential data makes them particularly well-suited for capturing user-item interaction sequences. They have been successfully adopted in practice by companies such as Alibaba \citep{chen2019behavior}, Amazon \citep{lake2019large}, Spotify \citep{moor2023exploiting}, and Wayfair \citep{mei2022lightweight}, among others. A large body of work has focused on designing specialized transformer-based architectures for recommendation tasks. These include self-attention-based sequential models like SASRec \citep{kang2018self,mei2022lightweight,wilm2023scaling}, single-layer attention models \citep{wang2018attention,chen2019behavior,bendada2023track,celikik2022reusable}, multi-head or multi-layer self-attention models \citep{yang2023going,zheng2023generative}, neural attention mechanisms \citep{chen2017attentive,fu2018attention,lake2019large}, recurrent attention models \citep{sukhbaatar2015end}, and sparse attention mechanisms \citep{li2023strec}. In contrast, our work focuses on the optimization task that follows once the transformer architecture has been established. Specifically, we consider architectures with a single attention layer--an approach that is both prominent and empirically successful in recommender systems--and aim to make fast and near-optimal recommendations based on these models.

\paragraph{Representational Power of Transformers.} Transformer architectures are built on the self-attention mechanism, which endows them with strong representational power in both theory and practice. Having discussed practical aspects above, we now review theoretical results. \citet{yun2019transformers} and \citet{wei2022statistically} proved universal approximation results, showing that sufficiently large transformers can approximate broad classes of functions, analogous to results for feedforward networks \citep{hornik1989multilayer}. \citet{perez2019turing} and \citet{wei2022statistically} further established the (approximate) Turing-completeness of transformers. \citet{sanford2024representational} gave both positive and negative results: they introduced a sparse averaging task where transformers scale logarithmically with input size (unlike recurrent or feedforward networks, which scale polynomially), but also a triple detection task where attention scales linearly. Similar negative results were shown for induction heads by \citet{sanford2024one,bietti2024birth,elhage2021mathematical}.

For recommender systems, the most relevant work connects representational power to choice modeling. \citet{ko2023modeling} showed that classic choice models such as Halo-MNL can be represented by a single attention layer, and \citet{wang2023transformer} developed a transformer architecture for learning and predicting many choice models. In our work, we demonstrate that a single attention layer can also capture user preference models involving sequential variety effects and pairwise complementarity/substitution.

\paragraph{Approximate Nearest Neighbor.} Approximate nearest neighbor (ANN) is a key tool for real-time personalization, where a system must quickly retrieve relevant items from embeddings. The classical nearest neighbor (NN) search is computationally intractable at scale, especially in high dimensions, motivating efficient approximation algorithms. ANN returns near-optimal results with much lower query time and storage, making it well-suited for real-time personalization.

Several algorithmic frameworks exist. Hashing-based methods, such as Locality-Sensitive Hashing (LSH) \citep{indyk1998approximate,andoni2008near,andoni2015practical}, provide sub-linear query time under specific metrics. Graph-based methods construct proximity graphs (e.g., navigable small-world graphs) for efficient traversal \citep{malkov2018efficient}. Tree-based methods such as KD-trees and Ball Trees are effective in low dimensions but deteriorate as dimensionality grows \citep{muja2014scalable}. In our work, ANN enables fast identification of items likely to interest a user, allowing our algorithm to run in real time.

\paragraph{Binary Quadratic Optimization.} With transformer architectures, the recommendation task can be cast as a binary quadratic optimization problem (quadratic knapsack), which is NP-hard and subsumes many difficult problems, including maximum clique and densest $k$-subgraph. Moreover, it is NP-hard to approximate within any finite factor \citep{rader2002quadratic}. Much of the literature therefore studies tractable special cases: \citet{rader2002quadratic} gave an FPTAS for series-parallel graphs, while \citet{taylor2016approximation} developed an FPTAS for bounded-treewidth graphs and a PTAS for planar graphs. For surveys, see \citet{pisinger1998knapsack,cacchiani2022knapsack}.

In our setting, we exploit the low non-negative rank of the softmax matrix to obtain a PTAS. Related work on low-rank optimization includes: an FPTAS for minimizing quasi-concave functions over convex sets \citep{goyal2013fptas}, an FPTAS for low-rank functions over polytopes \citep{mittal2013fptas}, and a PTAS for binary non-linear programs with low-rank objectives \citep{nguyen2021ptas}. Our problem differs in requiring low \emph{non-negative} rank, for which we refer to \citet{cohen1993nonnegative}. Owing to the specific structure of the softmax matrix, prior results do not directly apply. Our approach is also related to the multi-objective knapsack problem, where FPTAS algorithms compute the Pareto frontier \citep{elhage2021mathematical,bazgan2009implementing,bazgan2009solving}.

\section{Model}
\subsection{Real-Time Personalization without Transformers} \label{Section: before transformers}

\paragraph{Pure Embedding Models.}
Before introducing the transformer-based models that will be the focus of this paper, we begin with a brief overview of \emph{pure embedding models}, which are widely used in modern personalization systems. These models represent users and items as vectors in a shared low-dimensional space, enabling efficient modeling of interactions via simple operations such as the inner product.

Formally, let $[n] \equiv {1, \ldots, n}$ denote a set of items, and let $V \in \mathbb{R}^{n \times d}$ be a matrix whose $i$-th row $v_i^\top \in \mathbb{R}^d$ is the \emph{value vector} of item $i$. These value vectors are designed so that items with similar embeddings are likely to be perceived similarly by users. Each user is likewise represented by a \emph{user vector} $u \in \mathbb{R}^d$. Both the item and user vectors are typically learned from historical interaction data -- through matrix factorization, collaborative filtering, or more complex models trained on click or engagement feedback.

The utility (or reward) of recommending item $i$ to user $u$ is modeled as $f_i(v_i^\top u)$, where $f_i : \mathbb{R} \to \mathbb{R}$ is a non-decreasing reward function specific to item $i$. In many applications, the same function $f$ is used for all items. However, in some cases it is important to allow heterogeneous reward functions that reflect item- or category-specific behavior. For example, different product categories often exhibit different click-through-rate (CTR) saturation patterns: a small increase in relevance (e.g., $v_i^\top u$) might sharply increase CTR for breaking news articles, whereas product ads might exhibit more gradual, linear gains, and fashion items may plateau early due to browsing behavior. Modeling such differences requires different shapes of the function $f_i$, even if all are monotonic. To maintain full generality, we therefore allow each item to have its own reward function $f_i$. Moreover, in many applications the image of \( f_i \) is often \([0,1]\), in which case \( f_i(v_i^\top u) \) can be interpreted as the probability of a positive user action, such as a click or purchase. In these cases, \( f_i \) is sometimes chosen to resemble functions like the logistic function.

One of the key advantages of pure embedding models is the efficiency of real-time personalization. Given a user vector $u$ and a budget $k$, the goal is to select a subset $S \subset [n]$ of at most $k$ items that maximizes the total reward:
\begin{align}
\label{eqn:pure embedding}
\tag{$\mathsf{Pure \text{ } Embedding}$}
\max \quad & \sum_{i \in S} f_i(v_i^\top u) \\
\text{s.t.} \quad & S \subset [n],\ |S| \leq k. \notag
\end{align}

Since the objective function is additive and the items are treated independently, this problem can be solved exactly via \emph{nearest neighbor}. In our setting, this corresponds to identifying the $k$ items whose value vectors are most aligned with the user vector $u$ under inner product similarity. That is, for a query $u \in \mathbb{R}^d$, the goal is to find the top-$k$ items maximizing $v_i^\top u$.

Nearest neighbor arises in many applications across recommendation, vision, and language, where one often needs to retrieve items similar to a given input based on some feature representation. A naive solution evaluates all $n$ inner products $v_i^\top u$, which takes $O(nd)$ time. In our case, we compute $f_i(v_i^\top u)$ for every $i \in [n]$, and then select the top $k$ items with the largest values. The optimal solution $S^*$ can be described as follows: if there are at most $k$ items with positive reward (i.e., $f_i(v_i^\top u) > 0$), then $S^*$ includes all of them. Otherwise, $S^*$ consists of the $k$ items with the largest values of $f_i(v_i^\top u)$. This greedy procedure yields an exact solution in linear time with respect to the total number of items.

\paragraph{Algorithms: Approximate Nearest Neighbor.}
The greedy implementation of nearest neighbor becomes computationally prohibitive as the number of points $n$ grows. In most applications such as personalization tasks, the number of items $n$ is on the scale of millions, and an algorithm with runtime linear in $n$ cannot be performed in real-time. To address this problem, the notion of \textit{approximate nearest neighbor} (ANN) has been widely adopted. Rather than finding the exact nearest point, ANN aims to return a point whose distance to the query is within a factor of the true minimum. Formally, for an additive approximation factor $\epsilon > 0$, the objective of $\epsilon$-ANN is to find any point $v_{i^*}$, where $i^*\in [n]$, that approximately maximizes the inner product similarity to the query $u$:
\[
v_{i^*}^\top u \geq \arg\max_{i \in [n]} v_{i}^\top u-\epsilon.
\] 
This relaxation enables much more efficient data structures and algorithms -- often sub-linear in $n$ -- which can be implemented in real-time. \footnote{There are also other versions of $\epsilon$-ANN that concerns multiplicative approximation factors instead of additive ones. These two notions are equivalent under some boundedness assumptions on the points.} 
Various techniques have been applied to $\epsilon$-ANN algorithms. For example, \cite{andoni2015practical} gives an $\epsilon$-ANN algorithm by using Locality-Sensitive Hashing: 
\begin{proposition}[Corollary 1 in \cite{andoni2015practical}] 
    For any given $\epsilon>0$, $\epsilon$-ANN on the unit sphere $\mathbb{S}^{d-1}\subset\mathbb{R}^d$ can be solved with expected amortized runtime $O(dn^{\gamma(\epsilon)})$, where $\gamma(\epsilon)=\frac{1}{1+c\epsilon}+o(1)$, where $c>0$ is a universal constant.
\end{proposition} 
\cite{ailon2009fast} gives an $\epsilon$-ANN algorithm by applying fast Johnson-Lindenstrauss transform (FJLT):
\begin{proposition}[Theorem 1 in \cite{ailon2009fast}]
    For any given $\epsilon>0$, $\epsilon$-ANN on the unit sphere $\mathbb{S}^{d-1}\subset\mathbb{R}^d$ can be solved with expected amortized runtime  $O\left(d\log(d)+\epsilon^{-3}\log^2(n)\right)$.
\end{proposition}
\cite{arya1998optimal} gives an $\epsilon$-ANN algorithm by using a tree-based data structure:
\begin{proposition}[Theorem 3.1 in \cite{arya1998optimal}]
    For any given $\epsilon>0$, the $\epsilon$-ANN on the unit sphere $\mathbb{S}^{d-1}\subset\mathbb{R}^d$ can be solved with expected amortized runtime  $O(\log n + 1/\epsilon^d)$.
\end{proposition}

Because companies aim to recommend a set of at most $k$ items to a user in their personalization task, their objective is not merely to identify a single item that is attractive to the user, but rather to efficiently retrieve a set of $k$ items that are collectively among the most attractive to the user. Therefore, they consider the notion of \textit{$\epsilon$-Approximate $k$-Nearest Neighbor}, which returns a set of $k$ items whose inner product similarities are within an additive error of $\epsilon$ compared to the top-$k$ true nearest items. This is formally defined below.

\begin{definition}[$\epsilon$-Approximate $k$-Nearest Neighbor Algorithm]
    An $\epsilon$-Approximate $k$-Nearest Neighbor algorithm builds a data structure on any given set of points $\{v_1,\dots,v_n\}\subset\mathbb{R}^d$, and takes any query $u\in \mathbb{R}^{d}$, any $1\leq k\leq n$, and any $\epsilon>0$ as inputs. Let $\pi:[n]\to [n]$ be a permutation of the indices such that $v_{\pi(1)}^\top u\geq\cdots\geq v_{\pi(n)}^\top u$. Let $(i^*_1,\dots,i^*_k)=(\pi(1),\dots,\pi(k))$. The oracle outputs $k$ indices $i_1,\dots,i_k\in [n]$ such that $v_{i_j}^\top u\geq v_{i^*_j}^\top u-\epsilon$ for each $j=1,\dots,k$ with expected amortized runtime $k\textup{-ANN}(n,d,k,\epsilon)$.
\end{definition}


Many $\epsilon$-ANN algorithms can be naturally modified to $\epsilon$-Approximate $k$-Nearest Neighbor algorithms, with a similar expected amortized runtime, that is, sub-linear in $n$. As an example, we give an $\epsilon$-Approximate $k$-Nearest Neighbor algorithm in Appendix \ref{Appendix A}.
\begin{lemma}\label{lemma: nearest neighbor}
    There exists an $\epsilon$-Approximate $k$-Nearest Neighbor algorithm, which we give in Appendix \ref{Appendix A}, with expected amortized runtime $$k\textup{-ANN}(n,d,k,\epsilon)=O\left(k(d\log(d)+\epsilon^{-3}\log^2(n))\right).$$ 
\end{lemma}

Finally, companies apply any given $\epsilon$-Approximate $k$-Nearest Neighbor algorithm to approximately solve Problem~\eqref{eqn:pure embedding}. This can be done by first partitioning the items according to their reward functions \( f_i \). For each partition, they apply an \(\epsilon\)-Approximate \(k\)-Nearest Neighbor algorithm to identify \(k\) items that are collectively among the most attractive to the user within that partition. They then evaluate the rewards of all such candidate items across partitions and select the top-\(k\) items with the highest overall reward.

To analyze the approximation error incurred by this procedure, we introduce the following parametrization of the reward functions. For all \(i \in [n]\) and \(\epsilon > 0\), the function \(f_i\) satisfies:
\[
f_i(x - \epsilon) \geq (1 - g(\epsilon)) f_i(x) - h(\epsilon) \quad \text{for all } x,
\]
where \(0 \leq g(\epsilon) \leq 1\) and \(h(\epsilon) \geq 0\) are non-negative functions of \(\epsilon\). This parametrization captures the idea that a small additive error in the input \(x\) leads to a controlled multiplicative and additive error in the output \(f_i(x)\).

This parametrization captures the behavior of many reward functions commonly used in practice. In particular, we make the following observation:
\begin{observation} \label{obs: g and h}
    If \(f_i(\cdot)\) is non-decreasing and \(L\)-Lipschitz, then it satisfies the above condition with \(g(\epsilon) = 0\) and \(h(\epsilon) = L\epsilon\). In particular, many standard reward functions (also referred to as activation functions in machine learning), such as logistic, ReLU, leaky ReLU, PReLU, tanh, and softplus, satisfy this condition with \(g(\epsilon) = 0\) and \(h(\epsilon) = O(\epsilon)\).
\end{observation}

With this parametrization, we can now state our main result on applying ANN algorithms to real-time personalization with pure embedding models. We begin by introducing notations that will be used throughout the paper. For an optimization problem $\mathsf{P}$ with objective function $f_{\mathsf{P}}$, let $x^*_{\mathsf{P}}$ denote an optimal solution and define
$\textup{OPT}_{\mathsf{P}} = f_{\mathsf{P}}(x^*_{\mathsf{P}})$
as its optimal value. For an algorithm $\textup{ALG}$ applied to $\mathsf{P}$, let $x^{\textup{ALG}}_{\mathsf{P}}$ be the solution returned by $\textup{ALG}$, and set $\textup{ALG}_{\mathsf{P}} = f_{\mathsf{P}}(x^{\textup{ALG}}_{\mathsf{P}})$
to be the corresponding objective value.

\begin{proposition} \label{prop: pure embedding}
    Suppose we have an $\epsilon$-Approximate $k$-Nearest Neighbor  algorithm with expected amortized runtime $k\textup{-ANN}(n,d,k,\epsilon)$, and suppose $k\textup{-ANN}(n,d,k, \epsilon)$ is concave in $n$. Let $\tau
    $ be the number of distinct functions among $f_1,\dots,f_n$. Given any $\epsilon>0$, we give an algorithm $\textup{ALG} $ for solving Problem \eqref{eqn:pure embedding} that satisfies \begin{equation*}
            \textup{ALG}_{\mathsf{Pure \text{ }Embedding}}\geq (1-g(\epsilon))\textup{OPT}_{\mathsf{Pure \text{ }Embedding}}-kh(\epsilon).
        \end{equation*} with expected amortized runtime $$\tau\cdot k\text{-ANN}\left(\frac{n}{\tau},d,k,\epsilon\right)+\log_2(\tau k).$$
\end{proposition}

The proof of Proposition \ref{prop: pure embedding} appears in Appendix \ref{Appendix A}. Notice that many $\epsilon$-Approximate $k$-Nearest Neighbor algorithms have $k\textup{-ANN}\left(n,d,k, \epsilon\right)$ sub-linear and concave in $n$, and we have given an example in Lemma \ref{lemma: nearest neighbor}. Therefore Problem \eqref{eqn:pure embedding} can be approximately solved in sub-linear time, which is practical to implemented in real-time.

\paragraph{Modeling Limitations.}

Proposition \ref{prop: pure embedding} shows that pure embedding models can be optimized efficiently. However, pure embedding models have a fundamental limitation: they treat each item independently and fail to capture how the value of an item may depend on the context in which it is presented. That is, the reward associated with item $i$ is fixed once $u$ and $v_i$ are known, regardless of what other items are presented alongside it. In many personalization applications, this is an unrealistic assumption. For example, the value of a song recommendation may drop if a similar song is already in the playlist; or the likelihood a user clicks on a hotel result may depend on what other hotels appear in the same search result page and how they compare. These effects, often referred to as \emph{set effects}, are not captured by pure embedding models. Concretely, we present three common parametric models used in personalization that cannot be represented by pure embedding models.

First, we consider a famous parametric model of \textit{sequencial variety effects} in sequences. The concept of variety/diversity has been examined extensively in the marketing literature (see e.g. \cite{mcalister1982dynamic,hoch1999variety,rafieian2023optimizing}). This model proposes that the perceived utility of an item depends not only on its intrinsic quality but also on its novelty relative to previously seen items. In particular, repeated exposure to similar items leads to diminishing marginal utility, while introducing diverse or contrasting items can restore or amplify engagement. This intuition aligns closely with the notion of \emph{discounted utility} over sequences, where the utility derived from an item is multiplicatively reduced based on its similarity to past items (see, e.g. \cite{barbera2004handbook,baucells2007satiation}). Such models capture behavioral tendencies like satiation, boredom, and the desire for exploration, and have been influential in both economic theory and practical recommendation systems. Below we give a mathematically formulation of sequential variety effects.


\begin{model}[Sequential Variety Effects]\label{mod:variety}
Let $[n]:=\{1,\dots,n\}$ denote the set of items. Each item $i$ has a \emph{base utility} $\hat{u}_i>0$ and a \emph{similarity embedding} $x_i\in\mathbb{S}^{d-1}$. Define \emph{pairwise similarity score} between item $i$ and item $j$ by $s_{ij}:=x_i^\top x_j\in[-1,1]$. Fix a sequence length $k\ge 1$ and nonnegative \emph{lag weights} $\lambda_1,\dots,\lambda_{k-1}\ge 0$.

For a sequence $S=(i_1,\dots,i_k)$ of length $k$, the \emph{sequential–variety–adjusted utility} of the item at position $t$ is
\[
g(S,i_t)
= \hat{u}_{i_t}\,\exp\!\left(
\beta \sum_{\ell=1}^{t-1}\lambda_\ell\, s_{\,i_t,\,i_{t-\ell}}
\right),\qquad t=1,\dots,k,
\]
where $\beta\in\mathbb{R}$ controls the strength and sign of the variety effect. \footnote{By convention, the empty sum equals $0$, so $g(S,i_1)=\hat{u}_{i_1}$.}
\end{model}

Model~\ref{mod:variety} adjusts the context-free utility $\hat{u}_{i_t}$ multiplicatively according to the similarity between the current item and recently shown items, with older influences discounted by $\lambda_\ell$. When $\beta<0$ (preference for variety), similarity to recent items ($s_{ij}>0$) decreases the current utility, while dissimilarity ($s_{ij}<0$) increases it; when $\beta>0$ (preference for continuity), the effect reverses and thematic similarity boosts utility. The weights $\lambda_\ell$ specify the memory profile: choosing $\lambda_\ell=\rho^\ell$ with $\rho\in(0,1)$ yields exponential decay in influence with lag, whereas setting $\lambda_1>0$ and $\lambda_\ell=0$ for $\ell>1$ recovers a one-step effect based only on the immediate predecessor. The exponential factor ensures $g(S,i_t)>0$ and provides a smooth, multiplicative adjustment driven by recent sequence context.

Second, we consider a model of \emph{pairwise complementarity and substitution effects}.
In economics, these effects describe how the presence of one item influences the desirability of another: complementary items enhance each other’s attractiveness, while substitutes reduce it.
A common modeling approach represents items as vectors in a feature space, with pairwise interactions captured via inner products (see, e.g., \cite{lee2013direct,berry2014structural,ruiz2020shopper}).

\begin{model}[Pairwise Complementarity and Substitution Effects] \label{mod:com/sub} 
Let $[n]$ denote the set of items. Each item $i$ has a value vector $\hat{v}_i\in\mathbb{R}^d$. The pairwise complementarity and substitution effects is parametrized by a matrix $H\in\mathbb{R}^{n\times n}$ where $H_{ii}=0$ for every $i\in[n]$. For a user $\hat{u}\in\mathbb{R}^d$ and a subset of items $S\subset[n]$, the \textit{interaction-adjusted utility} of item $i\in S$ is defined as 
\[g(S,i)=\hat{v}_i^\top \hat{u}+\sum_{j\in S}H_{ij}.\]
\end{model}

In Model~\ref{mod:com/sub}, the value vector $\hat{v}_i$ captures the intrinsic preference of the user $u$ for item $i$, independent of any other items shown. The matrix $H$ encodes all pairwise interaction effects between items: a positive entry $H_{ij}$ means that the presence of item $j$ increases the utility of item $i$ (complementarity), while a negative entry $H_{ij}$ means that the presence of $j$ reduces the utility of $i$ (substitution). The diagonal entries are zero by definition, so an item does not directly influence its own utility. When $H_{ij}=0$, item $j$ has no effect on the utility of item $i$.

Model~\ref{mod:com/sub} is also closely related to \emph{choice models} that generate \emph{conversion probabilities}. Choice models are a fundamental input to many now-canonical optimization problems in the field of operations management, including assortment, inventory, and price optimization. At a high level, a  choice model maps an offer set \(S\subset[n]\) to conversion probabilities \(\{p(i\mid S)\}_{i\in S}\) that sum to one. Interpreting \(g(S,i)\) from Model~\ref{mod:com/sub} as an unnormalized log-weight, a choice model can naturally set purchase probabilities proportional to \(\exp(g(S,i))\). This precisely recovers a particular type of choice model, called the \textit{Halo Multinomial Logit} (Halo MNL) choice model \citep{maragheh2018customer}.

The Halo MNL choice model is a generalization of the simplest and most widely used \textit{multinomial logit}  (MNL) choice model, which assumes that customers choose among a set of items according to a fixed utility associated with each option. While the MNL choice model provides a tractable and interpretable framework, it has well-known limitations, such as the independence of irrelevant alternatives property, which can fail to capture context-dependent or irrational choice behaviors observed in practice. The Halo MNL choice model addresses this by allowing the utility of each item to depend on the presence or absence of other items in the offered set. In particular, certain items can impose a positive or negative ``halo effect'' on the utility of other items, enabling the model to capture behaviors that violate classical rationality assumptions. 

The halo effect is completely parametrized by a matrix $H\in\mathbb{R}^{n\times n}$, where each off-diagonal entry \(H_{ij}\) quantifies the halo effect that the presence of item \(j\) imposes on the utility of item \(i\): positive values capture complementarity (item \(j\) makes \(i\) more attractive), while negative values capture substitution (item \(j\) detracts from \(i\)). Crucially, these interactions are precisely represented by Model~\ref{mod:com/sub}: the interaction-adjusted scores serve as the weights of the conversion probabilities. Normalizing these weights over the offered set (e.g., via a softmax) yields the Halo MNL choice model with parameter $H$. Moreover, when \(H\) is the zero matrix, there are no interaction terms and the normalization reduces to the standard MNL choice model. Therefore the MNL choice model is a special case nested within our framework.



Because in both models the reward associated with an item can depend on the presence or absence of other recommended items, we make the (informal) observation that these models cannot be represented by pure embedding models. In contrast, we will later show that both models can be naturally expressed within the framework we study in this paper.

\subsection{Simple Transformers}

In this section, we formally state the model which will be our primary object of study: \textit{simple transformers}, or neural networks with a single {\em self-attention} layer. First, some preliminary definitions: the row-wise {\em softmax} operator, which we denote as $\text{softmax}:\mathbb{R}^{n\times d}\to \mathbb{R}^{n\times d}$,  is given by \begin{equation*}
     \text{softmax}(A)_{i,j}= \frac{\exp(A_{i,j})}{\sum_{j'=1}^d \exp(A_{i,j'})}.
\end{equation*} In particular, each row of $\text{softmax}(A)$ sums up to 1, and can be interpreted as a vector of weights.
The notion of attention is fundamental to transformer-based models (e.g. \cite{vaswani2017attention,sanford2024representational}). For input dimension $n$, output dimension $d_v$, embedding dimension $d_{kq}$, and matrices $Q,K\in \mathbb{R}^{n\times{d_{kq}}}$, and $V\in \mathbb{R}^{n\times{d_{v}}}$, a \textit{self-attention layer} is a  function $\mathrm{SA}_{Q,K,V}:\mathbb{R}^{n\times n}\to \mathbb{R}^{n\times{d_{v}}}$ given by
\begin{equation*}
\mathrm{SA}_{Q,K,V}(X) = \mathrm{softmax}((XQ) (XK)^\top) XV.
\end{equation*}
Here, the matrices $Q,K$, and $V$ are often called the {\em query}, {\em key}, and {\em value} matrix, respectively.

To better understand the self-attention layer, it is natural to view it as a function on subsets of items.  
Let $[n]$ denote a set of items. Then the self-attention layer is fully parameterized  
by an individual query, key, and value vector for each item (forming the rows of the respective matrices $Q$, $K$, and $V$).  
For any subset $S \subset [n]$,  
we define the \emph{set-membership matrix} $X_S \in \{0,1\}^{n\times n}$ by  
\[
(X_S)_{ii} = 
\begin{cases}
1, & i \in S,\\
0, & i \notin S.
\end{cases}
\]  
The function $\mathrm{SA}_{Q,K,V}(\cdot)$ is then applied to $X_S$.  
The output is an $n \times d_v$ matrix in which each row $i \in [n]$ is interpreted as follows:  

\begin{itemize}
    \item If $i \in S$, then the $i$-th row is a weighted average of the value vectors $\{v_i : i \in S\}$,  
    where the weight on each $v_i$ is given by the softmax-normalized dot product between the query vector $q_i$ and the key vector $k_j$ for $j \in S$.  

    \item If $i \notin S$, then the $i$-th row is the uniform average of the value vectors $\{v_{i} : i\in S\}$. This is because item $i$'s query vector is zeroed by $X_S$ -- that is, $(X_SQ)_i$ is the zero vector -- and thus assigns equal softmax weight to all items in $S$.  
\end{itemize}

\textit{Transformers} are a broad family of functions (also known as \textit{neural networks}) constructed by repeatedly applying self-attention layers along with simple transformations that operate independently on each item. These point-wise transformations typically consist of linear mappings followed by non-linear functions known as \textit{activation functions}, which introduce flexibility and allow the model to capture complex behaviors. In our context, these activation functions can be naturally interpreted as the reward functions $f_i$, representing how the utility of each item responds to its input.

\textit{Simple transformers} are a subclass of transformers with a single self-attention layer, followed by point-wise activation functions $f_i$. Formally, we define them as follows:

\begin{definition}[Simple Transformer]
For matrices $Q,K \in \mathbb{R}^{n\times d_{kq}}$ and $V \in \mathbb{R}^{n\times d_{v}}$, 
a vector $u \in \mathbb{R}^{d_{v}}$, and non-decreasing \textit{activation functions} 
$f_1,\dots,f_n : \mathbb{R} \to \mathbb{R}$, a \textit{simple transformer} is a function 
$\mathcal{T}_{Q,K,V,f_1,\dots,f_n,u} : \mathbb{R}^{n\times n} \to \mathbb{R}^{n}$ given by
\[
\mathcal{T}_{Q,K,V,f_1,\dots,f_n,u}(X) = 
\begin{bmatrix}
f_{1} \!\left( \mathrm{SA}_{Q,K,V}(X)_{1}^\top u \right) \\
f_{2} \!\left( \mathrm{SA}_{Q,K,V}(X)_{2}^\top u \right) \\
\vdots \\
f_{n} \!\left( \mathrm{SA}_{Q,K,V}(X)_{n}^\top u \right) \\
\end{bmatrix}.
\]
\end{definition}

As a side note, a simple transformer can also include multiple \textit{attention heads}, in which case several self-attention layers are computed in parallel -- each called an attention head -- using different learned $Q$, $K$, and $V$ matrices. The outputs of all attention heads are then concatenated and passed through point-wise transformations. Equivalently, a simple transformer with multiple attention heads can be written as one with a single head by arranging each $Q$, $K$, and $V$ in block form. All of our results extend naturally to this setting. For clarity of notation, however, throughout this paper we focus on simple transformers with a single attention head, as defined above.

\subsection{Modeling Power} \label{Section Model}
Simple transformers are already used extensively in personalization \citep{wang2018attention,chen2019behavior,bendada2023track,celikik2022reusable}, but other more-sophisticated transformers with multiple self-attention layers are also used. This raises a key question: to what extent can simple transformers effectively model user preferences?

To address this, we begin by building intuition for how simple transformers operate in the personalization context. A simple transformer can be interpreted as modeling the interaction between a user and a set of recommended items, while also capturing pairwise interactions among the items themselves -- commonly referred to as ``set effects''. To illustrate this, consider a concrete example from recommender systems, where the index set $[n]$ represents a set of $n$ items, and a subset $S \subset [n]$ is selected to be recommended to a user $u\in\mathbb{R}^{d_v}$. The matrix $V \in \mathbb{R}^{n \times d_v}$ contains the value vectors, where each row $v_i^\top$ encodes information relevant to item $i$ when considered in isolation. Indeed, when $S = \{i\}$ is a singleton, the output of the self-attention layer is simply $v_i$, and the final reward obtained by recommending item $i$ to the user $u$ is given by $f_i(v_i^\top u)$.

However, when $|S| > 1$, for each $i\in S$ the self-attention layer transforms each value vector $v_{i}$ into a convex combination of the value vectors $\{v_{j} : j \in S\}$, where the weights of the combination are determined by the similarity between the query vector $q_{i}$ and key vectors $\{k_{j} : j \in S\}$, computed via the softmax of inner products. That is, the weights are given by the vector $\mathrm{softmax}((X_SQ) (X_SK)^\top)_i$. Intuitively, the query vector $q_{i}$ for item $i$ captures how the presence of other items in $S$ influences the transformation of value vector $v_{i}$, while the key vector $k_{j}$ for item $j$ captures how item $j$ contributes to the transformation of other items' value vectors. The resulting vector, $\mathrm{SA}_{Q,K,V}(X_S)_i$, is then projected onto the user vector $u$ and passed through the item-specific reward function $f_{i}$. Therefore, final reward obtained by recommending item $i$ to the user $u$ is given by $f_{i}(\mathrm{SA}_{Q,K,V}(X_S)_i^\top u)$, which is precisely the $i$-th coordinate of $\mathcal{T}_{Q,K,V,f_1,\dots,f_n,u}(X_S)$.

Importantly, this formulation is inherently permutation invariant, as the set structure of $S$ does not impose any order on its elements. Sequence models, which require order-sensitive representations, can be accommodated by enriching the matrices $Q$, $K$, and $V$ with additional {\em positional encodings} that inject information about each item's position in an ordered sequence. These positional encodings act like tags that describe each item's position in the sequence, allowing the model to distinguish between, for example, an item that appears first and one that appears last -- even if the items themselves are otherwise identical. In this way, sequential structure is modeled within the same framework as set-based interactions.

Having discussed how simple transformers operate in the personalization setting, we now examine their ability to model user preferences. While our later experiments will empirically demonstrate that restricting the architecture to a single attention layer results in only a modest reduction in modeling/predictive power, we begin by supporting this claim through an analysis of the two common
parametric models used in personalization presented in the previous section. We have already seen that these two models cannot be represented by pure embedding models. On the contrary, below we show that both of them can be represented by simple transformers.

\begin{proposition} \label{prop: variety and subsitution}
 The sequential variety effects in Model \ref{mod:variety} and  the complementarity and substitution effects in Model \ref{mod:com/sub} can both be represented by a simple transformer.
\end{proposition}

The proof of Proposition \ref{prop: variety and subsitution} appears in Appendix \ref{Appendix B}.

\paragraph{Aside: Graph interpretation of self-attention layers.}
It is useful to view a self-attention layer through a graph lens. Let
\[
W = \mathrm{softmax}\big((XQ)(XK)^\top\big)\in\mathbb{R}^{n\times n}.
\]
Because each row of \(W\) sums up to one, \(W\) can be interpreted as the weighted adjacency matrix of a directed graph on \([n]\), with edge weight \(i\to j\) equal to \(W_{ij}\). Then a single self-attention layer performs one round of information passing:
\[
\mathrm{SA}_{Q,K,V}(X)_i = \sum_{j=1}^n W_{ij}\,(XV)_j,
\]
i.e., vertex \(i\) aggregates value vectors from its (out-)neighbors according to the edge weights in \(W\).

This perspective also clarifies what stacking self-attention layers does -- and does \emph{not} do. For \(\ell=1,2\), define
\[
W^{(\ell)} = \mathrm{softmax}\big((XQ^{(\ell)})(XK^{(\ell)})^\top\big),\qquad
H^{(1)} = W^{(1)}(XV^{(1)}),\qquad
H^{(2)} = W^{(2)}\big(H^{(1)}V^{(2)}\big),
\]
where we ignore point-wise activations functions for intuition. Expanding \(H^{(2)}\) shows that information propagate along \emph{two-hop} paths:
\[
H^{(2)}_i
=
\sum_{j=1}^n W^{(2)}_{ij}\,H^{(1)}_j
=
\sum_{j=1}^n\sum_{k=1}^n W^{(2)}_{ij}\,W^{(1)}_{jk}\,\big(XV^{(1)}V^{(2)}\big)_k,
\]
so the contribution of vertex \(k\) to vertex \(i\) through vertex \(j\) factorizes as \(W^{(2)}_{ij}\,W^{(1)}_{jk}\). This is a composition of two pairwise interactions, not an arbitrary \emph{triplet} interaction. Indeed, representing a general triplet interaction requires \(\Theta(n^3)\) free parameters, whereas two self-attention layers only provide \(2n^2\) free parameters. Therefore, stacked self-attention layers naturally model multi-hop effects rather than unrestricted set interactions.

\subsection{Optimization}
The primary purpose of this paper is to study the problem of personalizing a set (or sequence, equivalently) of items in real-time, where the underlying model of user preferences is given by a simple transformer. Following the setup and terminology in the previous subsections, let $[n]$ index a set of items, for which the query, key, and value matrices $Q, K \in \mathbb{R}^{n \times d_{kq}}$ and $V \in \mathbb{R}^{n \times d_v}$ are fixed in advance (these should be thought of as having been learned from prior data). Each user is represented by a vector $u \in \mathbb{R}^{d_v}$ in the same space as the value vectors. For a given set $S \subset [n]$, the simple transformer's output that corresponds to item $i \in S$ is given by $\mathcal{T}_{Q,K,V,f_1,\dots,f_n,u}(X_S)_i=f_{i}( \mathrm{SA}_{Q,K,V}(X_S)_{i}^\top u)$. Intuitively, this can be thought of as the 
reward obtained from recommending item $i$.
As a user arrives, our goal is to choose a set $S \subset [n]$ of at most $k$ items that maximizes the total reward. Formally, the optimization problem we study is:
\begin{definition}[Simple-Transformer Based Recommendations]
\begin{align}
    \label{eqn:problem}
\tag{$\mathsf{Main}$}     \max \quad & \sum_{i\in S} f_i(\mathrm{SA}_{Q,K,V}(X_S)_i^\top u) \\
    \text{s.t.} \quad & S \subset [n],\ |S| \leq k. \notag
\end{align}
    
\end{definition}


Let $\text{OPT}$ denote the optimal objective value of Problem \eqref{eqn:problem}. Notice that if $S=\emptyset$ -- that is, if nothing is recommended to the user --  then the objective value of of Problem \eqref{eqn:problem} equals to 0. Therefore $\text{OPT}\geq 0$. To ensure that Problem \eqref{eqn:problem} is meaningful, from now on we assume that $\textup{OPT} > 0$, since otherwise the best decision would be recommending nothing.



    


\subsection{Hardness}\label{Section Hardness}

As a starting point toward solving Problem~\eqref{eqn:problem}, observe that it can be solved exactly in $O(n^k k^2)$ time via brute-force evaluation of all feasible solutions. As mentioned earlier, we are motivated by settings in which the number of items $n$ is extremely large (possibly hundreds of millions), and Problem~\eqref{eqn:problem} must be solved in real time (possibly milliseconds). Thus, our goal will be to find an algorithm, potentially approximate rather than exact, whose runtime is {\em sub-linear} in $n$, i.e., $O(n^\gamma)$ for some $\gamma < 1$. Moreover, while the budget on the number of items to recommend, $k$, is typically moderate in practice (often around ten), exponential dependence on $k$ would still be impractical to be implemented in real time. Therefore, our algorithm's runtime should also be polynomial in $k$.

Before proceeding, it is useful to present some initial hardness results to temper our expectations. We provide two sets of results. The first addresses the requirement of sub-linear runtime in $n$, with hardness parametrized by $d_{kq}$, the dimension of the key and query vectors. The second addresses the requirement of polynomial runtime in $k$, with hardness parametrized by the \emph{non-negative rank} of the matrix $\mathrm{softmax}(Q K^\top)$, denoted as $\textup{rank}_+(W)$.

\paragraph{Hardness Parametrized by $d_{kq}$.} We first present a proposition that reduces Problem~\eqref{eqn:problem} to graph problems involving cliques, and then discuss its implications.

\begin{proposition}\label{prop: lower bound}
\text{ }
\begin{enumerate}
    \item[\textup{(a)}] If $d_{kq} = n$ and $k\geq 4$, then Problem \eqref{eqn:problem} subsumes the $(k-1)$-\textsc{Clique} problem \footnote{For a (undirected, unweighted) graph, the $k$-\textsc{Clique} problem requires deciding if a clique of size $k$ exists, and finding one if so.} on graphs with $n-1$ vertices.
    \item[\textup{(b)}] For any constant $M\geq 3$, there exists a number $c(M)>0$ for which the following holds. For any $d_{kq}$ such that $\exp(c(M)\cdot d_{kq})\leq n-1$ and any $k\geq M+1$, Problem \eqref{eqn:problem} subsumes the problem of finding a largest clique in a graph with
    \begin{itemize}
        \item $n-1$ vertices,
        \item $\exp(c(M) \cdot d_{kq})$ disjoint cliques,
        \item and all cliques have size at least $k-M$ and at most $k-1$.
    \end{itemize} 
\end{enumerate}
\end{proposition}


The proof of Proposition \ref{prop: lower bound} appears in Appendix \ref{Appendix C}. Proposition \ref{prop: lower bound} implies concrete limitations on the theoretical results we can expect for solving Problem \eqref{eqn:problem}:
\begin{itemize}
    \item By Proposition \ref{prop: lower bound} (a), Problem \eqref{eqn:problem} inherits the hardness of the $k$-\textsc{Clique} problem, which is known to be NP-hard (when $k$ is allowed to grow with $n$) \citep{karp2010reducibility}. Thus, we should not expect to find an exact algorithm which runs in $O(n^C)$ for some $C>0$ independent of $k$.
\item Even treating $k$ as a constant, it is known \citep{chen2006strong} that even an $O(n^{o(k)})$ exact algorithm cannot exist assuming \textit{Exponential Time Hypothesis} holds.\footnote{The Exponential Time Hypothesis (ETH) asserts that the 3\textsc{Sat} problem cannot be solved in sub-exponential time. For a Boolean formula in conjunctive normal form with exactly three literals per clause, the 3\textsc{Sat} problem requires deciding if there exists a truth assignment to the variables that satisfies all clauses, and finding one if so.} Thus, absent additional assumptions, we can only expect to {\em approximately} solve Problem \eqref{eqn:problem} in sub-linear time with respect to $n$.
    
    \item One natural assumption to make is that $d_{kq}$ is small (this is typically the case in practice), and indeed our main result will be parameterized by $d_{kq}$ and only non-trivial when $d_{kq} = o(\log n)$.
     By Proposition \ref{prop: lower bound} (b), if $d_{kq}=\Omega(\log n)$, Problem \eqref{eqn:problem} is at least as hard as finding the largest clique in a graph with $n$ vertices and $\Omega(n)$ disjoint cliques. In particular, each clique in such a graph is itself a candidate maximum clique, and any algorithm must effectively search over $\Omega(n)$ disjoint candidates to determine the largest, since there is no structural overlap between the cliques that an algorithm could exploit to narrow the search space. This indicates that an exact algorithm in sub-linear time with respect to $n$ cannot exist when $d_{kq}=\Omega(\log n)$.
\end{itemize}

\paragraph{Hardness Parametrized by $\textup{rank}_+(W)$.}
Following on the above discussion, we require some additional assumptions to ensure that Problem \eqref{eqn:problem} can be solved, even approximately, in sub-linear time with respect to $n$. One such ``assumption'' will be that $d_{kq}$ is small -- we do not state this as a formal assumption, but rather our main result will be parameterized by $d_{kq}$. 

Similarly, we also require some additional assumptions to ensure that Problem \eqref{eqn:problem} can be solved, even approximately, in polynomial time with respect to $k$. Our main result will be parameterized by a rank-type quantity pertaining to the matrix $W=\textup{softmax}(QK^T)$. Now $QK^\top$ is by definition of rank $d_{kq}$, and while the softmax operator does not preserve rank exactly, it is known that $W$ can be well-approximated by a matrix with rank polynomial in $d_{kq}$ (see \cite{han2023hyperattention,alman2024fast}). Thus, for example, if $d_{kq}$ is constant, then $W$ is approximately low-rank. 

Due to the softmax operator, the matrix $W$ is entry-wise non-negative. This allows us to consider a structural parameter known as the \textit{non-negative rank} of $W$, denoted $\mathrm{rank}_+(W)$. The non-negative rank of a matrix $W \in \mathbb{R}_{\ge 0}^{n \times n}$ is defined as the smallest integer $\mathrm{rank}_+(W)$ such that $W$ can be written as a product of two non-negative matrices: $$W=AB^\top, \text{ where } A,B\in\mathbb{R}_{\ge 0}^{n\times \mathrm{rank}_+(W)}.$$ Equivalently, $W$ can be expressed as the sum of $\mathrm{rank}_+(W)$ non-negative rank-one matrices. Such a representation is called a \emph{non-negative factorization} of $W$. Clearly, this notion is stronger than standard matrix rank; in particular, we always have $\mathrm{rank}(W)\leq \mathrm{rank}_+(W)\leq n$. For more properties on non-negative rank, see \cite{cohen1993nonnegative}. Non-negative rank has many applications in various fields, including data mining, combinatorial optimization, quantum mechanics, and more. In our case, it turns out that the non-negative rank of $W$ is a key parameter that quantifies the hardness of Problem \eqref{eqn:problem}. Formally, we present the following proposition:

\begin{proposition}\label{prop: hardness of k}
    If $\textup{rank}_+(W) = 2$, Problem~\eqref{eqn:problem} admits no $(1-\epsilon)$-approximation scheme
    \begin{itemize}
        \item with runtime $f(1/\epsilon)\,k^{O(1)}$ for any function $f$, assuming the $k$-\textsc{Clique} problem is not \emph{Fixed-Parameter Tractable} \footnote{A problem is said to be Fixed-Parameter Tractable (FPT) if it can be solved in time $f(k)n^{O(1)}$ for some function $f$, where $k$ is a chosen problem parameter. A widely held conjecture is that the $k$-\textsc{Clique} problem is not FPT.} (Corollary of Theorem 6 in \cite{kulik2010there}), and
        \item with runtime $f(1/\epsilon)\,k^{o(1/\epsilon)}$ for any function $f$, assuming Exponential Time Hypothesis holds (Corollary of Theorem 5.1 in \cite{jansen2016bounding}).
    \end{itemize}
    For general $\textup{rank}_+(W)$, Problem~\eqref{eqn:problem} admits no $(1-\epsilon)$-approximation scheme with runtime \[k^{o\left(\frac{\textup{rank}_+(W)}{\epsilon \log^2(\textup{rank}_+(W)/\epsilon)}\right)} \text{ or } k^{o(\sqrt{\textup{rank}_+(W)})},\] assuming \emph{Gap Exponential Time Hypothesis} holds \footnote{The Gap Exponential Time Hypothesis (Gap-ETH) asserts that, for some constant $\epsilon > 0$, distinguishing between satisfiable 3\textsc{Sat} formulas and those that are not even $(1 - \epsilon)$-satisfiable requires exponential time.} (Corollary of \cite{doron2024fine}).
\end{proposition}

The proof of Proposition \ref{prop: hardness of k} appears in Appendix \ref{Appendix C}. Our proof is based on a reduction from Problem \eqref{eqn:problem} to the well-known \textit{Multi-dimensional Knapsack Problem} ($\mathsf{MDKP}$), formally defined in Appendix \ref{Appendix C}. Proposition \ref{prop: hardness of k} shows that, when viewing $k$ as the parameter in Problem~\eqref{eqn:problem}, any algorithm achieving a \((1 - \epsilon)\)-approximation must incur runtime with exponential dependence on \( k \) that necessarily involves both $\textup{rank}_+(W)$ and \( 1/\epsilon \) in a non-trivial way.

\section{Main Result} \label{Section: Main Result}

From the discussions in the previous section, recall that while our algorithm's runtime is parametrized by the non-negative rank of $W$, achieving sub-linear dependence on $n$ and polynomial dependence on $k$ requires that the non-negative rank of $W$ be small. Importantly, our main result does not require $W$ itself to have low non-negative rank, but only that $W$ can be well approximated entry-wise by a low non-negative rank matrix. 

Formally, suppose there exists a non-negative matrix $W' \in \mathbb{R}_{\geq 0}^{n \times n}$ such that
\[
1 - \gamma \leq \frac{W_{ij}}{W'_{ij}} \leq 1 + \gamma \quad \text{for all } i,j,
\]
with $\mathrm{rank}_+(W') = r_+$.\footnote{Actually, we only require a particular sub-matrix to be of non-negative rank $r_+$. This sub-matrix is of significantly smaller size, corresponding to the entries which survives a certain pruning procedure.} Then our main result is parametrized by $r_+$ and $\gamma$. Just as in the case of the parameter $d_{kq}$, our guarantee is non-trivial when $r_+ = O(1)$.

Before stating our main result, we introduce some notation. For a vector $x$, let $x_{\max}$ and $x_{\min}$ denote the largest and smallest entry of $x$, respectively. Likewise, for a matrix $X$, let $X_{\max}$ and $X_{\min}$ denote its largest and smallest entry, respectively. If $X$ has rows $x_i^\top$, we define
\[
\|X\|_{2,\infty} \;=\; \max_i \|x_i\|_2,
\]
that is, the matrix norm induced by the vector $2$- and $\infty$-norms.

We are now prepared to state our main result, which is that our algorithm (to be described in the next section) achieves the following:
\setcounter{theorem}{0}

\begin{theorem} \label{thm: main theorem}
Let $k\textup{-ANN}(n,d,k,\epsilon)$ be the expected amortized runtime of an $\epsilon$-Approximate $k$-Nearest Neighbor algorithm, which we assume is concave in $n$ (e.g.~sub-linear suffices).

Let $\tau$ 
be the number of distinct functions among $f_1,\dots,f_n$. Suppose there exists $W'\in\mathbb{R}^{n\times n}_+$ such that $1-\gamma\leq W_{ij}/W'_{ij}\leq 1+\gamma$ for all $i,j$, and $W'$ has non-negative rank $r_+$ with an explicit non-negative factorization.  
    
    Given any $\epsilon>0$, our algorithm $\textup{ALG}$ achieves \[\textup{ALG}\geq (1-g(\epsilon)-2g(\gamma)-2g((1+\gamma)\epsilon)) \cdot \textup{OPT}-(\epsilon+h(\epsilon))k.\] Moreover, the expected amortized runtime of our algorithm is 
    \[O\left({\epsilon}^{-2d_{kq}} \cdot\tau \cdot k\textup{-ANN}\left(\frac{n
  }{\tau},d_v,k, {\epsilon}\right)+\epsilon^{-2d_{kq}r_+^2/\epsilon}(1+\gamma)^{r_+}(\tau k)^{r_+^2/\epsilon}\right),\]
  where the Big-Oh depends only on $\lVert Q\rVert_{2,\infty},\lVert K\rVert_{2,\infty},(Vu)_{\max},W_{\min},$ and $\max_{i\in[n]}\{f_i((Vu)_{\max})\}$.
\end{theorem}
A few remarks are in order:
\begin{itemize}
    \item Many $\epsilon$-Approximate $k$-Nearest Neighbor algorithms have $k\textup{-ANN}\left(n,d,k, \epsilon\right)$ sub-linear and concave in $n$, and we have given such an example in Lemma \ref{lemma: nearest neighbor}. 
    
    \item In practice, the number of distinct functions $\tau$ is typically very small -- often just one.

    \item In ALG's performance guarantee, the multiplicative dependence on $g$ and the additive dependence on $h$ arise from the parametrization of $f_i$:
\[
f_i(x - \epsilon) \;\geq\; (1 - g(\epsilon)) f_i(x) - h(\epsilon) 
\quad \text{for all } x.
\]
    The additional $\epsilon k$ additive term comes from the use of the 
$\epsilon$-Approximate $k$-Nearest Neighbor algorithm.

    \item By Observation \ref{obs: g and h}, if each $f_i$ is $L$-Lipschitz, our algorithm $\textup{ALG}$ achieves \[\textup{ALG}\geq \textup{OPT}-C\epsilon L k\] for some universal constant $C>0$.
    
    \item If $d_{kq}=o(\log n)$ and $r_+=O(1)$, then the amortized runtime of our algorithm can be simplified to \[\tilde{O}\left(kn^{1-c\epsilon/k}+k^{C/\epsilon}\right).\]
    
    \item Similar to how a small $d_{kq}$ implies low (standard) rank approximation, it also implies that $W$ can be approximated entry-wise by a matrix of low non-negative rank. In particular, we show that if the rows of $Q$ and $K$ can be grouped into $\ell$ clusters, then such a $W'$ can be constructed with $r_+ = \ell(\ell+1)/2$ and a small $\gamma$.
\begin{proposition} \label{prop: cluster}
Let $W = \textup{softmax}(QK^\top)$. Suppose $I_1, \dots, I_\ell$ form a partition of $[n]$ such that for every $\ell' \in [\ell]$ and $i, i' \in I_{\ell'}$, we have $\lVert q_i - q_{i'} \rVert_2 \leq \delta$ and $\lVert k_i - k_{i'} \rVert_2 \leq \delta$. Then we can construct $W' \in \mathbb{R}^{n \times n}_{\geq 0}$ such that
\begin{enumerate}
    \item $1-\gamma\leq W_{ij}/W'_{ij}\leq 1+\gamma$ for all $i,j$, where $\gamma = 17\delta \max\{ \lVert Q \rVert_{2,\infty}, \lVert K \rVert_{2,\infty} \}$,
    \item and $W'$ has non-negative rank $r_+ = \ell(\ell+1)/2$ with an explicit non-negative factorization.
\end{enumerate}
\end{proposition}
To obtain the clusters, one approach is to directly partition the row spaces of $Q$ and $K$, and group the rows of $Q$ and $K$ according to this partition. Since $d_{kq}$ is assumed to be small, the number of clusters is also small. This yields the following corollary.
\begin{corollary} \label{cor: rank partition}
    Let $W = \textup{softmax}(QK^\top)$. Given any $\delta>0$, we can construct $W'\in\mathbb{R}^{n\times n}$ such that $1-\gamma\leq W_{ij}/W'_{ij}\leq 1+\gamma$ for all $i,j$ where $\gamma=17\delta \max\{\lVert Q\rVert_{2,\infty},\lVert K\rVert_{2,\infty}\}$, and $W'$ has non-negative rank $$r_+=\lceil4\max\{\lVert Q\rVert_{2,\infty},\lVert K\rVert_{2,\infty}\}/\delta\rceil^{4d_{kq}}$$ with an explicit non-negative factorization.
\end{corollary}
The proofs of Proposition \ref{prop: cluster} and Corollary \ref{cor: rank partition} appear in Appendix \ref{Appendix D}. Notice that if $d_{kq}=O(1)$, Corollary \ref{cor: rank partition} shows that $r_+=O(1)$ for some small $\gamma$. In this case our main result is non-trivial. 

    \item Computing a non-negative matrix factorization is in general NP-hard \citep{vavasis2010complexity}: there is a rich literature on this subject that has yielded multiple algorithms (see \cite{lee2000algorithms,wang2012nonnegative} for surveys). We will not be concerned with this runtime in analyzing Problem \eqref{eqn:problem}, as $Q$ and $K$ are given beforehand, and thus we view this as amortized across multiple instances of Problem \eqref{eqn:problem}.

    \item In practice, the embedding dimension $d_{kq}$ of items is usually much smaller than the number of items $n$, so $d_{kq}=o(\log n)$ often holds. In addition under the condition that $r_+=O(1)$, for any fixed $\epsilon>0$, our algorithm gives a $1-\epsilon$ approximation algorithm with expected amortized runtime that is sub-linear in the total number of items $n$, and polynomial in the number of items to recommend $k$.
    
    
    
    \item Our algorithm operates under the same two-phase {\em retrieve} and {\em rank} paradigm that is used in many competition-winning personalization algorithms, where both phases are adapted specifically to simple transformers. We will discuss these two phases in detail later in this section. The two pieces in our algorithm's expected amortized runtime directly correspond to the expected amortized runtime of our algorithm's two phases: the expected amortized runtime of phase one is $$O\left({\epsilon}^{-2d_{kq}} \cdot\tau \cdot k\textup{-ANN}\left(\frac{n
  }{\tau},d_v,k, {\epsilon}\right)\right),$$ and the runtime of phase two is $$O\left(\epsilon^{-2d_{kq}r_+^2/\epsilon}(1+\gamma)^{r_+}(\tau k)^{r_+^2/\epsilon}\right).$$ 
    
\end{itemize}

Real-time personalization algorithms typically operate in two phases: \textit{retrieval} and \textit{ranking}. In phase one (retrieval), the algorithm efficiently selects a small subset of promising candidate items from a large item pool. This is typically achieved using ANN algorithms over learned item embeddings. In phase two (ranking), an optimization problem, which is determined by the specific personalization model, is approximately solved over the retrieved subset to produce the final recommendations. 

Our algorithm also follows the retrieval and ranking structure, but introduces novel methods in both the retrieval and ranking phases. At a high level, it operates as follows:

\paragraph{Phase One (Retrieval).}  
In the retrieval phase, our algorithm leverages ANN algorithms while explicitly incorporating the structure of the transformer architecture. Specifically, we begin by partitioning items offline based on their query vectors $q_i$, key vectors $k_i$, and reward functions $f_i$. This partitioning is designed so that items within the same partition produce similar outputs when processed through the self-attention layer, given a fixed user query. When a user $u$ arrives, the algorithm applies a given $\epsilon$-Approximate $k$-Nearest Neighbor algorithm (as defined in Section~\ref{Section: before transformers}) within each partition to retrieve up to $k$ items that are likely to be most attractive to the user. The final set of candidate items is obtained by taking the union of the retrieved items across all partitions.

\paragraph{Phase Two (Ranking).}  
In the ranking phase, our algorithm generalizes and improves upon the widely used \textit{Beam Search} heuristic. Beam Search explores a decision tree greedily, where each level corresponds to selecting the next item to add to the recommendation set $S$. Concretely, given a candidate budget, Beam Search generates a prescribed number of candidate solutions and returns the one with the highest reward. Each candidate solution is represented by a $k$-tuple $(b_1,\dots,b_k)$. To determine the $\ell$-th item in the tuple, the algorithm evaluates each item not yet selected by computing the incremental gain in the objective if the item were added, then chooses the item corresponding to the $b_\ell$-th largest increment. In particular, the tuple $(1,\dots,1)$ corresponds to the fully greedy solution. The values $b_1,\dots,b_k$ are tuned to control the number and diversity of candidate solutions. 

Our algorithm adopts a similar tree-exploration framework but introduces a refinement that yields provable performance guarantees. Specifically, it explores the decision tree only up to a limited number of levels, after which it applies linear programming (LP) rounding techniques to optimize over the remaining items. This creates a natural trade-off: exploring too many levels incurs excessive runtime, while exploring too few levels increases the error introduced by LP rounding. By carefully balancing this trade-off, our algorithm achieves both practical runtime -- by limiting the depth of exploration -- and near-optimal reward -- by minimizing the rounding loss over the unexplored levels.

\section{Algorithm and Proof of Main Theorem}

Before diving into the details of our algorithm, we rewrite Problem~\eqref{eqn:problem} with slightly different notations. Let $q_i^\top$ and $k_i^\top$ denote the $i$-th rows of $Q$ and $K$, respectively. Let $W = \textup{softmax}(QK^\top)$, and let $w_i^\top$ be the $i$-th row of $W$. For vectors $a,b\in\mathbb{R}^n$, let $a\odot b\in\mathbb{R}^n$ denote the element-wise product of $a$ and $b$, i.e., $(a\odot b)_i=a_ib_i$ for every $i\in[n]$. We make the following observation.

\begin{observation}\label{obs: simplification}
    Problem~\eqref{eqn:problem} is equivalent to the following Problem~\textup{($\mathsf{P}$)}:
    \begin{align}
       \tag{$\mathsf{P}$} \label{eqn: simplification}
        \max \quad & f_{\mathsf{P}}(x) = \sum_{i=1}^n x_i f_i\left(\frac{(w_i \odot Vu)^\top x}{w_i^\top x}\right) \\
        \textup{s.t.} \quad & x \in \{0,1\}^n,\quad 1 \leq e^\top x \leq k \notag
    \end{align}
    where $e \in \mathbb{R}^n$ is the all-ones vector, and $x_i = 1$ indicates that item~$i$ is selected into the set $S$.
\end{observation}

The proof of Observation~\ref{obs: simplification} appears in Appendix~\ref{Appendix E}. From this point onward, we will work with $\mathsf{P}$, using its notation in place of Problem~\eqref{eqn:problem}.

Below we give the pseudo-code of our main algorithm.

\begin{algorithm}[H]
\caption{Main Algorithm}
\label{alg: complete}
\DontPrintSemicolon
\SetAlgoLined

\KwIn{Number of items $n$, maximum number of recommended items $k$, key matrix $K\in\mathbb{R}^{n\times d_{kq}}$, query matrix $Q\in\mathbb{R}^{n\times d_{kq}}$, value matrix $V\in\mathbb{R}^{n\times d_v}$, reward functions $\{f_i\}_{i=1}^n$, user vector $u\in\mathbb{R}^{d_v}$, parameter $\epsilon>0$, $\epsilon$-Approximate $k$-Nearest Neighbor oracle (ANN), attention matrix $W=\textup{softmax}(QK^\top)\in\mathbb{R}_+^{n\times n}$, low-rank approximation $W'\in\mathbb{R}_+^{n\times n}$ with factorization $W'=AB^\top$ and element wise guarantee $1-\gamma \le W_{ij}/W'_{ij} \le 1+\gamma$}
\KwOut{Solution $x$ to Problem \eqref{eqn:problem}.}

\BlankLine
\tcp{Phase One: Preprocess (Algorithm \ref{alg: phase one preprocess})}
$(\delta,\mathcal{I},\mathcal{S},\text{preprocessed ANN}) \leftarrow$ \textbf{Run} Algorithm \ref{alg: phase one preprocess} with inputs $(n,k,K,Q,V,\{f_i\},\epsilon,\text{ANN})$\;

\BlankLine
\tcp{Phase One: Query (Algorithm \ref{alg: phase one query})}
$I \leftarrow$ \textbf{Run} Algorithm \ref{alg: phase one query} with inputs $(u,k,\delta,\mathcal{I},\mathcal{S},\text{preprocessed ANN})$\;

\BlankLine
\tcp{Phase Two (Algorithm \ref{alg: phase two})}
$x \leftarrow$ \textbf{Run} Algorithm \ref{alg: phase two} with inputs $(W,W'=AB^\top,\gamma,V,u,k,I,\epsilon)$\;

\BlankLine
\Return{$x$}
\end{algorithm}

\subsection{Phase One (Retrieval)} \label{Section Phase One}
In phase one, our algorithm aims to identify a small subset $I \subset [n]$ of items such that the optimal objective value of  \textup{$\mathsf{P}$} does not decrease significantly when restricted to $I$. To achieve this, our algorithm first partitions items offline based on their query vectors $q_i$, key vectors $k_i$, and reward functions $f_i$. Since both $q_i$ and $k_i$ lie in a space of dimension $d_{kq}$, the number of partitions is much smaller than $n$. Items in the same partition are designed to behave similarly under the self-attention layer. That is, they produce similar outputs when interacting with other items. When a user $u$ arrives, the algorithm applies an $\epsilon$-Approximate $k$-Nearest Neighbor algorithm (as defined in Section~\ref{Section: before transformers}) within each partition to select at most $k$ items whose value vectors have the highest approximate inner product similarity with $u$. Because items within the same partition respond similarly under attention, user preferences within each partition are primarily determined by similarity to the user vector. Thus, our algorithm retrieves a small subset $I$ containing high-reward items tailored to the user.

\begin{proposition}[Phase One]\label{prop: phase one}
       Suppose we have an $\epsilon$-Approximate $k$-Nearest Neighbor algorithm with expected amortized runtime $k\textup{-ANN}(n,d,k,\epsilon)$. Let $\tau
       $ be the number of distinct functions among $f_1,\dots,f_n$. Given any $\epsilon>0$, Algorithm \ref{alg: phase one query} returns an index set $I\subset [n]$ such that \begin{equation*}
|I|=\tau k\left\lceil\frac{140\max\{\lVert Q\rVert_{2,\infty},\lVert K\rVert_{2,\infty}\})^2 }{(Vu)_{\max}\cdot\epsilon}\right\rceil^{2d_{kq}},
\end{equation*} and the optimal value to the following Problem $\mathsf{P}(I)$ 
\begin{align}
    \tag{$\mathsf{P}(I)$} \label{eqn:probPI}
    \max \quad & f_{\mathsf{P}(I)}(x)=\sum_{i=1}^n x_i f_i\left(\frac{(w_i \odot Vu)^\top x}{w_i^\top x}\right) \\
    \textup{s.t.} \quad & x \in \{0,1\}^n,\ x_i = 0 \text{ for } i \notin I,\ 1 \leq e^\top x \leq k \notag
\end{align}
satisfies \begin{equation*}
            \textup{OPT}_{\mathsf{P}(I)}\geq (1-g(\epsilon))\textup{OPT}_{\mathsf{P}}-kh(\epsilon).
        \end{equation*} 
        Moreover, suppose $k\textup{-ANN}(n,d,k, \epsilon)$ is concave in $n$. Then the expected amortized runtime of Algorithm \ref{alg: phase one query} is \begin{equation*}
  \left\lceil\frac{140(\max\{\lVert Q\rVert_{2,\infty},\lVert K\rVert_{2,\infty}\})^2(Vu)_{\max}}{\epsilon}\right\rceil^{2d_{kq}} \tau \cdot k\textup{-ANN}\left(\frac{n
  }{\tau},d_v,k, \frac{\epsilon}{35\max\{\lVert Q\rVert_{2,\infty},\lVert K\rVert_{2,\infty}\}(Vu)_{\max}}\right).
\end{equation*}
\end{proposition}

We make the following remark: suppose $d_{kq}=o(\log n)$ and $\lVert Q\rVert_{2,\infty},\lVert K\rVert_{2,\infty}, (Vu)_{\max}$ are all viewed as constants, then for any given constant $\epsilon>0$, we have $|I|=\tau k n^{o(1)}$. Moreover, suppose $k\textup{-ANN}\left(n,d_v,k, \epsilon\right)$ is sub-linear in $n$ when $d_v,k,\epsilon$ are fixed, then  the expected amortized runtime of Algorithm \ref{alg: phase one query} is also sub-linear in $n$.

\subsection{Phase Two (Ranking)} \label{Section Phase Two}
In phase two, our algorithm approximately solves $\mathsf{P}(I)$, which is $\mathsf{P}$ over the retrieved subset of items $I\subset [n]$. Without loss of generality, assume $I=[m]$. By the first remark of Proposition \ref{prop: phase one}, we may treat $m=kn^{o(1)}$ under mild assumptions. Then since $x_i=0$ for $i>m$, we may consider only the first $m$ entries of $Vu$ and the top-left $m\times m$ principal sub-matrix of $W$. Therefore, with slight abuse of notation, we redefine $Vu\in\mathbb{R}^m$ to include its first $m$ entries, and $W,W'\in\mathbb{R}^{m\times m}_{+}$ to be the top-left $m\times m$ principal sub-matrices of the corresponding matrices, respectively. Moreover, note that the quantity $$\frac{(w_i\odot Vu)^\top x}{ w_i^\top x}$$ remains unchanged if the vector $w_i$ is multiplied by a non-zero constant. Thus, rescaling the rows of $W$ does not change $\mathsf{P}(I)$. Because $W\in\mathbb{R}^{m\times m}_{+}$ is now the $m\times m$ principal sub-matrix, the sum of its rows is not normalized to 1. So for simplicity of exposition, we assume each row of $W$ is rescaled so that $\sum_{j=1}^m W_{ij}=1$, and each row of $W'$ is rescaled accordingly so that $1-\gamma\leq W_{ij}/W'_{ij}\leq1+\gamma$ for all $i,j$. Then we may rewrite  $\mathsf{P}(I)$ as \begin{align}
    \tag{$\mathsf{P}(I)$} \label{eqn:problem2}
    \max \quad & f_{\mathsf{P}(I)}(x)=\sum_{i=1}^m x_i f_i\left(\frac{(w_i \odot Vu)^\top x}{w_i^\top x}\right) \\
    \textup{s.t.} \quad & x \in \{0,1\}^m,\ 1 \leq e^\top x \leq k. \notag
\end{align}
From this point onward, we will work with this new form of $\mathsf{P}(I)$.

Our algorithm begins by replacing $W$ with a low non-negative rank surrogate $W'$ and showing that solving the problem under this approximation is sufficient. Rather than exhaustively enumerating all possible solutions, our algorithm then focuses on a restricted collection of partial solutions that retain the key structural information. The nonlinear terms in the objective are handled by introducing a family of auxiliary linearized problems, which can be further simplified through discretization. This reduction ensures that only a small number of auxiliary linearized problems need to be solved.  

To address each auxiliary linearized problem, our algorithm employs a rounding procedure that converts fractional linear-programming solutions into valid discrete ones. At this stage, the central trade-off emerges: exploring too many partial solutions increases runtime beyond practical limits, while exploring too few places excessive burden on the rounding step, leading to higher approximation error. By carefully balancing this trade-off, the ranking phase achieves both computational efficiency -- through controlled exploration -- and strong accuracy -- by minimizing the loss introduced during rounding.

\begin{proposition}\label{prop: phase two}
        
Suppose there exists $W'\in\mathbb{R}^{n\times n}_+$ such that $1-\gamma\leq W_{ij}/W'_{ij}\leq 1+\gamma$ for all $i,j$, and $W'$ has non-negative rank $r_+$ with an explicit non-negative factorization. Given any $\epsilon>0$, our algorithm $\textup{ALG}$ achieves \begin{align*}
\textup{ALG}_{\mathsf{P}(I)}
    &\geq (1-g(2\gamma(Vu)_{\max}))^2(1-g(c_{\epsilon,\gamma,W'_{\min}}))^2\textup{OPT}_{\mathsf{P}(I)}\\&\quad\text{ }-k(1-g(2\gamma(Vu)_{\max}))(2\epsilon(1-g(c_{\epsilon,\gamma,W'_{\min}}))+g(c_{\epsilon,\gamma,W'_{\min}})h(c_{\epsilon,\gamma,W'_{\min}})\\&\quad\text{ }+(1-g(c_{\epsilon,\gamma,W'_{\min}}))^2h(2\gamma(Vu)_{\max})+h(2\gamma(Vu)_{\max})),
\end{align*} where $$c_{\epsilon,\gamma,W'_{\min}}=\frac{(1+\gamma)\epsilon}{W'_{\min}},$$ with runtime \begin{align*}
    &\quad \text{ } \text{ } r_+ m \log_2 m +\lambda m^\lambda+ \lambda r_+^2m^{\lambda r_+}\\&+\left\lceil\left(\frac{4(1+\gamma)k}{\epsilon W'_{\min}}\right)^{r_+}\right\rceil \cdot \left\lceil \frac{(Vu)_{\max}-\min\{0,(Vu)_{\min}\}}{\epsilon}\right\rceil^{r_+}\cdot\left\lceil \frac{\max_{i\in[m]}\{f_i((Vu)_{\max})\}}{\epsilon}\right\rceil\cdot\lambda'm^{\lambda r_++\lambda'}T_{\textup{LP}},
\end{align*} where $\lambda=\lceil(2r_++2)(Vu)_{\max}/\epsilon\rceil$ and $\lambda'=\lceil(2r_++2)\max_{i\in[m]}\{f_i((Vu)_{\max})\}/\epsilon\rceil$. Here, $T_{\textup{LP}}=\textup{LP}(m,3m+r_++2)+\textup{LP}(m,2m+2r_++2)$ and $\textup{LP}(m,n)$ is the runtime of solving a linear program with $m$ variables and $n$ constraints.
\end{proposition}

Our proof of Proposition \ref{prop: phase two} appears in Appendix \ref{Appendix F}. The proof is completed according to the following steps:
\begin{enumerate}
    \item {\bf Low Non-negative Rank Approximation:} In Lemma \ref{lemma: replace W'}, we prove that in order to approximately solve $\mathsf{P}(I)$, it is sufficient to approximately solve $\mathsf{P}'(I)$, where $\mathsf{P}'(I)$ is obtained by replacing $W$ with $W'$.
    
    \item {\bf Enumeration of Partial Solutions:} We took guesses on some index sets $X_1,\dots, X_{r_+}$, which corresponds to the non-negative factorization of $W'$. Let $X=(X_1,\dots, X_{r_+})$ and let $\mathsf{P}(X)$ denote the problem where $\mathsf{P}'(I)$ has additional constraints that $x_i = 1$ for all $i\in\cup_j X_j$. We showed that the total number of guesses $X$ is bounded above, so it is sufficient to solve $\mathsf{P}(X)$ for each guess $X$.
    
    \item {\bf Linearization of Fractional Objective Terms:} In order to solve $\mathsf{P}(X)$, we linearize the fractional terms in the objective function by defining a set of auxiliary problems $\mathsf{P}(X,t)$ for each $t\in\mathbb{R}^m_+$. These problems are parameterized by the denominator terms in the objective function of $\mathsf{P}(X)$. 
    In Lemma \ref{lemma: auxiliary problem}, we prove it suffices to find a $t^*$ for which $\mathsf{P}(X,t^*)$  has the highest optimal value among all $\mathsf{P}(X,t)$'s, as the corresponding optimal $x^*$ 
    is an optimal solution to $\mathsf{P}(X)$. 
    
    \item {\bf Dimensionality Reduction and Discretization of Auxiliary Problems:} In order to approximately solve $\mathsf{P}(X,t)$ for all $t\in\mathbb{R}^m_+$, we discretize $t$-space and show in Lemma \ref{lemma: auxiliary oracle is enough} that it suffices to solve $\mathsf{P}(X,t)$ for a small number of $t$'s.
    
    \item {\bf Complete Linearization of Auxiliary Problems:} Fix a given $t$, the objective functions of $\mathsf{P}(X,t)$ inside $f$ has rank $r_+$. We discretized the value space of those objective functions. In Lemma \ref{lemma: auxiliary oracle on each tuple is enough}, we showed that in order to solve $\mathsf{P}(X,t)$, it is sufficient to give an oracle that, for each discretization of the value space, identify whether there exists a feasible solution to $\mathsf{P}(X,t)$ with objective values that are approximately inside the discretization.
    
    \item {\bf Approximation of Linearized Auxiliary Problems via LP Rounding:} Finally, we gave such an oracle by a rounding procedure. Lemma \ref{lemma: constant fractional components} and Lemma \ref{lemma: z satisfies auxiliary oracle} proved that the oracle is correct by using the properties of our guess $X$.
    
\end{enumerate}

\section{Experimental Results}

We performed two sets of experiments, which demonstrate the following.

\paragraph{Representation.} Simple transformers empirically captured user preferences nearly as accurately as more sophisticated transformer-based models. In particular, in a machine learning task involving learning from user behaviors and predicting user preferences, simple transformers achieved an average accuracy that was 14.1\% higher than the best among various non-attention machine learning models, such as logistic regression, random forest, and support vector machines. Compared to more general transformer models  --  that is, transformers with more self-attention layers  --  the accuracy of simple transformers was only 2.5\% lower on average. These results demonstrate that simple transformers can learn from user behaviors and predict user preferences with much higher accuracy than non-attention models and with accuracy nearly matching that of general transformers.

\paragraph{Optimization.} Our algorithm completes simple-transformer-based recommendation tasks both efficiently and accurately. Using the parameters learned from simple transformers in the initial set of experiments, we conducted an optimization task involving the recommendation of a set of items to each arriving user. Recall that our algorithm operates in two phases: retrieval and ranking. We compared each phase of our algorithm against a natural benchmark: $k$-Nearest Neighbors for retrieval and Beam Search for ranking. Combining the two retrieval methods with the two ranking methods yields four total algorithmic variants. We evaluated performance based on the best candidate solution produced from a fixed total number of candidate solutions. On average, our algorithm achieved an objective value that was 20.56\% higher than that of the algorithm using our retrieval method and Beam Search for ranking, and 20.86\% higher than that of the algorithm using $k$-Nearest Neighbors for retrieval and our ranking method, with the same fixed total number of candidate solutions. These results demonstrate that our algorithm outperforms the natural benchmarks in both phases of the recommendation process. 

We used two dataset. The first dataset was the Spotify Million Playlist Dataset \citep{chen2018recsys}. Spotify is one of the largest music streaming platforms, with over 640 million monthly active users, including 252 million paying subscribers. The dataset comprises one million user-generated playlists created on the Spotify platform between January 2010 and October 2017. Each playlist includes features such as the playlist title and the titles of the tracks it contains.

The second dataset was the Trivago 
 Session-based Hotel Recommendations Dataset \citep{Knees_etal:RSC:2019}. Trivago is a global hotel search platform that operates 55 localized websites across more than 190 countries, providing access to over two million hotels. The dataset consists of user sessions related to hotel search and booking, encompassing approximately 730,000 unique users and 340,000 unique hotels across roughly 900,000 sessions. Each session includes information on user interactions with hotels, such as clicks and checkouts. In addition, the dataset contains various hotel attributes, including price, city, and other relevant features.

\subsection{Representation}

In this set of experiments, our goal was to show that simple transformers were able to learn from user behaviors and predict user preferences with high accuracy. More specifically:

\paragraph{Spotify.} In the Spotify experiment, we designated playlists containing 20 songs as “true” playlists. To construct “fake” playlists, we took the first 15 songs from each true playlist and appended 5 randomly selected songs. The number of true and fake playlists was balanced to be equal. Given a playlist, the task was to classify it as either true or fake. The performance of an algorithm was evaluated based on the average classification accuracy.
\paragraph{Trivago.} In the Trivago experiment, each session provided information on a user's interactions with the first 15 hotels. The task was to predict the user's interactions with the subsequent 5 hotels in the same session. The performance of an algorithm was evaluated based on the average prediction accuracy.

We compared three classes of machine-learning algorithms on these prediction tasks:

\begin{itemize}
\item \textbf{Non-Attention Models:} This class included well-known machine learning algorithms that do not incorporate self-attention mechanisms, such as random guessing, logistic regression, support vector machines, and nearest neighbors. These models disregarded any potential sequential structure in the data.
\item \textbf{Simple Transformers:} This class consisted of transformer architectures with a single self-attention layer,  followed by linear layers and activation functions.
\item \textbf{General Transformers:} This class contained more complex transformer architectures with multiple self-attention layers and potentially deeper network structures.
\end{itemize}

Below we give the architecture of the simple transformers used in both experiments. 
\paragraph{Architecture used in Spotify.} Let $x_i$ denote the word2vec embedding of the $i$-th song in a playlist, after being processed by a linear layer. Each user vector is modeled as the average of the embeddings of the first 15 songs in the playlist, that is, $u = \sum_{j=1}^{15} x_j$. Let $\text{SA}_{Q,K,V}(\cdot)$ be the  self-attention layer with learned parameters $Q, K, V$, and let $f_i(\cdot) = f(\cdot)$ be an activation function composed of a linear transformation and a logistic function, both parameterized and learned during training. Let $S$ be the set of indices corresponding to the 16th through 20th songs in the playlist. Then, the output of the simple transformer is given by
\[
\sum_{i=1}^{5} f(\mathrm{SA}_{Q,K,V}(X_S)_i^\top u),
\]
where $X_S$ denotes the input embeddings corresponding to songs in set $S$. General transformer models extend this architecture by incorporating additional self-attention layers. 

\begin{figure}[h!]
    \centering
    \includegraphics[width=0.8\linewidth]{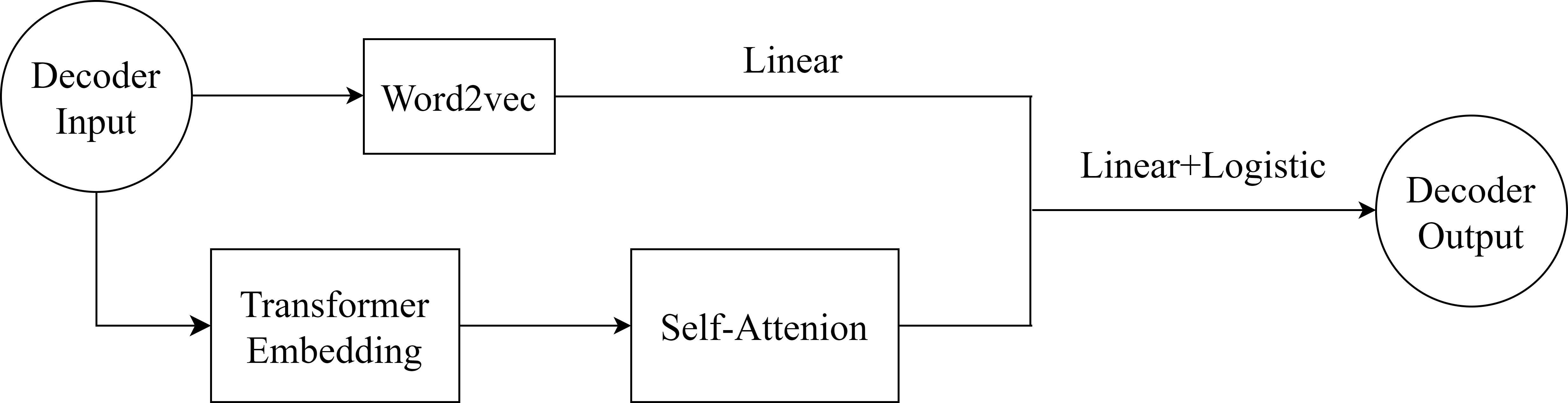}
    \caption{Architecture of the simple transformer used in the Spotify experiment.}
    \label{fig: Sportify arc}
\end{figure}
    
\paragraph{Architecture used in Trivago.}
Let $x_i$ denote the learned embedding of the $i$-th hotel. Each user vector is modeled as the average of the embeddings of the first 15 hotels the user engaged within a given session, that is, $u = \sum_{j=1}^{15} x_j$. Let $\text{SA}_{Q,K,V}(\cdot)$ denote a self-attention layer with learned parameters $Q$, $K$, and $V$, and let $f_i(\cdot) = f(\cdot)$ be an activation function composed of a linear transformation and a logistic function, both parameterized and learned during training. Let $S$ be the set of indices corresponding to the 16th through 20th hotels in the session. Then, the output of the simple transformer is given by
\[
\sum_{i=1}^{5}f(\mathrm{SA}_{Q,K,V}(X_S)_i^\top u),
\]
where $X_S$ denotes the input embeddings corresponding to hotels in set $S$. The simple transformer architecture consisted solely of a decoder with one self-attention layer. General transformer models extended this architecture by incorporating an encoder, as well as additional self-attention layers in both the encoder and decoder.

\begin{figure}[h!]
    \centering
    \includegraphics[width=0.9\linewidth]{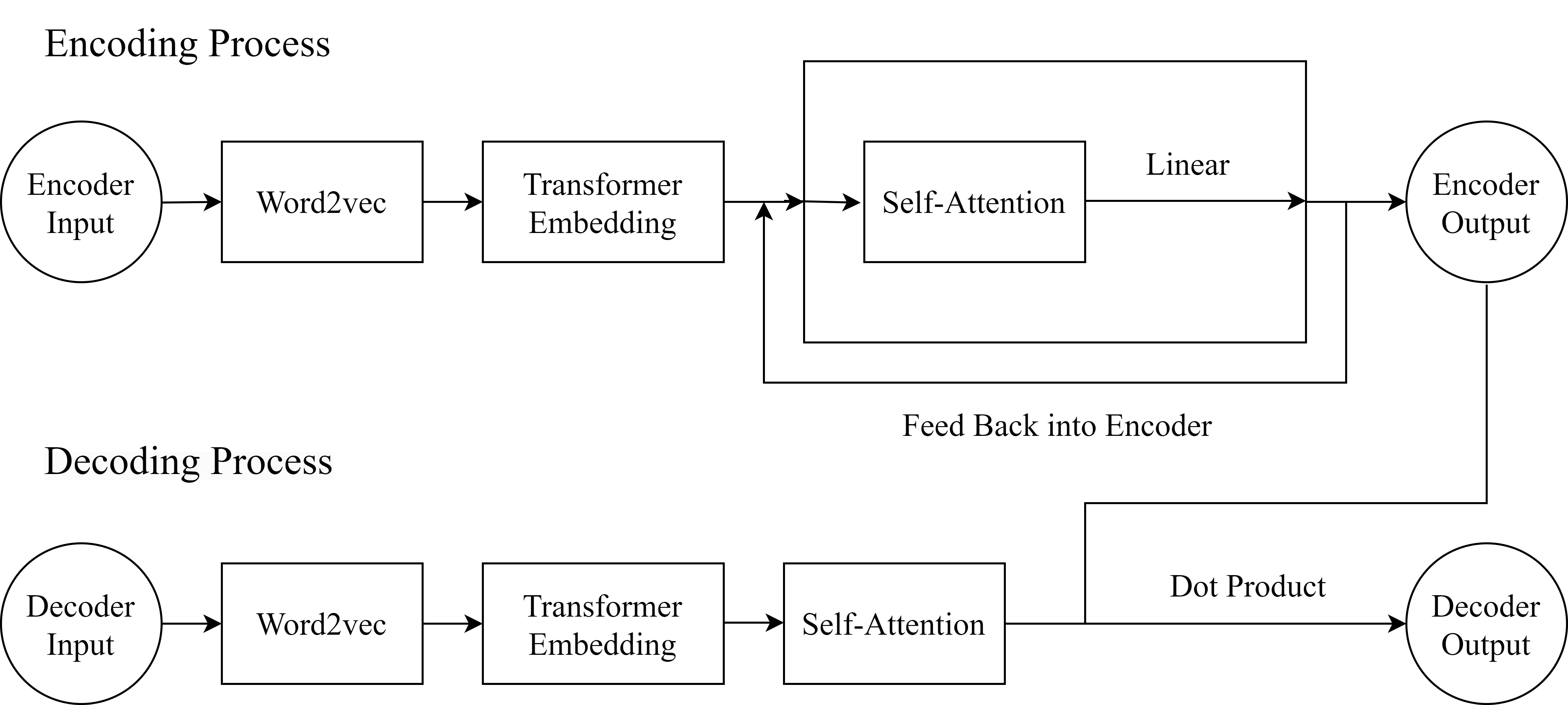}
    \caption{Architecture of the transformer used in the Trivago experiment. The simple transformer only contained the decoder, that is, a single self-attention layer.}
    \label{fig: Trivago arc}
\end{figure}

The experimental results are presented in the tables below. 

\begin{table}[h!]
\centering
\begin{tabular}{@{}lccccc@{}}
\toprule
 & \begin{tabular}[c]{@{}c@{}}Random\\ Forest\end{tabular} 
 & \begin{tabular}[c]{@{}c@{}}Logistic\\ Regression\end{tabular} 
 & \begin{tabular}[c]{@{}c@{}}Support Vector\\ Machine\end{tabular}
 & \begin{tabular}[c]{@{}c@{}}Simple\\ Transformer\end{tabular} 
 & \begin{tabular}[c]{@{}c@{}}General\\ Transformers\end{tabular} \\
\midrule
Spotify & 0.518 & 0.520 & 0.334 & 0.702 & 0.726 \\ 
Trivago & 0.271 & 0.530 & 0.531  & 0.631 & 0.742 \\
\bottomrule
\end{tabular}

\caption{Average accuracy of different machine-learning models on Spotify and Trivago.}
\label{table: spotify-trivago}
\end{table}

\begin{table}[h!]
\centering
\begin{tabular}{@{}cccc@{}}
\toprule
\begin{tabular}[c]{@{}c@{}}Dec. Layers \textbackslash\ Enc. Layers\end{tabular} & 1 & 2 & 4 \\
\midrule
1 & 0.590 & 0.602 & 0.596 \\
2 & 0.654 & 0.692 & 0.700 \\
4 & 0.724 & 0.742 & 0.675 \\
\bottomrule
\end{tabular}
\caption{Average accuracy of general (full encoder–decoder) transformers with various numbers of self-attention layers on Trivago.}
\label{Trivago table 2}
\end{table}

In Table~\ref{table: spotify-trivago}, the simple transformer outperformed the non-attention models by an average accuracy margin of 0.182 on the Spotify dataset and 0.2 on the Trivago dataset, while achieving performance comparable to that of general transformers with additional self-attention layers. Moreover, on the Trivago dataset, compared to the various general transformer architectures in Table~\ref{Trivago table 2}, the simple transformer outperformed some more complex architectures and performed only slightly worse than the best-performing ones. Specifically, its accuracy was only 2.4\% lower than the average accuracy across all general transformers. It is important to note that, in practice, the optimal architecture is not known a priori. Therefore, the simple transformer represents a strong and robust choice for this prediction task. In summary, simple transformers effectively learned from user behaviors and predicted user preferences with substantially higher accuracy than non-attention models, and with accuracy nearly matching that of general transformers.

\subsection{Optimization}

In the previous set of experiments, we demonstrated that simple transformers could empirically capture user preferences on both datasets. In this set of experiments, we turned to the task of personalized recommendation based on simple transformers. We treated the parameters $Q$, $K$, and $V$ learned in the previous experiments as ground truth and solved Problem (\ref{eqn:problem}). Each instance corresponded to an arriving user. More specifically:

\paragraph{Spotify.} In the Spotify experiment, we were given 15 songs as input, and the task was to recommend an additional 5 songs to complete a 20-song playlist. The given 15 songs were treated as a representation of the user. The corresponding user vector $u$ was computed by first averaging the word2vec embeddings of the 15 songs, and then applying a linear transformation to project the result into the value space. The key, query, and value vectors of each song were obtained from the parameters learned in the previous experiments.

\paragraph{Trivago.} In the Trivago experiment, we were given 15 user interactions as input, and the task was to recommend 5 additional hotels to maximize the booking rate. The user vector $u$ was computed by averaging the learned embeddings of the 15 hotels with which the user had interacted. The key, query, and value vectors of each hotel were similarly obtained from the parameters learned in the previous experiments.

Recall that our algorithm operates in two phases: retrieval and ranking. Our algorithm performs these two phases as described below:

\begin{itemize}
    \item \textbf{Phase One (Retrieval):} Our algorithm partitions the row space of the query matrix $Q$ and the key matrix $K$ offline, with the number of partitions treated as a tunable hyperparameter. Then, when a user vector $u$ arrives, for each group of items whose query and key vectors both belong to the same corresponding partition, we apply the $k$-Nearest Neighbor algorithm and retain the $k$ items with the highest base reward $v^\top u$. All retained items are then passed to the ranking phase.

\item \textbf{Phase Two (Ranking):} Our algorithm first enumerates all possible combinations of the top $c$ highest-reward items to include in a candidate solution -- referred to as \emph{valid tuples} in the proof of Proposition~\ref{prop: phase two} -- where $c$ is a tunable hyperparameter. For the remaining items, the algorithm solves the residual subproblem by evaluating a fixed number of auxiliary problems, denoted as $\mathsf{P}(X, t)$ in the same proof. The number of auxiliary problems solved was tuned to control the total number of candidate solutions generated.

\end{itemize}

We compared each phase of our algorithm against a natural benchmark:

\begin{itemize}
    \item \textbf{$k$-Nearest Neighbor (Retrieval):} The $k$-Nearest Neighbor algorithm served as the benchmark for the retrieval phase. It ignored any potential sequential effects and instead ranked items solely based on their base rewards, computed as $v^\top u$. The algorithm then greedily selected the same number of items as our algorithm's retrieval phase to pass on to the ranking phase.

    \item \textbf{Beam Search (Ranking):} Beam Search served as the benchmark for the ranking phase. Beam Search is a greedy-type heuristic commonly used in practice. Each candidate solution generated by Beam Search is specified by a $k$-tuple $(b_1,\dots,b_k)$. To select the $\ell$-th item in the candidate solution, the algorithm evaluated each item not yet included by computing the incremental gain in the objective value if the item were added. It then selected the item corresponding to the $b_\ell$-th highest increment and added it to the candidate solution. In particular, the tuple $(1,\dots,1)$ corresponds to the fully greedy solution. The values of $b_1, \dots, b_k$ were tuned based on the desired number of candidate solutions.
\end{itemize}

\begin{figure}[h!]
 \captionsetup[subfigure]{font=scriptsize,labelfont=scriptsize}
  \centering
      \begin{subfigure}{0.4\textwidth}
    \centering
    \includegraphics[width=\linewidth]{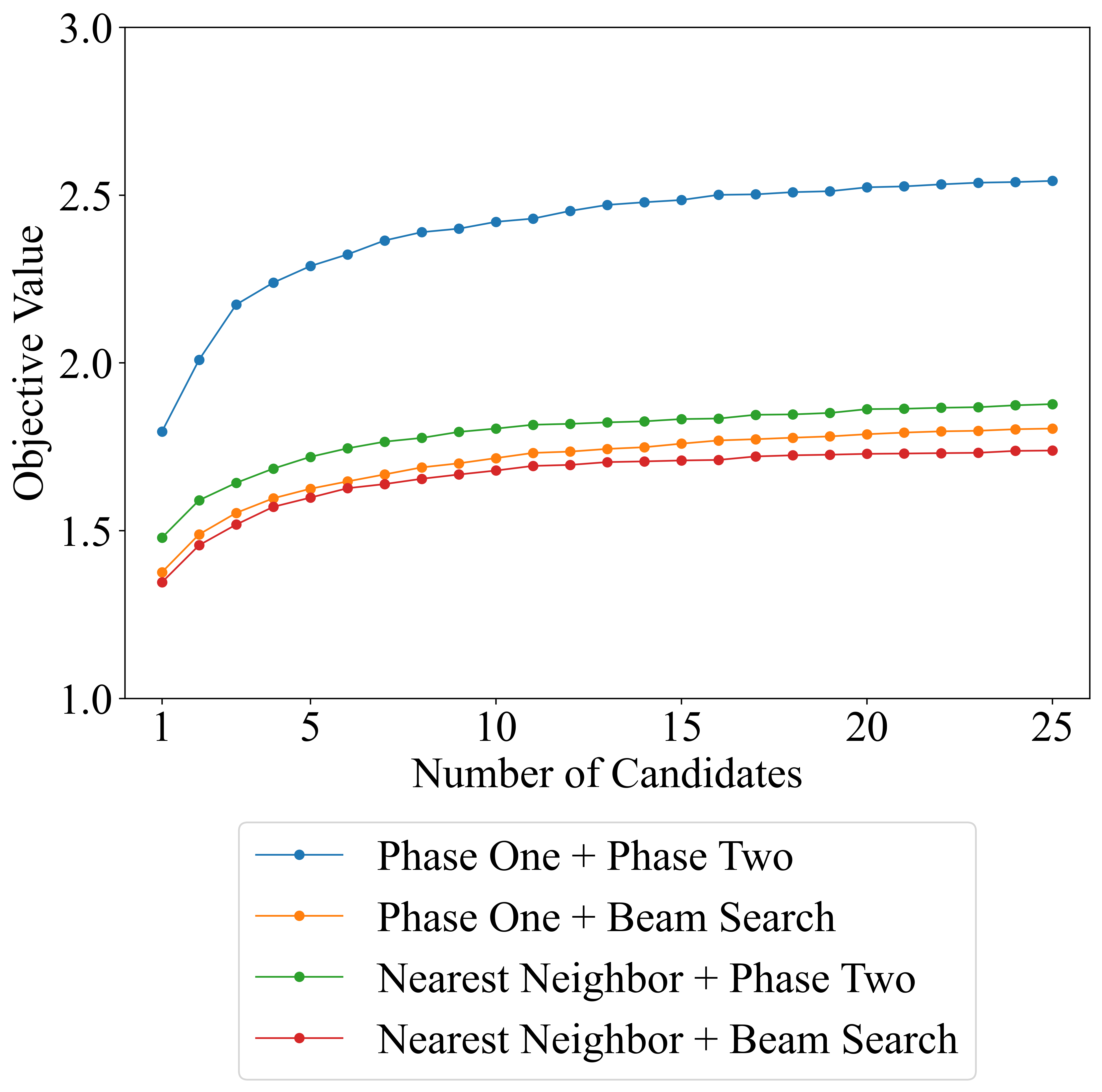}
    \caption{Spotify}
  \end{subfigure}
    \hspace{0.1\textwidth} 
  \begin{subfigure}{0.4\textwidth}
    \centering
    \includegraphics[width=\linewidth]{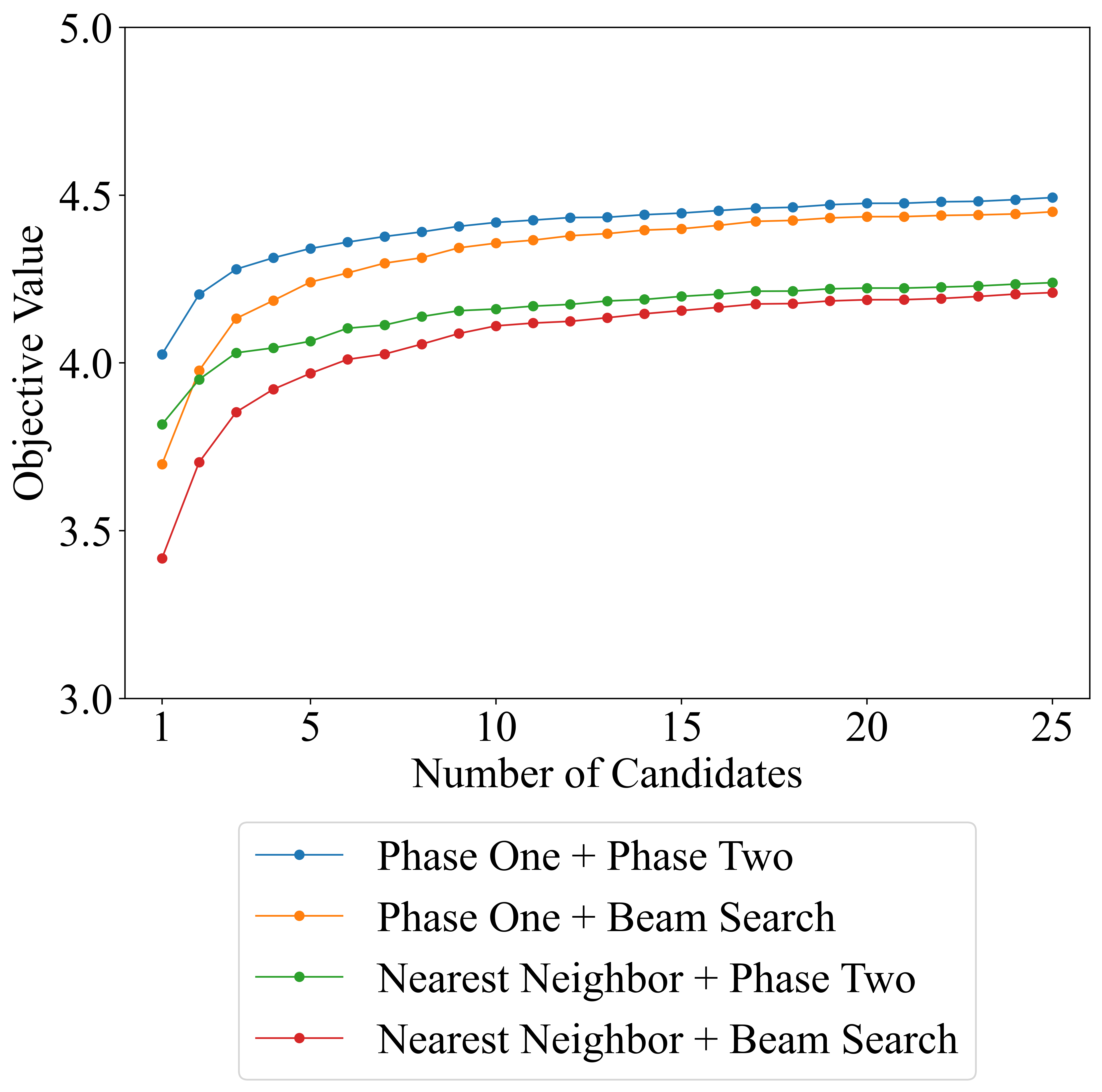}
    \caption{Trivago}
  \end{subfigure}
  \caption{Performances of four algorithms. The $x$-axis is the number of candidate solutions generated by each algorithm, and the $y$-axis is the objective value of the current best candidate solution. Each figure is averaged across 100 instances.}
  \label{fig:opt}
\end{figure}

Combining our algorithm’s retrieval and ranking phases with the two benchmark algorithms yields four distinct algorithms. The experimental results are presented in Figure~\ref{fig:opt}. Our complete algorithm consistently outperformed the other three across all settings for the number of candidate solutions generated.

In particular, when the ranking phase is fixed -- either to our algorithm’s Phase Two or to Beam Search -- our algorithm's Phase One outperformed $k$-Nearest Neighbor by an average margin of 20.86\%. This demonstrates the effectiveness of our algorithm's Phase One in selecting a superior subset of items to be passed on to the ranking phase. On the other hand, when the retrieval phase is fixed -- either to our algorithm's Phase One or to $k$-Nearest Neighbor -- our algorithm's Phase Two outperformed Beam Search by an average margin of 20.56\%. This shows our algorithm's Phase Two in ranking a given set of items is effective. In summary, our algorithm achieved both high efficiency and high accuracy in solving the personalized recommendation task, with strong empirical performance in both the retrieval and ranking phases.

We further present scatter plots to directly compare our algorithm's Phase Two with Beam Search. We fixed our algorithm’s Phase One as the retrieval phase. For the ranking phase, we proceeded as follows: we fixed the top two highest-reward items to be included in each candidate solution. Our algorithm’s Phase Two then generated candidate solutions by greedily solving a fixed number of auxiliary problems. In contrast, Beam Search generated candidate solutions by branching the same number of times and selecting the best resulting solution. As a result, the candidate solutions produced by our algorithm and Beam Search involved the same number of “iterations” and thus incurred roughly the same computational cost.

\begin{figure}[h!]
 \captionsetup[subfigure]{font=scriptsize,labelfont=scriptsize}
  \centering
      \begin{subfigure}{0.4\textwidth}
    \centering
    \includegraphics[width=\linewidth]{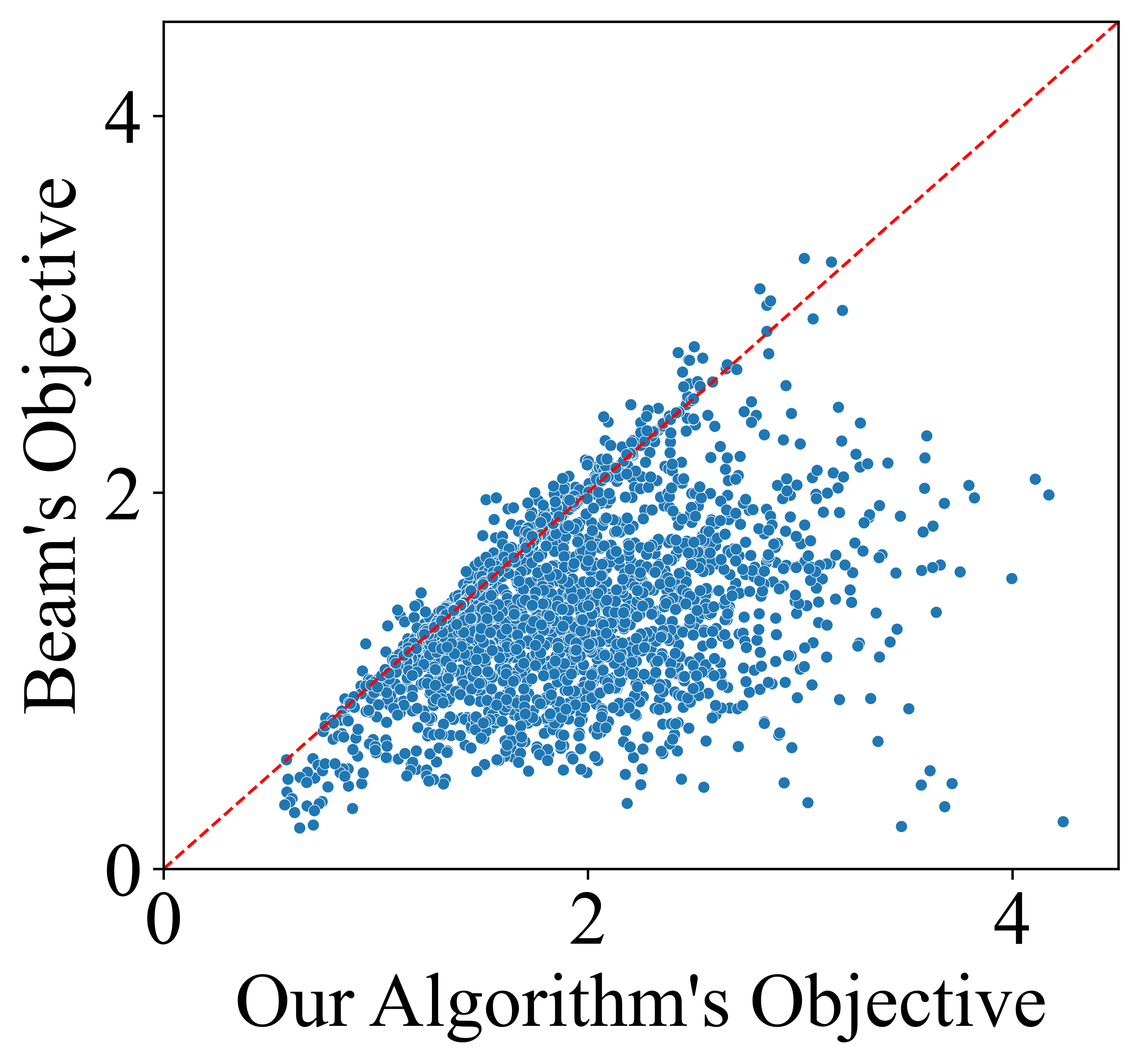}
    \caption{Spotify}
  \end{subfigure}
    \hspace{0.1\textwidth} 
  \begin{subfigure}{0.4\textwidth}
    \centering
    \includegraphics[width=\linewidth]{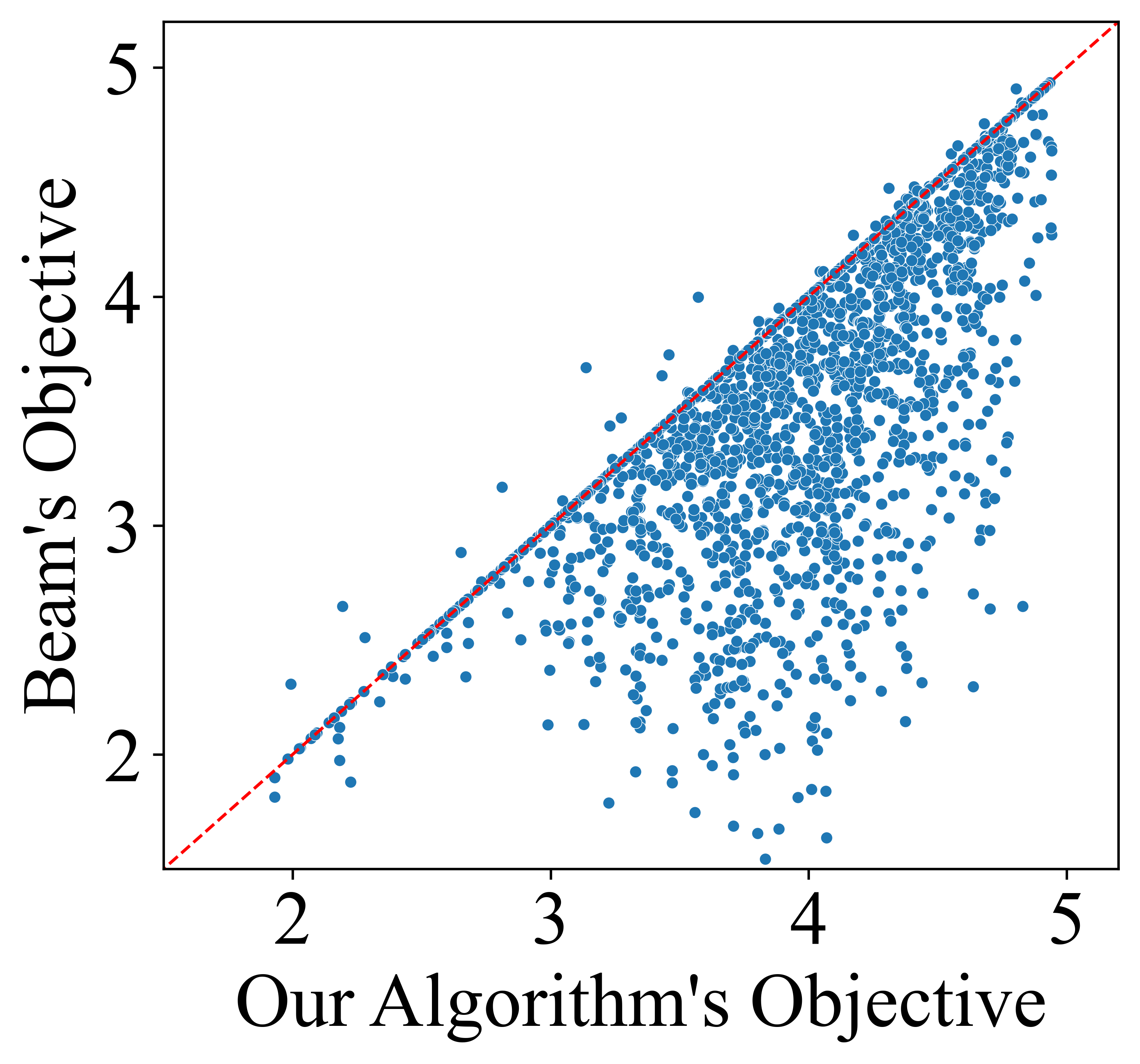}
    \caption{Trivago}
  \end{subfigure}
  \caption{Scatter plots of candidate solutions. Each point corresponds to a pair of matched candidate solutions produced by our algorithm and Beam Search. The $x$-axis represents our algorithm's objective value and the $y$-axis represents the Beam Search candidate solution's objective value. Each plot contains 2500 data points given by 25 candidate solutions in each of the 100 instances.}
  \label{fig:scatter}
\end{figure}

Each point in Figure~\ref{fig:scatter} corresponds to a pair of matched candidate solutions produced by our algorithm and Beam Search under this setup, where the $x$-axis represents the objective value of our algorithm’s solution, and the $y$-axis represents that of the Beam Search solution. Figure~\ref{fig:scatter} shows that our algorithm’s candidate solution outperformed the corresponding Beam Search solution in 90.92\% of instances, with an average improvement of 29.01\%. These results demonstrate that our algorithm’s Phase Two yields substantial improvements over Beam Search in the ranking phase.

\section{Conclusion}
In conclusion, our paper concerned the problem of real-time personalization. Traditional embedding-based machine learning models are provably unable to model certain user preferences, and recent transformer-based models are difficult to optimize in practice.  We considered a specific transformer architecture called simple transformers, which are transformers with a single self-attention layer. We proved that simple transformers were able to capture complex user preferences, such as sequential variety effect, and pairwise complementarity and substitution effects, which are essential for accurate recommendations. We then presented an algorithm that optimizes simple-transformer-based recommendation tasks, which achieves near-optimal performance with sub-linear runtime. Empirical results demonstrated that simple transformers outperformed non-transformer models in accuracy and were competitive compared to more complex transformers, and our algorithm optimized the recommendation problem with higher object value than standard benchmark algorithms like $k$-Nearest Neighbor and Beam Search.



%
%
%
\small
\setstretch{1.0}
\bibliographystyle{informs2014} 
\bibliography{Reference} 

\newpage
\appendix
\small
\setstretch{1.3}
\section{Proofs in Section \ref{Section: before transformers}} \label{Appendix A}

\begin{myproof}[Proof of Lemma \ref{lemma: nearest neighbor}]
    We prove that this can be done by applying a $\theta$-Approximate Nearest Neighbor ($\theta$-ANN) algorithm from \cite{ailon2009fast}. Given a set of $n$ ``data points'' $P\subset \mathbb{R}^{d_v}$, the goal of the $\theta$-Approximate
Near Neighbor algorithm is to build a data structure that, given a query $q \in \mathbb{R}^{d_v}$, returns a data point $p\in P$ such that $\lVert q-p\rVert_2\leq (1+\theta)\lVert q-p'\rVert_2$ for all $p'\in P$. 

\begin{lemma}[Theorem 1 in \cite{ailon2009fast}] \label{cor: nearest neighbor}
    Algorithm 1 in \cite{ailon2009fast} solves the $\theta$-ANN problem in query time $O(d\log(d)+\epsilon^{-3}\log^2(n))$. 
\end{lemma}

Moreover, the Algorithm 1 in \cite{ailon2009fast} can be easily modified to return $k$ data points $p_1,\dots,p_k\in P$ such that $\lVert q-p_i\rVert_2\leq (1+\theta)\lVert q-p'\rVert_2$ for all $i=1,\dots,k$ and all $p'\in P\setminus\{p_1,\dots,p_k\}$. This can be done by running their Algorithm 1 for $k$ times, where each time we ignore the data points that have already been returned in previous runs. More specifically, in their Algorithm 1 we change the `else if $D_T (x, p'
) > (1 + 1/2 )k '$' to `else if $D_T (x, p'
) > (1 + 1/2 )k '$ or $p'$ has been returned in previous rounds', which gives our desired performance.

Because the guarantee of Lemma \ref{lemma: nearest neighbor} is a multiplicative guarantee on the $\ell_2$ distance, below we show how to convert it into an additive guarantee on the inner products. 
We first reduce $v_1,\dots, v_n$ and $u$ to the unit sphere. Let $v_{\max}=\max_{i} \lVert v_i\rVert_2$. Let $$v'_i=\left(\sqrt{v_{\max}-\lVert v_i\rVert_2^2},v_i^\top\right)^\top/v_{\max}\in \mathbb{R}^{d+1},$$ for each $i=1,\dots,n$, where we append a scalar to $v_i$ and rescale the vector.
Let $u'=(0,u^\top)^\top/\lVert u\rVert_2$. Then $v'_i,u'\in \mathbb{S}^{d}$. Moreover, 
\begin{eqnarray*}
     \lVert u'-v_i'\rVert_2^2&=& \lVert u'\rVert_2^2+\lVert v'_i\rVert_2^2-{2v'_i}^\top u' \\
     &=&2-2  v_i^\top u/(v_{\max}||u||_2).
\end{eqnarray*}
Therefore $v_i^\top u= \frac{v_{\max}||u||_2}{2} (2-\lVert u'-v_i'\rVert_2^2)$. Hence if $|\lVert u'-v_i'\rVert_2^2-\lVert u'-v_j'\rVert_2^2|\leq\frac{2\epsilon}{v_{\max}||u||_2}$, then $|v_i^\top u-v_j^\top u|\leq\delta$. Hence, we can run the Algorithm 1 in \cite{ailon2009fast} with inputs $v_1',\dots,v_n',u'$, and $\theta=\frac{8\sqrt{2}\epsilon}{(Vu)_{\max}}$, which outputs $k$ indices $i_1,\dots,i_k$ such that $\lVert u'-v'_{i_j}\rVert_2\leq (1+\theta)\lVert u'-v'_{i_j^*}\rVert_2$ for each $j=1,\dots,k$. Then \begin{eqnarray*}
    |\lVert u'-v_{i_j}'\rVert_2^2-\lVert u'-v_{i_j^*}'\rVert_2^2|&=&\lVert u'-v_{i_j^*}'\rVert_2^2-\lVert u'-v_{i_j}'\rVert_2^2\\
    &=& (\lVert u'-v_{i_j^*}'\rVert_2+\lVert u'-v_{i_j}'\rVert_2)(\lVert u'-v_{i_j^*}'\rVert_2-\lVert u'-v_{i_j}'\rVert_2)\\
    &\leq& 4(\lVert u'-v_{i_j^*}'\rVert_2-\lVert u'-v_{i_j}'\rVert_2)\\
    &\leq&4\theta \lVert u'-v'_{i_j^*}\rVert_2\\
    &\leq& 4\sqrt{2}\epsilon\\
    &=& \frac{2\epsilon}{v_{\max}||u||_2},
\end{eqnarray*}
where the first equality follows from the definition of $i_j^*$, the second equality and the second inequality follows since $v_{i_j},v_{i^*_j},u'\in\mathbb{S}^{d_v}$. Therefore we have $v_{i_j}^\top u\geq v_{i^*_j}^\top u-\epsilon$ for each $j=1,\dots,k$ as desired.
\end{myproof}

\begin{algorithm}
\caption{$\epsilon$-Approximate $k$-Nearest Neighbor in Lemma \ref{lemma: nearest neighbor}}
{\sffamily PREPROCESS}\\

{\bf Input:} $v_1,\dots,v_n\in\mathbb{R}^{d}$ and $\epsilon>0$\;
$v_{\max} \gets \max_{i}\lVert v_i\rVert_2$\;
$v'_i \gets \left(\sqrt{v_{\max}-\lVert v_i\rVert_2^2},v_i\right)/v_{\max}$ for each $i=1,\dots,n$\;
$\theta\gets 8\sqrt{2}\epsilon/c$, where $c$ is an upper bound on $v_{\max}||u||_2$ for all possible $u$\;
Run \textsf{PREPROCESS} of the $\theta$-ANN algorithm in \cite{ailon2009fast} with inputs $v'_1,\dots,v'_n\in\mathbb{R}^{d}$ and $\theta>0$.

~\\
{\sffamily QUERY}\\

{\bf Input:} $u\in\mathbb{R}^d$\;
$u'\gets(0,u)/\lVert u\rVert_2$\;
$j\gets 1$\;
\While{$j\leq k$}{Run \textsf{PREPROCESS} of the $\theta$-ANN algorithm in \cite{ailon2009fast} with input $u'$, with the modification of changing `else if $D_T (x, p'
) > (1 + 1/2 )k '$' to `else if $D_T (x, p'
) > (1 + 1/2 )k '$ or $p'\notin\{v_{i_{j'}}\}_{j'< j}$' \;
$i_j\gets$ the index of the output of the $\theta$-ANN algorithm in \cite{ailon2009fast}\;}
Return $i_1,\dots,i_k$.
\label{alg: nearest neighbor}
\end{algorithm}

\begin{myproof}[Proof of Proposition \ref{prop: pure embedding}.] 
The algorithm \textup{ALG} operates as follows:\begin{enumerate}
    \item Partition the items according to their reward functions \( f_i \). 
    \item For each partition, apply the given \(\epsilon\)-Approximate \(k\)-Nearest Neighbor algorithm to identify \(k\) items that are collectively among the most attractive to the user within that partition.
    \item Evaluate the rewards of all such candidate items across partitions and select the top-\(k\) items with the highest overall reward.
\end{enumerate}

We analyze the performance of \textup{ALG}.  Let $\{S_1,\dots, S_\tau\}$ be a partition of the index set $[n]$ such that $f_i=f_j$ for every $\tau'\in[\tau]$ and $i,j\in S_{\tau'}$. We first construct an index set $J_{\tau'}\subset S_{\tau'}$ for each $\tau'\in[\tau]$, and then combine them to obtain the set of all candidate items $I$. 

For each index set $S_{\tau'}$, we only choose $k$ indices out of it to include in $J_{\tau'}$, namely the $k$ indices that are approximately the $k$ highest indices in $\{v_i^\top u\}_{i\in S_{\tau'}}$.\footnote{For simplicity we assume $\lvert  S_{\tau'}\rvert\geq k$. Otherwise we simply choose all indices in $S_{\tau'}$.} Specifically, for each $\ell'\in[\ell]$, we run the given $\epsilon$-Approximate $k$-Nearest Neighbor oracle with given set of points $\bigcup_{i\in S_{\tau'}}\{V_i\}\subset\mathbb{R}^{d}$, and query $u$, numbers $k$ and $\epsilon$ as inputs. We let $J_{\tau'}$ be the collection of all output indices for each $\ell'\in[\ell]$. Then $|J_{\tau'}|\leq k$ for each $\tau'$. Let $I=\bigcup_{\tau'\in[\tau]}J_{\tau'}$, then $|I|\leq \tau k$.

Let $S^*=\{i^*_1,\dots,i_k^*\}$ be the optimal solution to Problem \eqref{eqn:pure embedding} (for simplicity we assume $|S^*|=k$, and other cases can be handled similarly). We show that $\textup{ALG}_{\mathsf{Pure \text{ }Embedding}}\geq (1-g(\epsilon))\textup{OPT}_{\mathsf{Pure \text{ }Embedding}}-kh(\epsilon)$. For each $m=1,\dots,k$, let $i_m$ be the index such that $i_m$ and $i^*_m$ are in the same $S_{\tau'}$ and $v^{\top}_{i_m}u\geq v^{\top}_{i^*_m}u-\epsilon$. Then $f_{i_m}=f_{i^*_m}$. Let $S=\{i_1,\dots,i_k\}$. Because \textup{ALG} returns the $k$-highest value in $\bigcup_{i\in I}\{f_{i}(v_i^\top u)\}$, we have \begin{align*}
    \textup{ALG}_{\mathsf{Pure \text{ }Embedding}}&\geq \sum_{i \in S} f_i(v_i^\top u)\\
    &\geq \sum_{i \in S^*} f_i(v_i^\top u-\epsilon)\\
    &\geq (1-g(\epsilon))\sum_{i \in S^*} f_i(v_i^\top u)-kh(\epsilon)\\
    &=(1-g(\epsilon))\textup{OPT}_{\mathsf{Pure \text{ }Embedding}}-kh(\epsilon).
\end{align*}

Finally we analyze the expect amortized runtime of \textup{ALG}. The expect amortized runtime of constructing each $J_{\tau'}$ is 
    $k\text{-ANN}(|S_{\tau'}|,d,k,\epsilon).$ Because $|I|\leq \tau k$, the runtime of finding the $k$-highest value in $\bigcup_{i\in I}\{f_{i}(v_i^\top u)\}$ is $\log_2(\tau k)$. Therefore the expect amortized runtime of \textup{ALG} is \begin{align*}
    \sum_{\tau'=1}^\tau k\text{-ANN}(|S_{\tau'}|,d,k,\epsilon)+\log_2(\tau k)\leq \tau\cdot k\text{-ANN}\left(\frac{n}{\tau},d,k,\epsilon\right)+\log_2(\tau k).
\end{align*} The inequality follows since $\sum_{\tau'=1}^\tau |S_{\tau'}|=n$ and $k\textup{-ANN}(n,d,k,\epsilon)$ is concave in $n$.
    
\end{myproof}

\section{Proofs in Section \ref{Section Model}} \label{Appendix B}
\begin{myproof}[Proof of Proposition \ref{prop: variety and subsitution} (Model \ref{mod:variety}).]

We construct a simple transformer with the following parameters:

\begin{itemize}
    \item The input dimension is set to
    \[
    N = nk + nkd + 1 + nk,
    \]
    and is indexed as follows. Define index sets
    \[
    \mathcal{I}=\{(i,t): i\in[n],\ t\in[k]\},\quad
    \mathcal{M}=\{(j,m,a): j\in[n],\ m\in[k],\ a\in[d]\},
    \]
    \[
    \mathcal{D}=\{\odot\}\quad\text{(one dummy row)},\qquad
    \mathcal{B}=\{(i,t)^{\mathrm{base}}: i\in[n],\ t\in[k]\},
    \]
    and order the \(N\) rows of \(Q\) and \(K\) as \([\mathcal{I};\,\mathcal{M};\,\mathcal{D};\,\mathcal{B}]\).

    Before proceeding, we give intuitions on how these index sets are used:
    \begin{itemize}
\item $\mathcal I=\{(i,t)\}$:
These rows are the only rows with nonzero queries $Q$. Each of them represents the effects of the past items to the current item.

\item $\mathcal M=\{(j,m,a)\}$:
These rows have non-zero keys $K$ and scalar values $V$. They have zero queries $Q$. They encode the item $i_m$ at
position $m$ and its similarity embedding $x_{i_m,a}$.

\item $\mathcal D=\{\odot\}$:
This is a dummy row that dominates the softmax
denominator, so that the softmax vector behaves like division by a constant.

\item $\mathcal B=\{(i,t)^{\mathrm{base}}\}$:
These rows encode the base utility $\hat{u}_i$ of each item $i$.
\end{itemize}

    \item The embedding dimension is set to \(d_{kq}=k+d+2\). This is split into a \emph{position} block of length \(k\), a \emph{component} block of length \(d\), a single \emph{dummy} column (denoted \(\Delta\)), and a single \emph{base} column (denoted \(\Theta\)).

    \item Let \(M>0\) be a sufficiently large constant and \(b_0\in\mathbb{R}\) a fixed constant. Later we will take \(M\) large enough so that the simple transformer approximates \(g(S,i_t)\) to arbitrary precision.

    \item For each position \(t\in[k]\), define the row vector \(r_t\in\mathbb{R}^{1\times k}\) by
    \[
    (r_t)_m \;=\;
    \begin{cases}
    \log \lambda_{\,t-m}, & m<t,\\
    -M, & m\ge t,
    \end{cases}
    \qquad m\in[k],
    \]
    i.e.
    \[
    r_t
    =\big[ \log\lambda_{t-1},\;\; \log\lambda_{t-2},\;\; \cdots\;\; \log\lambda_{1},\;\;
    \underbrace{-M,\;\; -M,\;\; \cdots,\;\; -M}_{k-t+1} \big].
    \]
    For each \(t\in[k]\), let \(R_t\in\mathbb{R}^{n\times k}\) be the matrix with all rows equal to \(r_t\):
    \[
    R_t \;=\;
    \begin{bmatrix}
    r_t\\
    r_t\\
    \vdots\\
    r_t
    \end{bmatrix}
    \quad\text{($n$ rows).}
    \]
    Then the \emph{position block} of the query matrix is the vertical stacking of these \(k\) blocks:
    \[
    Q_{\mathrm{pos}}
    \;=\;
    \begin{bmatrix}
    R_1\\
    R_2\\
    \vdots\\
    R_k
    \end{bmatrix}
    \;\in\mathbb{R}^{(nk)\times k}.
    \]

    Next, define the \(n\times d\) matrix \(G\) collecting the (component-wise) logarithms of the similarity embeddings
    \(x_i=(x_{i,1},\dots,x_{i,d})\in S^{d-1}\):
    \[
    G \;=\;
    \begin{bmatrix}
    (\log x_1)^\top\\
    (\log x_2)^\top\\
    \vdots\\
    (\log x_n)^\top
    \end{bmatrix}
    \;=\;
    \begin{bmatrix}
    \log x_{1,1} & \log x_{1,2} & \cdots & \log x_{1,d}\\
    \log x_{2,1} & \log x_{2,2} & \cdots & \log x_{2,d}\\
    \vdots & \vdots & \ddots & \vdots\\
    \log x_{n,1} & \log x_{n,2} & \cdots & \log x_{n,d}
    \end{bmatrix}.
    \]
    The \emph{component block} of the query matrix is \(k\) identical copies of \(G\) stacked vertically:
    \[
    Q_{\mathrm{cmp}}
    \;=\;
    \begin{bmatrix}
    G\\
    G\\
    \vdots\\
    G
    \end{bmatrix}
    \;\in\mathbb{R}^{(nk)\times d}.
    \]

    Let \(\mathbf{1}_{p}\in\mathbb{R}^{p\times 1}\) be the all-ones column, and $0_{p}\in\mathbb{R}^{p\times 1}\) be the all-zeross column. Define the \emph{dummy} and \emph{base} query columns as
    \[
    Q_{\Delta} \;=\;
    \begin{bmatrix}
    M^{3}\,\mathbf{1}_{nk}\\[2pt]
    0_{nkd}\\[2pt]
    0\\[2pt]
    0_{nk}
    \end{bmatrix}
    ,\qquad
    Q_{\Theta} \;=\;
    \begin{bmatrix}
    b_0\,\mathbf{1}_{nk}\\[2pt]
    0_{nkd}\\[2pt]
    0\\[2pt]
    0_{nk}
    \end{bmatrix}
    \ \in\ \mathbb{R}^{N\times 1}.
    \]

    Putting the query blocks together, we define \(Q\) to be
    \[
    Q
    \;=\;
    \begin{bmatrix}
    Q_{\mathrm{pos}} & Q_{\mathrm{cmp}} & Q_{\Delta} & Q_{\Theta}\\[2pt]
    0 & 0 & 0 & 0\\[2pt]
    0 & 0 & 0 & 0\\[2pt]
    0 & 0 & 0 & 0
    \end{bmatrix}
    \;\in\;\mathbb{R}^{N\times (k+d+2)}.
    \]

    \item 
    The \emph{position block} of the key matrix is defined as 
    \[
    K_{\mathrm{pos}}
    \;=\;
    \begin{bmatrix}
    \mathbf{1}_{nd} & 0 & \cdots & 0\\
    0 & \mathbf{1}_{nd} & \cdots & 0\\
    \vdots & \vdots & \ddots & \vdots\\
    0 & 0 & \cdots & \mathbf{1}_{nd}
    \end{bmatrix}\;\in\mathbb{R}^{(nkd)\times k},
    \]
    i.e., a block-diagonal matrix with \(k\) diagonal blocks, each block the column \(\mathbf{1}_{nd}\).

    Let \(I_d\in\mathbb{R}^{d\times d}\) be the identity matrix. The \emph{component block} of the key matrix is defined as 
    \[
    K_{\mathrm{cmp}}
    \;=\;
    \begin{bmatrix}
    I_{d}\\
    I_d\\
    \vdots\\
    I_d
    \end{bmatrix}
    \;\in\mathbb{R}^{(nkd)\times d},
    \]
    Thus, in the block associated with position \(m\) and item \(j\), the \(d\) consecutive rows
    equal the identity \(I_d\).

    The \emph{dummy} and \emph{base} key columns are
    \[
    K_{\Delta} \;=\;
    \begin{bmatrix}
    0_{nk}\\[2pt]
    0_{nkd}\\[2pt]
    1\\[2pt]
    0_{nk}
    \end{bmatrix}
    ,\qquad
    K_{\Theta} \;=\;
    \begin{bmatrix}
    0_{nk}\\[2pt]
    0_{nkd}\\[2pt]
    0\\[2pt]
    \mathbf{1}_{nk}
    \end{bmatrix}
    \ \in\ \mathbb{R}^{N\times 1}.
    \]

  Putting the key blocks together, we define \(K\) to be
    \[
    K
    \;=\;
    \begin{bmatrix}
    0 & 0 & 0 & 0\\[2pt]
    K_{\mathrm{pos}} & K_{\mathrm{cmp}} & 0 & 0\\[2pt]
    0 & 0 & K_{\Delta} & 0\\[2pt]
    0 & 0 & 0 & K_{\Theta}
    \end{bmatrix}
    \;\in\;\mathbb{R}^{N\times (k+d+2)}.
    \]

    \item 
    Let \(v\in\mathbb{R}^{nd\times 1}\):
    \[
    v \;=\;
    \begin{bmatrix}
    x_1^\top\\
    x_2^\top\\
    \vdots\\
    x_n^\top
    \end{bmatrix}
    \;=\;
    \begin{bmatrix}
    x_{1,1}\\ \vdots\\ x_{1,d}\\ x_{2,1}\\ \vdots\\ x_{n,d}
    \end{bmatrix}.
    \]
    We define \(V\) (one scalar per row) by placing zeros on \(\mathcal{I}\) and \(\mathcal{D}\),
    the item–component entries on \(\mathcal{M}\), and the base utilities on \(\mathcal{B}\):
    \[
    V \;=\;
    \begin{bmatrix}
    0_{nk}\\[2pt]
    \beta\,e^{M^3}\,v\\[2pt]
    0\\[2pt]
    e^{M^3}\,\big(\tfrac{1}{\beta}\log\hat u_1,\dots,\tfrac{1}{\beta}\log\hat u_n\big)^\top\\[2pt]\vdots\\[2pt]
    e^{M^3}\,\big(\tfrac{1}{\beta}\log\hat u_1,\dots,\tfrac{1}{\beta}\log\hat u_n\big)^\top
    \end{bmatrix}
    \ \in\ \mathbb{R}^{N\times 1},
    \] where $e^{M^3}\,\big(\tfrac{1}{\beta}\log\hat u_1,\dots,\tfrac{1}{\beta}\log\hat u_n\big)^\top$ is repeated $k$ times.
    Set \(u=1\) so that \(V_i^\top u=V_i\) for each row \(i\).
\item We set \(f_i(x)=\exp(\beta x)\) for all $i\in[nk+nkd+1+nk]$.
    
    \item For a length-$k$ sequence $S=(i_1,\dots,i_k)$, let $S'\subset [nk+nkd+1+nk]=[\mathcal{I};\,\mathcal{M};\,\mathcal{D};\,\mathcal{B}]$ such that $S'$ contains $\{(i_t,t):t\in[k]\}\subset[\mathcal{I}]$ from the first set of rows,
$\{(i_m,m,a):m\in[k],\,a\in[d]\}\subset[\mathcal{M}]$ from the second set of rows, the dummy row $\mathcal{D}$,
and $\{(i_t,t)^{\mathrm{base}}:t\in[k]\}\subset\mathcal{B}$ from the set of base rows. 
\end{itemize}

We show that the simple transformer above approximates $g(S,i_t)$ to arbitrary precision. The intuition of our construction is given below:
\begin{itemize}
    \item $Q$ encodes lags ($\lambda$) and similarity embeddings ($\log x_i$).
    \item $K$ encodes positions ($t$) in the sequence.
    \item $V$ encodes similarity embeddings ($\log x_i$) and base utilities ($\hat{u}_{i}$). 
    \item The dummy row makes softmax act like a constant divider.
\end{itemize}

Fix a length-$k$ sequence $S=(i_1,\dots,i_k)$. We first compute $X_{S'}Q(X_{S'}K)^\top$ for each $(i_t,t)\in S'$:
\[
\begin{aligned}
\text{(context)}\quad &(X_{S'}Q(X_{S'}K)^\top)_{(i_t,t),(i_m,m,a)} \;=\; r_t(m) + \log x_{i_t,a}
=
\begin{cases}
\log\lambda_{t-m} + \log x_{i_t,a}, & m<t,\\
-\,M + \log x_{i_t,a}, & m\ge t,
\end{cases}\\[2pt]
\text{(dummy)}\quad &(X_{S'}Q(X_{S'}K)^\top)_{(i_t,t),\odot} \;=\; M^3,\\[2pt]
\text{(base)}\quad &(X_{S'}Q(X_{S'}K)^\top)_{(i_t,t),(i_t,t)^{\mathrm{base}}} \;=\; 0.
\end{aligned}
\]
Hence the denominator of the softmax operation on row $X_{S'}Q(X_{S'}K)^\top$ is
\[
Z_t \;=\; e^{M^3} + 1
+ \sum_{m<t}\sum_{a=1}^d \lambda_{t-m}x_{i_t,a}
+ \sum_{m\ge t}\sum_{a=1}^d e^{-M}x_{i_t,a}.
\]
Define the following (finite) constants
\[
\Lambda_t=\sum_{m<t}\lambda_{t-m},\qquad
S^{(\le)}_t=\sum_{m<t}\sum_{a} \lambda_{t-m}x_{i_t,a}\ ,
\qquad
S^{(\ge)}_t=\sum_{m\ge t}\sum_{a} e^{-M}x_{i_t,a}\ .
\]
Set \[
\delta_t = e^{-M^3}+ e^{-M^3}S_t^{(\le)}+ e^{-M^3}S_t^{(\ge)}.\]
Then
\[
Z_t = e^{M^3}(1+\delta_t).
\]
Therefore the attention weights on row \((i_t,t)\) are
\[
\textup{softmax}(X_{S'}Q(X_{S'}K)^\top)_{(i_t,t),(i_m,m,a)}
=\frac{\exp(r_t(m)+\log x_{i_t,a})}{Z_t}
=\frac{\lambda_{t-m}\,x_{i_t,a}}{e^{M^3}(1+\delta_t)}\quad(m<t),
\]
and \[
\textup{softmax}(X_{S'}Q(X_{S'}K)^\top)_{(i_t,t), (i_t,t)^{\mathrm{base}}}
=\frac{\exp(0)}{Z_t}=\frac{e^{-M^3}}{1+\delta_t}.
\]
Recall that
\[
V_{(i_m,m,a)} = e^{M^3}\,x_{i_m,a},\qquad
V_\odot = 0,\qquad
V_{(i_t,t)^{\mathrm{base}}} = e^{M^3}\cdot \frac{1}{\beta}\log\hat u_{i_t}.
\]
So we have
\[
\begin{aligned}
\mathrm{SA}_{Q,K,V}(X_{S'})_{(i_t,t)}
&=\sum_{m<t}\sum_{a=1}^d \textup{softmax}(X_{S'}Q(X_{S'}K)^\top)_{(i_t,t), (i_m,m,a)}\,V_{(i_m,m,a)}
\\&\quad \text{ }+\textup{softmax}(X_{S'}Q(X_{S'}K)^\top)_{(i_t,t), (i_t,t)^{\mathrm{base}}}\,V_{(i_t,t)^{\mathrm{base}}}
\\[2pt]
&=\frac{1}{1+\delta_t}
\left(
\sum_{m<t}\sum_{a=1}^d \frac{\lambda_{t-m}x_{i_t,a}}{e^{M^3}}\,e^{M^3}x_{i_m,a}
+ e^{-M^3}\, e^{M^3}\cdot \frac{1}{\beta}\log\hat u_{i_t}
\right)
\\[2pt]
&=\frac{1}{1+\delta_t}
\left(
\sum_{m<t}\lambda_{t-m}\sum_{a=1}^d x_{i_t,a}x_{i_m,a}
\;+\; \frac{1}{\beta}\log\hat u_{i_t}
\right)
\\[2pt]
&=\frac{1}{1+\delta_t}
\left(
\sum_{\ell=1}^{t-1}\lambda_{\ell}\, x_{i_t}^\top x_{i_{t-\ell}}
\;+\; \frac{1}{\beta}\log\hat u_{i_t}
\right).
\end{aligned}
\]

Fix any $\epsilon>0$. Since \(S_t^{(\le)}\le \Lambda_t\le \Lambda_{k-1}\) and \(S_t^{(\ge)}\le d\,e^{-M}\),
we can set \(M\) large enough so that
\[
e^{-M^3}\big(1+\Lambda_{k-1}\big)+d\,e^{-(M^3+M)} \le \epsilon,
\] which gives $\delta_t\leq \epsilon$ for every $t$. Finally, we have 
\[
\begin{aligned}
\mathcal{T}_{Q,K,V,f_1,\dots,f_n,u}(X_{S'})_{(i_t,t)}&=f_{(i_t,t)}\big(\mathrm{SA}_{Q,K,V}(X_{S'})_{(i_t,t)}\big)\\
&=\exp\!\left(
\frac{\beta}{1+\delta_t}\sum_{\ell=1}^{t-1}\lambda_{\ell}\, x_{i_t}^\top x_{i_{t-\ell}}
+\frac{1}{1+\delta_t}\log\hat u_{i_t}
\right)
\\[2pt]
&=\hat u_{i_t}^{\,\frac{1}{1+\delta_t}}\;
\exp\!\left(\frac{\beta}{1+\delta_t}\sum_{\ell=1}^{t-1}\lambda_{\ell}\, x_{i_t}^\top x_{i_{t-\ell}}\right).
\end{aligned}
\]
Hence, as \(M\to\infty\) (so \(\delta_t\to 0\) uniformly in \(t\)), we have
\[(1-\epsilon)g(S,i_t)\leq\mathcal{T}_{Q,K,V,f_1,\dots,f_n,u}(X_{S'})_{(i_t,t)}\leq \mathcal{T}_{Q,K,V,f_1,\dots,f_n,u}(X_{S'})_{(i_t,t)}.
\] Therefore the simple transformer approximate $g(S,i_t)$ to arbitrary precision.

\end{myproof}

\begin{myproof}[Proof of Proposition \ref{prop: variety and subsitution} (Model~\ref{mod:com/sub}).]

We construct a simple transformer with three self-attention heads. The first head and the second head represent the complementarity effects and the substitution effects, respectively. They have input dimensions $n+1$, output dimension $1$, and embedding dimension $n+1$. The third head represents the base utility of each item. It has input dimensions $n$, output dimension $d$, and embedding dimension $n$. Note that this multi-head construction can be equivalently represented as a single-head simple transformer by arranging the $Q$, $K$, $V$, and $u$ matrices for all three heads into block form. However, for clarity, we present the proof by describing each head separately.

Let $A,B\in\mathbb{R}^{n\times n}_{+}$ be two matrices with positive entries such that $H_{ij}=\exp(A_{ij})-\exp(B_{ij})$ for every $i,j\in [n]$. Let $M$ be a sufficiently large constant. Later we will set $M$ to be large enough so that the simple transformer approximates $g(S,i)$ to arbitrary precision.

\paragraph{Head 1.} The first self-attention head has the following parameters:
\begin{itemize}
    \item $Q^{(1)},K^{(1)}\in \mathbb{R}^{(n+1)\times (n+1)}$ are set such that \[ Q^{(1)}(K^{(1)})^\top =
\begin{bNiceArray}{cw{c}{1cm}c|c}[margin]
\Block{3-3}{A} & & & M \\
& & & \Vdots \\
& & & M \\
\hline
0 & \Cdots& 0 & M
\end{bNiceArray}.
\]

\item $V^{(1)}\in \mathbb{R}^{n+1}$ is set to $V^{(1)}_i=1$ for each $i\in[n]$ and  $V^{(1)}_{n+1}=0$. Set $u^{(1)}=1$ so that $({V_i^{(1)}})^\top u^{(1)}=V_i^{(1)}$.

\item Set $f^{(1)}_i(x)=\exp(M)x$ for every $i\in[n^2+1]$.

\item For a subset $S\subset [n]$, we set $S^{(1)}\subset[n+1]$ where $S^{(1)}=S\cup\{n+1\}$.

\end{itemize}

We have 
\[
\textup{softmax}\big(X_{S^{(1)}}Q^{(1)} (X_{S^{(1)}} K^{(1)})^\top\big)_{ij}
=
\begin{cases}
\dfrac{\exp({A_{ij}})}{\sum_{j\in S} \exp({A_{ij}}) + \exp(M) + (n-|S|)}, & i\in S,\ j\in S,\\[10pt]
\dfrac{\exp({M})}{\sum_{j\in S} \exp({A_{ij}}) + \exp(M) + (n-|S|)}, & i\in S,\ j=n+1,\\[10pt]
\dfrac{1}{\sum_{j\in S} \exp({A_{ij}}) + \exp(M) + (n-|S|)}, & i\in S,\ j\in[n]\setminus S,\\[10pt]
\dfrac{1}{n+\exp(M)}, & i=n+1,\ j\in[n],\\[10pt]
\dfrac{\exp(M)}{n+\exp(M)}, & i=n+1,\ j=n+1,\\[10pt]
\dfrac{1}{n+1}, & i\in[n]\setminus S,\ j\in[n{+}1].
\end{cases}
\]
Therefore, for every $i\in S$, we have
\[\textup{SA}_{Q^{(1)},K^{(1)},V^{(1)}}(X_{S^{(1)}})_i^\top u^{(1)}=\dfrac{\sum_{j\in S}\exp({A_{ij}})}{\sum_{j\in S} \exp({A_{ij}}) + \exp(M) + (n-|S|)}.\]
For any given $\epsilon>0$, we can take $M$ to be sufficiently large such that 
\[\dfrac{\sum_{j\in S}\exp({A_{ij}})-\frac{\epsilon}{2}}{\exp(M) }\leq\textup{SA}_{Q^{(1)},K^{(1)},V^{(1)}}(X_{S^{(1)}})_i^\top u^{(1)}\leq \dfrac{\sum_{j\in S}\exp({A_{ij}})}{\exp(M) }.\]
Finally, we get
\[\sum_{j\in S}\exp({A_{ij}})-\frac{\epsilon}{3}\leq \mathcal{T}_{Q^{(1)},K^{(1)},V^{(1)},f^{(1)}_1,\dots,f^{(1)}_n,u^{(1)}}(X_{S^{(1)}})_i\leq \sum_{j\in S}\exp({A_{ij}}).\]

\paragraph{Head 2.} The second self-attention head has exactly the same parameters, expect we replace $A$ with $B$ and we flip the sign of $f_i$:
\begin{itemize}
    \item $Q^{(2)},K^{(2)}\in \mathbb{R}^{(n+1)\times (n+1)}$ are set such that \[ Q^{(2)}(K^{(2)})^\top =
\begin{bNiceArray}{cw{c}{1cm}c|c}[margin]
\Block{3-3}{B} & & & M \\
& & & \Vdots \\
& & & M \\
\hline
0 & \Cdots& 0 & M
\end{bNiceArray}.
\]

\item $V^{(2)}\in \mathbb{R}^{n+1}$ is set to $V^{(2)}_i=1$ for each $i\in[n]$ and  $V^{(2)}_{n+1}=0$. Set $u^{(2)}=1$ so that $({V_i^{(1)}})^\top u^{(2)}=V_i^{(2)}$.

\item Set $f^{(2)}_i(x)=-\exp(M)x$ for every $i\in[n^2+1]$.

\item For a subset $S\subset [n]$, we set $S^{(2)}\subset[n+1]$ where $S^{(2)}=S\cup\{n+1\}$.

\end{itemize}

Then, similar to head 1, we get 
\[-\sum_{j\in S}\exp({B_{ij}})\leq \mathcal{T}_{Q^{(2)},K^{(2)},V^{(2)},f^{(2)}_1,\dots,f^{(2)}_n,u^{(2)}}(X_{S^{(2)}})_i\leq -\sum_{j\in S}\exp({B_{ij}})+\frac{\epsilon}{3}.\]

\paragraph{Head 3.} The third self-attention head has the following parameters:
\begin{itemize}
    \item $Q^{(3)},K^{(3)}\in \mathbb{R}^{n\times n}$ are set such that \[Q^{(3)}(K^{(3)})^\top =\begin{bmatrix}
0 & -M & \cdots & -M \\
-M & 0 & \cdots & -M \\
\vdots & \vdots & \ddots & \vdots \\
-M & -M & \cdots & 0
\end{bmatrix}.\]

\item $V^{(3)}\in \mathbb{R}^{n\times d}$ is set to $V^{(3)}_i=\hat{v}_i$ for each $i\in[n]$. 

\item $u^{(3)}\in\mathbb{R}^d$ is set to $u^{(3)}=\hat{u}$.

\item Set $f^{(3)}_i(x)=x$ for every $i\in[n]$.

\item For a subset $S\subset [n]$, we set $S^{(3)}=S$.
\end{itemize}

For every $i\in S$, we have
\[\textup{SA}_{Q^{(3)},K^{(3)},V^{(3)}}(X_{S^{(3)}})_i^\top u^{(3)}=\dfrac{(\hat{v_i})^\top \hat{u}+(n-|S|)\exp(-M)}{1 + (n-|S|)\exp(-M)}.\]
For any given $\epsilon>0$, we can take $M$ to be sufficiently large such that 
\[(\hat{v_i})^\top \hat{u}-\frac{\epsilon}{3}\leq \mathcal{T}_{Q^{(3)},K^{(3)},V^{(3)},f^{(3)}_1,\dots,f^{(3)}_n,u^{(3)}}(X_{S^{(3)}})_i\leq (\hat{v_i})^\top \hat{u}.\]

\paragraph{Complete the Proof.} Because $\exp(A_{ij})-\exp(B_{ij})=H_{ij}$ Combining the above three self-attention heads, we have
\[\hat{v}_i^\top \hat{u}+\sum_{j\in S}H_{ij}-\epsilon\leq \sum_{\ell=1}^3\mathcal{T}_{Q^{(\ell)},K^{(\ell)},V^{(\ell)},f^{(\ell)}_1,\dots,f^{(\ell)}_n,u^{(\ell)}}(X_{S^{(\ell)}})_i\leq \hat{v}_i^\top \hat{u}+\sum_{j\in S}H_{ij}.\]
Therefore the three attention heads approximate $g(S,i)$ to arbitrary precision.

\end{myproof}

\section{Proofs in Section \ref{Section Hardness}.} \label{Appendix C}

\begin{myproof}[Proof of Proposition \ref{prop: lower bound} (a).]
Fix any instance of the $(k-1)$-\textsc{Clique} problem: let $G$ be the (undirected, unweighted) graph with $n-1$ vertices, 
and let $A\in \{0,1\}^{(n-1)\times (n-1)}$ be its adjacency matrix (where we follow the convention that diagonal entries are set equal to 1). We will create an explicit instance of Problem \eqref{eqn:problem} 
in the following way (using the formulation in \ref{eqn: simplification}):
\begin{itemize}
    \item $K=I_{n\times n}$, i.e.~the $n\times n$ identity matrix, so that  $W=\textup{softmax}(Q K^\top)=\textup{softmax}(Q)$.
    \item $Q\in \mathbb{R}^{n\times n}$ is set according to 
$$Q_{ij}=\begin{cases}
0& \text{for } i=n \text{ or }j=n\\
-N A_{ij} & \text{for } i,j \neq n ,
\end{cases}$$
where $N > 1$ is a constant large enough that $(k-1)\exp(-N)\leq 1/2$. 
In block notation, this is \[ Q =
\begin{bNiceArray}{cw{c}{1cm}c|c}[margin]
\Block{3-3}{-N \cdot A} & & & 0 \\
& & & \Vdots \\
& & & 0 \\
\hline
0 & \Cdots& 0 & 0
\end{bNiceArray}.
\]

\item $d_v = 1$, and we set $u=1$, so that $Vu=V$.

\item $V\in\mathbb{R}^{n\times 1}$ is set according to $V_i=1$ for $i =1,\dots,n-1$, and $V_n=M$, where we take any constant $M\geq 2$.

\item For $i=1,\dots, n-1$, we set $f_i(\cdot)$ to be:
$$f_i(x)=\begin{cases}
0& \text{if }0\leq x\leq  \frac{M+1}{2}\\
\frac{6x - 3(M+1)}{M - 1} & \text{if } \frac{M+1}{2}<x< \frac{2M+1}{3}\\
1& \text{if } x\geq \frac{2M+1}{3}.
\end{cases}$$ 
Because $M\geq 2$, we have $(2M+1)/3>(M+1)/2$, so $f_i(\cdot)$ is continuous piece-wise linear for every $i=1,\dots,n-1$. We set $f_n(x)=0$ for every $x$.

\end{itemize}

Note that the quantity $$\frac{(w_i\odot Vu)^\top x}{ w_i^\top x}$$ remains unchanged if the vector $w_i$ is multiplied by a non-zero constant. Thus, rescaling the rows of $W$ does not change $\mathsf{P}$. So for simplicity of exposition, we replace $W$ with $W'$, defined to be
$$W'_{ij}=\begin{cases}
1& \text{for } i=n \text{ or }j=n\\
\exp(-N A_{ij})& \text{for } i,j \neq n.
\end{cases}$$ As a sanity check, $W$ is simply $W'$ with each row rescaled to sum to one.

Now notice that because $f_i(x)\leq 1$ for every $i=1,\dots,n-1$, and $f_n(x)=0$, the optimal value of this instance of $\mathsf{P}$ is at most $k-1$. It will suffice to show that $G$ has a clique of size $k-1$ if and only if the optimal value is $k-1$. We prove both directions separately.


\paragraph{If $G$ has a clique of size $k-1$, then the optimal value is $k-1$:} Suppose $G$ has a clique of size $k-1$, and let $C\subset [n-1]$ be the vertex set of one such clique. 
Consider the solution $x^*$, where $x^*_i=1$ if and only if $i\in C\cup\{n\}$. We will show that the objective value at $x^*$ is $k-1$ (and thus the optimal value is $k-1$): 

Because $W'_{ij}=\exp(-N)$ whenever $i,j\in C$, we have \begin{eqnarray*}
    \sum_{i=1}^{n} x^*_i f_i\left(\frac{(w'_i\odot Vu)^\top x^*}{{w'_i}^\top x^*}\right)&=&
    \sum_{i\in C} f_i\left(\frac{(w'_{i}\odot Vu)^\top x^*}{{w'_{i}}^\top x^*}\right)\\
    &=&
    \sum_{i\in C}f_i\left(\frac{M+(k-1)\exp(-N)}{1+(k-1)\exp(-N)}\right)
    .
\end{eqnarray*}
Since $N$ was chosen to be large enough that $(k-1)\exp(-N)\leq 1/2$, it follows that $$\frac{M+(k-1)\exp(-N)}{1+(k-1)\exp(-N)}\geq \frac{M+\frac{1}{2}}{1+\frac{1}{2}}=\frac{2M+1}{3}.$$ 
Therefore, 
\begin{eqnarray*}
    \sum_{i=1}^{n} x^*_i f_i\left(\frac{(w'_i\odot Vu)^\top x^*}{{w'_i}^\top x^*}\right)
    =
    \sum_{i\in C}f_i\left(\frac{M+(k-1)\exp(-N)}{1+(k-1)\exp(-N)}\right)
    =|C|=k-1
    .
\end{eqnarray*}

\paragraph{If the optimal value is $k-1$, the $G$ has a clique of size $k-1$:} Suppose the optimal value is $k-1$. Let $x^*$ be an optimal solution, and let $C\subset [n-1]$ be the index set such that $x^*_{i}=1$ for $i\in C$. Notice that if a solution $x$ has $x_n=0$, then since $(w'_i\odot Vu)_j=W'_{ij}$ for all $j\in[n-1]$, we have $$\sum_{i=1}^{n} x_i f_i\left(\frac{(w'_i\odot Vu)^\top x}{{w_i'}^\top x}\right)=\sum_{i=1}^{n} x_i f_i\left(1\right)=0.$$ 
Therefore, it must be the case that $x^*_n=1$. Also, because $$k-1=\sum_{i=1}^{n} x^*_i f_i\left(\frac{(w'_i\odot Vu)^\top x^*}{{w_i'}^\top x^*}\right)\leq \sum_{i=1}^{n-1}x^*_i=|C|,$$ we must have $|C|=k-1$ and $$f_i\left(\frac{(w'_i\odot Vu)^\top x^*}{{w_i'}^\top x^*}\right)=1$$ for every $i\in C$.

Let $d(i)$ be the degree of vertex $i$ in the induced subgraph of $G$ with vertex set $C$. Then \begin{eqnarray*}
    \sum_{i=1}^{n} x^*_i f_i\left(\frac{(w'_i\odot Vu)^\top x^*}{{w'_i}^\top x^*}\right)&=&
    \sum_{i\in C} f_i\left(\frac{(w'_{i}\odot Vu)^\top x^*}{{w'_{i}}^\top x^*}\right)\\
    &=&
    \sum_{i\in C}f_i\left(\frac{M+(k-2-d(i))+\exp(-N)(d(i)+1)}{1+(k-2-d(i))+\exp(-N)(d(i)+1)}\right)\\
    &\leq&\sum_{i\in C}f_i\left(\frac{M+(k-2-d(i))}{1+(k-2-d(i))}\right),
\end{eqnarray*} where the last inequality follows since $M\geq 2$.
Because $|C|=k-1$, we have $d(i)\leq k-2$. Suppose for the sake of contradiction that $d(i)\leq k-3$ for some $i\in I$. Then $$f_i\left(\frac{M+(k-2-d(i))}{1+(k-2-d(i))}\right)\leq f_i\left(\frac{M+1}{2}\right)=0,$$ which contradicts that $$f_i\left(\frac{(w'_i\odot Vu)^\top x^*}{{w_i'}^\top x^*}\right)=1$$ for every $i\in C$. Therefore we must have $d(i)=k-2$ for every $i\in C$. Hence $C$ corresponds to the vertex set of a clique of size $k-1$.
\end{myproof}

\begin{myproof}[Proof of Proposition \ref{prop: lower bound} (b).]
Fix any constant $M\geq 1$. We construct an instance of Problem \eqref{eqn:problem} (using the fornulation in $\mathsf{P}$) by applying the Johnson-Lindenstrauss Lemma:
\begin{lemma}[Johnson-Lindenstrauss Lemma] \label{lemma: JL} For any $0<\epsilon<1$ and any set $S$ of $m$ points in $\mathbb{R}^n$, there exists a universal constant $c>0$ and a linear function $f:\mathbb{R}^n\to\mathbb{R}^d$ with $d=c\epsilon^{-2}\log(m)$ such that $$(1-\epsilon)\lVert x_i\rVert_2^2\leq \lVert f(x_i)\rVert_2^2\leq (1+\epsilon)\lVert x_i\rVert_2^2$$
for all $x_i\in S$ and $$(1-\epsilon)\lVert x_i-x_j\rVert_2^2\leq \lVert f(x_i)-f(x_j)\rVert_2^2\leq (1+\epsilon)\lVert x_i-x_j\rVert_2^2$$ for all $x_i,x_j\in S$.
\end{lemma}

Take $N$ to be large enough such that $\exp(-N)\leq 1/2M$. Take $0<\delta<1$ to be small enough such that $\exp(-\delta N)\geq1-\exp(-N)$. Then $N$ and $\delta$ can be chosen to be both only depend on $M$. We obtain the following corollary:
\begin{corollary}[Corollary of Lemma \ref{lemma: JL}.]\label{cor: JL}
    There exists a number $c(M)>0$ and a set $S$ of $\exp(c(M)\cdot d_{kq})$ unit vectors in $\mathbb{R}^{d_{kq}}$ such that  $|u_i^\top u_j|\leq \delta$ for every $u_i,u_j\in S$ and $u_i\neq u_j$.
\end{corollary}
\begin{myproof}[Proof of Corollary \ref{cor: JL}.]
    Let $c>0$ be the universal constant in Lemma \ref{lemma: JL} and let $n=\exp(c^{-1}(\delta/4)^2\cdot d_{kq})$. Let $c(M)=c^{-1}(\delta/4)^2$, then $n=\exp(c(M)\cdot d_{kq})$. Because $\delta>0$ only depends on $M$, we have $c(M)>0$ also only depends on $M$. Consider $\{e_i\}_{i=1}^n\subset \mathbb{R}^n$ where $e_i\in\mathbb{R}^n $ is the unit vector where the $i$-th entry equals to 1 and all other entries equal to 0. Then we have $\lVert e_i-e_j\rVert_2^2=2$ for every $e_i\neq e_j$. By Lemma \ref{lemma: JL}, there exists a linear function $f:\mathbb{R}^n\to\mathbb{R}^{d_{kq}}$ such that $$\left(1-\frac{\delta}{4}\right)\leq \lVert f(e_i)\rVert_2^2\leq \left(1+\frac{\delta}{4}\right)$$
for all $e_i$ and $$2\left(1-\frac{\delta}{4}\right)\leq \lVert f(e_i)-f(e_j)\rVert_2^2\leq 2\left(1+\frac{\delta}{4}\right)$$ for all $e_i\neq e_j$.

For every $e_i\neq e_j$, because $$
\|f(e_i) - f(e_j)\|_2^2 = \|f(e_i)\|_2^2 + \|f(e_j)\|_2^2 - 2f(e_i)^\top f(e_j),
$$ we have \begin{align*}
    f(e_i)^\top f(e_j)&=(\|f(e_i)\|_2^2 + \|f(e_j)\|_2^2-\|f(e_i) - f(e_j)\|_2^2)/2\\
    &\leq \left(\left(1+\frac{\delta}{4}\right) + \left(1+\frac{\delta}{4}\right)-2\left(1-\frac{\delta}{4}\right)\right)/2\\
    &=\frac{\delta}{2},
\end{align*}
and \begin{align*}
    f(e_i)^\top f(e_j)&=(\|f(e_i)\|_2^2 + \|f(e_j)\|_2^2-\|f(e_i) - f(e_j)\|_2^2)/2\\
    &\geq \left(\left(1-\frac{\delta}{4}\right) + \left(1-\frac{\delta}{4}\right)-2\left(1+\frac{\delta}{4}\right)\right)/2\\
    &=-\frac{\delta}{2}.
\end{align*}
Let $u_i=f(e_i)/\|f(e_i)\|_2$ for all $i\in[n]$, then each $u_i$ is a unit vector. Moreover, for every $u_i\neq u_j$, \begin{align*}
    u_i^\top u_j=\frac{f(e_i)^\top f(e_j)}{\|f(e_i)\|_2\|f(e_j)\|_2}\leq\frac{\frac{\delta}{2}}{1-\frac{\delta}{4}}=\delta\left(\frac{2}{4-\delta}\right)<\delta,
\end{align*} where the last inequality follows since $0<\delta<1$. Similarly, \begin{align*}
    u_i^\top u_j=\frac{f(e_i)^\top f(e_j)}{\|f(e_i)\|_2\|f(e_j)\|_2}\geq\frac{-\frac{\delta}{2}}{1-\frac{\delta}{4}}=-\delta\left(\frac{2}{4-\delta}\right)>-\delta.
\end{align*} Therefore  $|u_i^\top u_j|\leq \delta$ for every  $u_i\neq u_j$. Take $S=\{u_i\}_{i=1}^n$ gives the desired set of unit vectors.
\end{myproof}

Corollary \ref{cor: JL} states that there exists a set of unit vectors in $\mathbb{R}^{d_{kq}}$, with size exponential in $d_{kq}$, where all unit vectors in the set are approximately orthonormal. 

Fix any $d_{kq}$ such that $\exp(c(M)\cdot d_{kq})\leq n-1$ and any $k\geq M+1$. Let $G$ be a graph with $n-1$ vertices $v_1,\dots,v_{n-1}$ and $\ell=\exp(c(M)\cdot d_{kq})$ disjoint cliques, each of size at least $k-M$ and at most $k-1$. Let $I_1,\dots,I_\ell\subset [n-1]$ be the index sets of vertices corresponding to these $\ell$ cliques. That is, $\{v_i\}_{i\in I_{\ell'}}$ forms a clique for each $\ell'\in[\ell]$. Without loss of generality we assume $\ell'\in I_{\ell'}$ for every $\ell'\in[\ell]$. 

By Corollary \ref{cor: JL}, there exists $\ell$ unit vectors $u_1,\dots, u_{\ell}\in\mathbb{R}^{d_{kq}}$ such that $|u_i^\top u_j|\leq \delta$ for $i\neq j$. Let $u_{\ell+1},\dots,u_{n-1}$ be unit vectors such that $u_i=u_{\ell'}$ for $i\in I_{\ell'}$. That is, for indices $i,j$ such that $v_i$ and $v_j$ are in the same clique, we have $u_i=u_j$. Let $U\in\mathbb{R}^{(n-1)\times d_{kq}}$ where the $i$-th row of $U$ is $u_i^\top$. Let $A=UU^\top\in\mathbb{R}^{(n-1)\times (n-1)}$,  then $A$ has rank at most $d_{kq}$. 

Because $A$ has rank at most $d_{kq}$, there exists $Q,K\in\mathbb{R}^{n\times d_{kq}}$ such that \[ Q K^\top =
\begin{bNiceArray}{cw{c}{1cm}c|c}[margin]
\Block{3-3}{-N \cdot A} & & & 0 \\
& & & \Vdots \\
& & & 0 \\
\hline
0 & \Cdots& 0 & 0
\end{bNiceArray}.
\] We create an
explicit instance of $\mathsf{P}$ similar to the instance in the proof of Proposition \ref{prop: lower bound} (a):

\begin{itemize}
\item $W=\textup{softmax}(Q K^\top)$. For simplicity of exposition, we replace $W$ with $W'$, defined to be
$$W'_{ij}=\begin{cases}
1& \text{for } i=n \text{ or }j=n\\
\exp(-N|u_i^\top u_j|)& \text{for } i,j \neq n\text{ and }A_{ij}\neq 1\\
\exp(-N)& \text{for } i,j \neq n\text{ and }A_{ij}=1,
\end{cases}$$ As a sanity check, $W$ is simply $W'$ with each row rescaled to sum to one. 

\item $d_v = 1$, and we set $u=1$, so that $Vu=V$.

\item $V\in\mathbb{R}^{n\times 1}$ is set according to $V_i=1$ for $i =1,\dots,n-1$, and $V_n=2$.

\item For $i=1,\dots, n-1$, we set $f_i(\cdot)$ to be:
$$f_i(x)=\begin{cases}
0& \text{if }0\leq x\leq   \frac{2+\exp(-N)k+\frac{1}{4}}{1+\exp(-N)k+\frac{1}{4}}\\
\frac{x-\frac{2+\exp(-N)k+\frac{1}{4}}{1+\exp(-N)k+\frac{1}{4}}}{ \frac{2+\exp(-N)k}{1+\exp(-N)k}-\frac{2+\exp(-N)k+\frac{1}{4}}{1+\exp(-N)k+\frac{1}{4}}} & \text{if } \frac{2+\exp(-N)k+\frac{1}{4}}{1+\exp(-N)k+\frac{1}{4}}<x< \frac{2+\exp(-N)k}{1+\exp(-N)k}\\
1& \text{if } x\geq \frac{2+\exp(-N)k}{1+\exp(-N)k}.
\end{cases}$$ 
Because $\frac{2+\exp(-N)k}{1+\exp(-N)k}>\frac{2+\exp(-N)k+\frac{1}{4}}{1+\exp(-N)k+\frac{1}{4}}$, we have $f_i(\cdot)$ is continuous piece-wise linear for every $i=1,\dots,n-1$. We set $f_n(x)=0$ for every $x$.
\end{itemize}

Now notice that because $f_i(x)\leq 1$ for every $i=1,\dots,n-1$, and $f_n(x)=0$, the optimal value of this instance of Problem $\mathsf{P}$ is at most $k-1$. It will suffice to show that the largest clique $G$ has size $k'$ if and only if the optimal value is $k'$. We prove both directions separately.


\paragraph{If $G$ has a clique of size $k'$, then the optimal value is at least $k'$:} Suppose $G$ has a clique of size $k'$. Without loss of generality, let $I_1\subset [n-1]$ correspond the vertex set of a clique of size $k'$ in $G$. Consider the solution $x^*$, where $x^*_i=1$ for $i\in I_1\cup\{n\}$. Because $|I_1\cup\{n\}|\leq k$, the solution $x^*$ is feasible. We will show that the objective value of $\mathsf{P}$ at $x^*$ is $k'$ (and thus the optimal value is at least $k'$).

Because $W'_{ij}=\exp(-N)$ whenever $i,j\in I_1$, we have \begin{eqnarray*}
    \sum_{i=1}^{n} x^*_i f_i\left(\frac{(w'_i\odot Vu)^\top x^*}{{w'_i}^\top x^*}\right)&=&
    \sum_{i\in I_1} f_i\left(\frac{(w'_{i}\odot Vu)^\top x^*}{{w'_{i}}^\top x^*}\right)\\
    &=&
    \sum_{i\in I_1}f_i\left(\frac{2+k'\exp(-N)}{1+k'\exp(-N)}\right)
    .
\end{eqnarray*}
Since $k'\leq k-1$. It follows that $$\frac{2+k'\exp(-N)}{1+k'\exp(-N)}> \frac{2+k\exp(-N)}{1+k\exp(-N)}.$$ Therefore \begin{eqnarray*}
    \sum_{i=1}^{n} x^*_i f_i\left(\frac{(w'_i\odot Vu)^\top x^*}{{w'_i}^\top x^*}\right)
    =
    \sum_{i\in I_1}f_i\left(\frac{2+k'\exp(-N)}{1+k'\exp(-N)}\right)
    =|I_1|=k'
    .
\end{eqnarray*}

\paragraph{If the largest clique in $G$ has size $k'$, then the optimal value is at most $k'$:} 
By our assumption we have $k'\geq k-M$. Suppose for the sake of contradiction that the optimal value is $k^*>k'$. Let $x^*$ be an optimal solution to $\mathsf{P}$. Notice if $x^*_n=0$, then since $(w'_i\odot Vu)_j=W'_{ij}$ for all $j\in[n-1]$, we have $$\sum_{i=1}^{n} x^*_i f_i\left(\frac{(w'_i\odot Vu)^\top x^*}{{w_i'}^\top x^*}\right)=\sum_{i=1}^{n} x^*_i f_i\left(1\right)=0.$$ Therefore we must have $x^*_n=1$. 

Let $I\subset [n-1]$ be the index set where $x^*_{i}=1$ for $i\in I$. Because $f_i(x)\leq 1$ for every $i=1,\dots,n-1$, we have $|I|\geq k^*$. Let $d(i)$ be the degree of vertex $i$ in the induced subgraph of $G$ with vertex set $\{v_i\}_{i\in I}$. Because $G$ consists of disjoint cliques with size at most $k'$, we have $d(i)\leq k'-1$ for every $i\in C$. 

Fix any $i\in I$ such that $$f_i\left(\frac{(w'_i\odot Vu)^\top x^*}{{w_i'}^\top x^*}\right)>0,$$ or, equivalently, $$\frac{(w'_i\odot Vu)^\top x^*}{{w_i'}^\top x^*}> \frac{2+\exp(-N)+\frac{1}{4}}{1+\exp(-N)+\frac{1}{4}}.$$ 
Because $\exp(-N|u_i^\top u_j|)\geq\exp(-\delta N)$ for every $i,j \neq n\text{ and }A_{ij}\neq 1$, we have \begin{eqnarray*}
    \frac{(w'_i\odot Vu)^\top x^*}{{w'_i}^\top x^*}&=&
    \frac{2+\sum_{j\in I, A_{ij}\neq 1}\exp(-N|u_i^\top u_j|)+(d(i)+1)\exp(-N)}{1+\sum_{j\in I, A_{ij}\neq 1}\exp(-N|u_i^\top u_j|)+(d(i)+1)\exp(-N)}
      \\ &\leq&
    \frac{2+(|I|-d(i)-1)\exp(-\delta N)+(d(i)+1)\exp(-N)}{1+(|I|-d(i)-1)\exp(-\delta N)+(d(i)+1)\exp(-N)}
    \\ &\leq&
    \frac{2+(k^*-d(i)-1)\exp(-\delta N)+(d(i)+1)\exp(-N)}{1+(k^*-d(i)-1)\exp(-\delta N)+(d(i)+1)\exp(-N)}
    \\ &=&
    \frac{2+\exp(-\delta N)k^*-(\exp(-\delta N)-\exp(-N))(d(i)+1)}{1+\exp(-\delta N)k^*-(\exp(-\delta N)-\exp(-N))(d(i)+1)}
    \\ &\leq&
    \frac{2+\exp(-\delta N)k^*-(\exp(-\delta N)-\exp(-N))k'}{1+\exp(-\delta N)k^*-(\exp(-\delta N)-\exp(-N))k'},
\end{eqnarray*}
where the second inequality follows since $|I|\geq k^*$, and the last inequality follows since $\exp(-\delta N)-\exp(-N)>0$ and $d(i)\leq k'-1$.
Therefore $$\frac{2+\exp(-\delta N)k^*-(\exp(-\delta N)-\exp(-N))k'}{1+\exp(-\delta N)k^*-(\exp(-\delta N)-\exp(-N))k'}>\frac{2+\exp(-N)k+\frac{1}{4}}{1+\exp(-N)k+\frac{1}{4}}.$$ So we have $$\exp(-\delta N)k^*-(\exp(-\delta N)-\exp(-N))k'<\exp(- N)k+\frac{1}{4}.$$ Rearrange the above inequality gives  $$\exp(-\delta N)(k^*-k')<\exp(- N)(k-k')+\frac{1}{4}.$$ Because $k^*\geq k'+1$ and $k'\geq k-M$, we get $$\exp(-\delta N)<\exp(- N)M+\frac{1}{4}.$$ However, since $\exp(- N)\leq 1/2M$ and $\exp(-\delta N)\geq 1-1/2M$, we have $$\exp(-\delta N)\geq1-\frac{1}{2M}\geq\frac{5}{6}> \exp(- N)M+\frac{1}{4}.$$ A contradiction. Therefore $k^*\leq k'$.
\end{myproof}

\begin{myproof}[Proof of Proposition \ref{prop: hardness of k}]
Our proof is based on a reduction from Problem \eqref{eqn:problem} to the well-known \textit{Multi-dimensional Knapsack Problem} ($\mathsf{MDKP}$), defined as follows:

\begin{definition}[Multi-dimensional Knapsack Problem]
The Multi-dimensional Knapsack Problem ($\mathsf{MDKP}$) with $c$ items and $d$ dimensions is defined as:
\begin{align}
\max \quad & f_{\mathsf{MDKP}}(x) = p^\top x \tag{$\mathsf{MDKP}$} \label{eqn:MDK} \\
\textup{s.t.} \quad & y_i^\top x \leq t_i \quad \forall i \in [d], \notag \\
& x \in \{0,1\}^c. \notag
\end{align}
Here, $p \in \mathbb{N}^c$, and for all $i \in [d]$, we have $y_i \in \mathbb{Z}_{\geq 0}^c$ and $t_i \in \mathbb{Z}_{\geq 0}$.%
\footnote{We assume $p_{\min} > 0$ without loss of generality. If $p_j = 0$ for some $j \in [c]$, there always exists an optimal solution with $x_j = 0$. Hence, we can safely ignore the $j$-th entry of $p$, $x$, and each $y_i$.}
\end{definition}

We present the reduction in the following proposition: 

\begin{proposition} \label{prop: MDK lower bound} 
Consider Problem \eqref{eqn:problem}, written in the equivalent form $\mathsf{P}$ as given in Observation \ref{obs: simplification}, reproduced below, which has parameters $n$, $k$, and $r_+$, where $r_+$ is the non-negative rank of $W$:
\begin{align}
\tag{$\mathsf{P}$}
\max \quad & f_{\mathsf{P}(I)}(x) = \sum_{i=1}^n x_i f_i\left(\frac{(w_i \odot Vu)^\top x}{w_i^\top x}\right) \\
\textup{s.t.} \quad & x \in \{0,1\}^n, \quad 1 \leq e^\top x \leq k. \notag
\end{align}

Now consider an instance of $\mathsf{MDKP}$ with parameters $c$ and $d$. Suppose there exists an algorithm $\textup{ALG}$ for solving $\mathsf{P}$ with parameters $n = k = c + d + 1$ and $\min\{c, d\} \leq r_+ \leq c + d + 1$, such that for any sufficiently small $\epsilon>0$, the algorithm satisfies 
\[
\textup{ALG}_{\mathsf{P}} \geq \left(1 - \epsilon\right) \textup{OPT}_{\mathsf{P}}
\]
with runtime $T$. Then we can construct an algorithm $\textup{ALG}'$ for solving $\mathsf{MDKP}$ with parameters $c$ and $d$, such that for the same $\epsilon$, it satisfies
\[
\textup{ALG}'_{\mathsf{MDKP}} \geq \left(1 - \epsilon\right) \textup{OPT}_{\mathsf{MDKP}}
\]
with runtime $O(T)$.
\end{proposition}

\begin{myproof}[Proof of Proposition \ref{prop: MDK lower bound}.]

Fix an $\mathsf{MDKP}$ instance. We create an instance of $\mathsf{P}$ as follows:
\begin{itemize}
    \item $n=k=d+c+1$. Here, out of the $n$ total variables, the first $d$ variables will correspond to the $d$ constraints of $\mathsf{MDKP}$, the $c$ variables after that will correspond to the $c$ variables of $\mathsf{MDKP}$, and the last variable will be a dummy variable that must be selected in any optimal solution of $\mathsf{P}$.
    \item For each $i=1,\dots,d$, set $$w_{ij}=\begin{cases}
0& \text{for } j=1,\dots, d\\
y_{i,j-d}& \text{for } j=d+1,\dots, d+c\\
1& \text{for } j=d+c+1.
\end{cases}$$ That is, in block notation, 
\[
w_{i}^\top =
\begin{bNiceArray}{c|c|c}[margin]
0 \cdots 0 & y_i^\top & 1
\end{bNiceArray}.
\]Then we have $$w_i^\top 
x=\sum_{j=1}^cy_{ij}x_{d+j}+x_{d+c+1}.$$ Moreover, since $y_i\in\mathbb{Z}_{\geq 0}^c$ for each $i\in[d]$, the non-negative rank of the sub-matrix of $W$ consisting of its first $d$ rows is at most $\min\{c,d\}$.

\item For each $i=d+1,\dots,d+c+1$, set $w_{i}\in\mathbb{R}_{\geq 0}^{d+c+1}$'s to be any non-zero vectors such that $W$ has the desired non-negative rank $r_+$. This is possible since $\min\{c,d\}\leq r_+\leq c+d+1$.

\item $d_v = 1$, and we set $u=1$, so that $Vu=V$.

\item $V\in\mathbb{R}^{(d+c+1)\times 1}$ is set to be $V_i=0$ for $i=1,\dots,d$, and $V_i=1$ for $i=d+1,\dots,d+c$, and $V_{d+c+1}=2$.

\item For each $i=1,\dots,d$, set $$f_i(x)=\begin{cases}
-\left(\sum_{j=1}^cp_j+1\right)& \text{for } x\leq   \frac{t_i+3}{t_i+2}\\
\dfrac{x - \frac{t_i+2}{t_i+1}}{\frac{t_i+2}{t_i+1} - \frac{t_i+3}{t_i+2}} \cdot \left( \sum\limits_{j=1}^c p_j+1  \right)& \text{for } \frac{t_i+3}{t_i+2}<x< \frac{t_i+2}{t_i+1}\\
0& \text{for } x\geq \frac{t_i+2}{t_i+1}.\\
\end{cases}$$ Then $f_i(x)$ is continuous piecewise-linear.

\item For each $i=d+1,\dots,d+c$, set $f_{i}(x)=p_{i-d}$. Set $f_{d+c+1}(x)=0$.
\end{itemize}

Because $p_{\min}>0$ and $\vec{0}\in\{0,1\}^c$ is a feasible solution to $\mathsf{MDKP}$, we have $\mathrm{OPT}_{\mathsf{MDKP}}=0$ if and only if $\vec{0}\in\{0,1\}^c$ is the only feasible solution to $\mathsf{MDKP}$. Moreover, if $\vec{0}\in\{0,1\}^c$ is not the only feasible solution to $\mathsf{MDKP}$, then $\mathrm{OPT}_{\mathsf{MDKP}}\geq p_{\min}$.

Our proof relies on the following lemma.

\begin{lemma} \label{lemma: construct MDKP solution}
  Let $x\in \{0,1\}^n$ be a feasible solution to $\mathsf{P}$ such that $f_{\mathsf{P}}(x)> 0.$ Then we can construct a feasible solution $z\in \{0,1\}^c$ to $\mathsf{MDKP}$ such that $f_{\mathsf{P}}(x)= f_{\mathsf{MDKP}}(z).$

  Conversely, let $z\in \{0,1\}^c$ be a feasible solution to $\mathsf{MDKP}$ such that $f_{\mathsf{MDKP}}(z)>0.$ Then we can construct a feasible solution $x\in \{0,1\}^n$ to $\mathsf{P}$ such that $f_{\mathsf{P}}(x)= f_{\mathsf{MDKP}}(z).$
\end{lemma}

Lemma \ref{lemma: construct MDKP solution} gives a correspondence between solutions to $\mathsf{P}$ and solutions to $\mathsf{MDKP}$. In particular, as a corollary, Lemma \ref{lemma: construct MDKP solution} implies the relationship between the optimal objective values of $\mathsf{P}$ and $\mathsf{MDKP}$.



\begin{corollary}[Corollary of Lemma \ref{lemma: construct MDKP solution}.] \label{cor: OPT and MDKP}
    $\textup{OPT}_{\mathsf{P}}\leq0$ if and only if $\mathrm{OPT}_{\mathsf{MDKP}}=0$. Moreover, suppose $\mathrm{OPT}_{\mathsf{MDKP}}>0$. Then $\textup{OPT}_{\mathsf{P}}=\mathrm{OPT}_{\mathsf{MDKP}}.$
\end{corollary}
\begin{myproof}[Proof of Corollary \ref{cor: OPT and MDKP}.]
    Suppose $\textup{OPT}_{\mathsf{P}}> 0.$ Let $x^*$ be an optimal solution to $\mathsf{P}$. By Lemma \ref{lemma: construct MDKP solution} there exists a feasible solution $z$ to $\mathsf{MDKP}$ such that $f_{\mathsf{MDKP}}(z)>0$ and $$\textup{OPT}_{\mathsf{P}}=f_{\mathsf{P}}(x^*)=f_{\mathsf{MDKP}}(z)\leq\mathrm{OPT}_{\mathsf{MDKP}}.$$

    Conversely, suppose $\mathrm{OPT}_{\mathsf{MDKP}}>0.$ Let $z^*$ be an optimal solution to $\mathsf{MDKP}$. By Lemma \ref{lemma: construct MDKP solution} there exists a feasible solution $x$ to $\mathsf{MDKP}$ such that $$\mathrm{OPT}_{\mathsf{MDKP}}=f_{\mathsf{MDKP}}(z^*)=f_{\mathsf{P}}(x)\leq\textup{OPT}_{\mathsf{P}}.$$ We also get $$f_{\mathsf{P}}(x)\geq f_{\mathsf{MDKP}}(z^*)>0.$$
\end{myproof}

Before proving Lemma \ref{lemma: construct MDKP solution}, we first show that it is sufficient to prove Lemma \ref{lemma: construct MDKP solution}. Assume we have an algorithm $\textup{ALG}$ for solving $\mathsf{P}$ such that, for any sufficiently small $\epsilon>0$, we have $\textup{ALG}$ satisfies $$\textup{ALG}_{\mathsf{P}}\geq \left(1-\epsilon\right)\textup{OPT}_{\mathsf{P}}.$$ Then our proposed algorithm $\textup{ALG}'$ for solving $\mathsf{MDKP}$ works as follows:
\begin{itemize}
    \item If $\textup{ALG}_{\mathsf{P}}\leq  0,$ $\textup{ALG}'$ outputs $\vec{0}\in\{0,1\}^c$.
    \item If $\textup{ALG}_{\mathsf{P}}> 0,$ let $x$ be the solution to $\mathsf{P}$ given by $\textup{ALG}_{\mathsf{P}}$. Then $\textup{ALG}'$ outputs the solution $z\in \{0,1\}^c$ to $\mathsf{MDKP}$ given by Lemma \ref{lemma: construct MDKP solution}. That is, $z\in \{0,1\}^c$ that satisfies $$f_{\mathsf{P}}(x)=f_{\mathsf{MDKP}}(z).$$
\end{itemize}

\paragraph{Performance Guarantee of $\textup{ALG}'$:} Suppose $\mathrm{OPT}_{\mathsf{MDKP}}=0$. Then by Corollary \ref{cor: OPT and MDKP}, $$\textup{ALG}_{\mathsf{P}}\leq \textup{OPT}_{\mathsf{P}}\leq 0.$$ Therefore $\textup{ALG}'$ correctly outputs $\vec{0}\in\{0,1\}^c$.

On the other hand, suppose $\mathrm{OPT}_{\mathsf{MDKP}}>0$. Then by Corollary \ref{cor: OPT and MDKP}, $$\textup{ALG}'_{\mathsf{MDKP}}=f_{\mathsf{MDKP}}(z)=f_{\mathsf{P}}(x)\geq (1-\epsilon)\textup{OPT}_{\mathsf{P}}=(1-\epsilon)\textup{OPT}_{\mathsf{MDKP}}.$$ 

This proves the desired performance guarantee of $\textup{ALG}'$. To finish the proof, we prove Lemma \ref{lemma: construct MDKP solution} and analyze the runtime of $\textup{ALG}'$. 

\begin{myproof}[Proof of Lemma \ref{lemma: construct MDKP solution}.]
    Recall for each $i=1,\dots,d$, we set $$w_{ij}=\begin{cases}
0& \text{for } j=1,\dots, d\\
y_{i,j-d}& \text{for } j=d+1,\dots, d+c\\
1& \text{for } j=d+c+1.
\end{cases}$$ Also, $(Vu)_i=0$ for $i=1,\dots,d$, and $(Vu)_i=1$ for $i=d+1,\dots,d+c$, and $(Vu)_{d+c+1}=2$. Therefore for $i=1,\dots,d$ we have $$\frac{(w_i\odot Vu)^\top x}{w_i^\top x}=\frac{\sum_{j=1}^cy_{ij}x_{d+j}+2x_{d+c+1}}{\sum_{j=1}^cy_{ij}x_{d+j}+x_{d+c+1}}.$$
Recall for each $i=d+1,\dots,d+c$ we set $f_{i}(x)=p_{i-d}$, and we set $f_{d+c+1}(x)=0$. Therefore we have \begin{eqnarray*}
    f_{\mathsf{P}}(x)&=&\sum_{i=1}^{c+d+1} x_i f_i\left(\frac{(w_i\odot Vu)^\top x}{w_i^\top x}\right)\\
    &=&\sum_{i=1}^{d} x_if_i\left(\frac{\sum_{j=1}^cy_{ij}x_{d+j}+2x_{d+c+1}}{\sum_{j=1}^cy_{ij}x_{d+j}+x_{d+c+1}}\right)+\sum_{i=1}^{c} x_{d+i}p_i.
\end{eqnarray*}

First, let $x\in \{0,1\}^n$ be a feasible solution to $\mathsf{P}(I)$ such that $f_{\mathsf{P}}(x)> 0.$ Let $z\in \{0,1\}^c$ where $z_j=x_{d+j}$ for $j\in[c]$. We claim that $z$ is a feasible solution to $\mathsf{MDKP}$ and $$f_{\mathsf{P}}(x)= f_{\mathsf{MDKP}}(z).$$
Because $\sum_{i=1}^cx_{d+i}p_i<\sum\limits_{j=1}^c p_j+1$, we must have $x_i=1$ for every $i\in[d]$. Moreover, if $x_{d+c+1}=0$, then $$f_i\left(\frac{\sum_{j=1}^cy_{ij}x_{d+j}+2x_{d+c+1}}{\sum_{j=1}^cy_{ij}x_{d+j}+x_{d+c+1}}\right)=f_i(1)<0$$ for every $i\in[d]$. Hence we must have $x_{d+c+1}=1$. Because $y_i\in\mathbb{Z}^c_{\geq 0}$ and $ t_i\in\mathbb{Z}_{\geq 0}$ for all $i\in[d]$, if $\sum_{j=1}^cy_{ij}x_{d+j}>t_i$, we must have $\sum_{j=1}^cy_{ij}x_{d+j}\geq t_i+1$. Then $$f_i\left(\frac{\sum_{j=1}^cy_{ij}x_{d+j}+2x_{d+c+1}}{\sum_{j=1}^cy_{ij}x_{d+j}+x_{d+c+1}}\right)\leq f_i\left(\frac{t_i+3}{t_i+2}\right)=-\left(\sum\limits_{j=1}^c p_j+1\right).$$
Therefore we must have  $\sum_{j=1}^cy_{ij}x_{d+j}\leq t_i$ for all $i\in[d]$. Hence $\sum_{j=1}^cy_{ij}z_{j}\leq t_i$ for all $i\in[d]$, which shows $z$ is a feasible solution to $\mathsf{MDKP}$. Finally, because $$\frac{\sum_{j=1}^cy_{ij}x_{d+j}+2x_{d+c+1}}{\sum_{j=1}^cy_{ij}x_{d+j}+x_{d+c+1}}\geq\frac{t_i+2}{t_i+1}$$ for all $i\in[d]$, we have $$f_i\left(\frac{\sum_{j=1}^cy_{ij}x_{d+j}+2x_{d+c+1}}{\sum_{j=1}^cy_{ij}x_{d+j}+x_{d+c+1}}\right)=0$$ for all $i\in[d]$. Therefore \begin{align*}
    f_{\mathsf{P}}(x)
    &=\sum_{i=1}^{d} x_if_i\left(\frac{\sum_{j=1}^cy_{ij}x_{d+j}+2x_{d+c+1}}{\sum_{j=1}^cy_{ij}x_{d+j}+x_{d+c+1}}\right)+\sum_{i=1}^{c} x_{d+i}p_i\\
    &=\sum_{i=1}^{c} z_{i}p_i\\
    &=f_{\mathsf{MDKP}}(z).
\end{align*}

Conversely, let $z\in \{0,1\}^c$ be a feasible solution to $\mathsf{MDKP}$ such that $f_{\mathsf{MDKP}}(z)>0.$ Let $x$ be a solution to $\mathsf{P}$ where $x_i=1$ for $i=1,\dots,d$ and $i=d+c+1$, and $x_{d+i}=z_i$ for $i=1,\dots,c$. We claim that $x$ is a feasible solution to $\mathsf{MDKP}$ and $$f_{\mathsf{P}}(x)= f_{\mathsf{MDKP}}(z).$$ Because $z$ is a feasible solution to $\mathsf{MDKP}$, we have $\sum_{j=1}^cy_{ij}z_{j}\leq t_i$ for all $i\in[d]$. Therefore $\sum_{j=1}^cy_{ij}x_{d+j}\leq t_i$ for all $i\in[d]$. Then $$\frac{\sum_{j=1}^cy_{ij}x_{d+j}+2x_{d+c+1}}{\sum_{j=1}^cy_{ij}x_{d+j}+x_{d+c+1}}\geq\frac{t_i+2}{t_i+1}$$ for all $i\in[d]$. Therefore we have $$f_i\left(\frac{\sum_{j=1}^cy_{ij}x_{d+j}+2x_{d+c+1}}{\sum_{j=1}^cy_{ij}x_{d+j}+x_{d+c+1}}\right)=0$$ for all $i\in[d]$. Hence \begin{align*}
    f_{\mathsf{P}}(x)
    &=\sum_{i=1}^{d} x_if_i\left(\frac{\sum_{j=1}^cy_{ij}x_{d+j}+2x_{d+c+1}}{\sum_{j=1}^cy_{ij}x_{d+j}+x_{d+c+1}}\right)+\sum_{i=1}^{c} x_{d+i}p_i\\
    &=\sum_{i=1}^{c} z_{i}p_i\\
    &=f_{\mathsf{MDKP}}(z).
\end{align*}
\end{myproof}

Given an $\mathsf{MDKP}$ instance, our construction of the corresponding $\mathsf{P}$ instance takes $O(1)$ runtime. Also, the procedure of the construction of a feasible solution $z\in \{0,1\}^c$ to $\mathsf{MDKP}$ given a feasible solution to $\mathsf{P}$ described in Lemma \ref{lemma: construct MDKP solution} takes $O(1)$ runtime: we simply set $z\in \{0,1\}^c$ where $z_j=x_{d+j}$ for $j\in[c]$. Therefore, if $\textup{ALG}$ has runtime $T$, our $\textup{ALG}'$ has runtime $O(T)$.

\end{myproof}

By Proposition~\ref{prop: MDK lower bound}, any $(1-\epsilon)$-approximation scheme for $\mathsf{P}$ can be directly translated into a $(1-\epsilon)$-approximation scheme for $\mathsf{MDKP}$ with comparable runtime. Therefore, many hardness results for $\mathsf{MDKP}$ carry over directly to $\mathsf{P}$. We cite several such results below:
\begin{proposition}[Theorem 6 in \cite{kulik2010there}]
     If $d = 2$, $\mathsf{MDKP}$ admits no $(1-\epsilon)$-approximation scheme with runtime $f(1/\epsilon)\,c^{O(1)}$ for any function $f$, assuming the $k$-\textsc{Clique} problem is not Fixed-Parameter Tractable.
\end{proposition}

\begin{proposition}[Theorem 5.1 in \cite{jansen2016bounding}]
     If $d = 2$, $\mathsf{MDKP}$ admits no $(1-\epsilon)$-approximation scheme with runtime $f(1/\epsilon)\,c^{O(1)}$ for any function $f$, with runtime $f(1/\epsilon)\,c^{o(1/\epsilon)}$ for any function $f$, assuming Exponential Time Hypothesis holds.
\end{proposition}

\begin{proposition}[Corollary of \cite{doron2024fine}]
    For general $d$, $\mathsf{MDKP}$ admits no $(1-\epsilon)$-approximation scheme with runtime \[k^{o\left(\frac{d}{\epsilon \log^2(d/\epsilon)}\right)} \text{ or } k^{o(\sqrt{d})},\] assuming Gap Exponential Time Hypothesis holds.
\end{proposition}
These results, along with Proposition~\ref{prop: MDK lower bound}, gives our desired lower bound statement in Proposition \ref{prop: hardness of k}.

\end{myproof}

\section{Proofs in Section \ref{Section: Main Result}.}\label{Appendix D}

\begin{myproof}[Proof of Proposition \ref{prop: cluster}]
    Let $i_1\in I_1,\dots, i_\ell\in I_\ell$. Let $\pi:\{(m,m')\mid m,m'\in [\ell],m\leq m'\}\to [\ell(\ell+1)/2]$ be any bijection. 
    Let $A,B\in\mathbb{R}^{n\times \ell(\ell+1)/2}_{\geq 0}$ where $$A_{i,\pi(m,m')}=\begin{cases}
\sqrt{\exp(q_i^\top k_{m'})/\sum_{s=1}^n\exp(q_i^\top k_{s})}& \text{for } i=m\\
0& \text{otherwise,}
\end{cases}$$ and $$B_{j,\pi(m,m')}=\begin{cases}
\sqrt{\exp(q_m^\top k_j)/\sum_{s=1}^n\exp(q_m^\top k_{s})}& \text{for } j=m'\\
0& \text{otherwise.}
\end{cases}$$
Let $W'=AB^\top$. Then \begin{eqnarray*}
    W'_{ij}&=&A_i^\top B_j\\
    &=&\sum_{\pi(m,m')} A_{i,\pi(m,m')}B_{j,\pi(m,m')}\\
    &=&A_{i,\pi(i,j)}B_{j,\pi(i,j)}\\
    &=& \exp(q_i^\top k_j)/\sum_{s=1}^n\exp(q_i^\top k_{s}).
\end{eqnarray*}
Finally, because for every $i,i'\in I_{\ell'}$ we have $\lVert q_i-q_{i'}\rVert_2\leq \delta$ and $\lVert k_i-k_{i'}\rVert_2\leq \delta$, we have exactly the same setup as in the proof of Proposition \ref{prop: phase one}. Therefore the exact same proof gives $1-\gamma\leq W_{ij}/W'_{ij}\leq 1+\gamma$ for all $i,j$, where $\gamma=17\delta \max\{\lVert Q\rVert_{2,\infty},\lVert K\rVert_{2,\infty}\}$.

\end{myproof}

\begin{myproof}[Proof of Corollary \ref{cor: rank partition}.]
     In the proof of Proposition \ref{prop: phase one}, we constructed a partition $\mathcal{I}=\{I_1,\dots,I_\ell\}$ with $\ell=(4\max\{\lVert Q\rVert_{2,\infty},\lVert K\rVert_{2,\infty}\}/\delta)^{2d_{kq}}$ that satisfies for every $i,i'\in I_{\ell'}$ we have $\lVert q_i-q_{i'}\rVert_2\leq \delta$ and $\lVert k_i-k_{i'}\rVert_2\leq \delta$. The result then follows from Proposition \ref{prop: cluster}.
\end{myproof}

\section{Proofs in Section \ref{Section Phase One}.} \label{Appendix E}

\begin{myproof}[Proof of Observation \ref{obs: simplification}]
    Fix $S=\{s_1,\dots,s_m\}$ and let $x\in\{0,1\}^n$ where $x_{s_1}=\cdots=x_{s_m}=1$. We prove that the objective values of Problem \eqref{eqn:problem} and $\mathsf{P}$ are the same. By our construction of $X_S$, we get $XQ\in\mathbb{R}^{n\times d_{kq}}$ where $(XQ)_{s_j}=q_{s_j}$, and similarly for $XK\in\mathbb{R}^{n\times d_{kq}}$ and $XV\in\mathbb{R}^{n\times d_{v}}$. Then for $i\in S$ we get $\text{softmax}((XQ)(XK)^\top)_{i,s_\ell}=\exp(q_{i}^\top k_{s_\ell})/\sum_{s_{\ell'}\in S}\exp(q_{i}^\top k_{s_{\ell'}})$. Therefore, for $i\in S$ we have \begin{eqnarray*}
        \frac{(w_i\odot Vu)^\top x}{w_i^\top x}&=&\frac{\sum_{s_{\ell}\in S}w_{i,s_{\ell}}{Vu}_{s_{\ell}}}{\sum_{s_{\ell}\in S}w_{i,s_{\ell}}}\\
        &=&\frac{\sum_{s_{\ell}\in S}(\exp(q_{i}^\top k_{s_\ell})/\sum_{j\in [n]}\exp(q_{i}^\top k_{j}))(V_{s_\ell}^\top u)}{\sum_{s_{\ell}\in S}(\exp(q_{i}^\top k_{s_\ell})/\sum_{j\in [n]}\exp(q_{i}^\top k_{j}))}\\
        &=&\frac{\sum_{s_{\ell}\in S}(\exp(q_{i}^\top k_{s_\ell})/\sum_{s_{\ell'}\in S}\exp(q_{i}^\top k_{s_{\ell'}}))(V_{s_\ell}^\top u)}{\sum_{s_{\ell}\in S}(\exp(q_{i}^\top k_{s_\ell})/\sum_{s_{\ell'}\in S}\exp(q_{i}^\top k_{s_{\ell'}}))}\\
        &=&\frac{\sum_{s_{\ell}\in S}\text{softmax}((XQ)(XK)^\top)_{i,s_{\ell}}(V_{s_\ell}^\top u)}{\sum_{s_{\ell}\in S}\text{softmax}((XQ)(XK)^\top)_{i,s_{\ell}}}\\
        &=&((\text{softmax}((XQ)(XK)^\top) XV)_i)^\top u\\
        &=& \mathrm{SA}_{Q,K,V}(X_S)_i^\top u,
    \end{eqnarray*}
    where the fifth equality follows since $\sum_{j}\text{softmax}(A)_{i,j}=1$ for every matrix $A$. Because this holds for every $i\in S$, the objective values of Problem \eqref{eqn:problem} and $\mathsf{P}$ are the same.
\end{myproof}

\begin{myproof}[Proof of Proposition \ref{prop: phase one}.]
Fix any $\epsilon>0$. Let $\delta>0$ be a parameter such that \begin{equation}\label{eqn:phase one delta}
\delta \leq \frac{1}{140\max\{\lVert Q\rVert_{2,\infty},\lVert K\rVert_{2,\infty}\}}
\quad\text{ and }\quad
\epsilon \geq 34\delta\max\{\lVert Q\rVert_{2,\infty},\lVert K\rVert_{2,\infty}\}(Vu)_{\max}+\delta.
\end{equation}

First we partition the rows of $Q$ and $K$. Because we can cover a ball with radius $\lVert Q\rVert_{2,\infty}$ in $\mathbb{R}^{d_{kq}}$ using $\lceil4\lVert Q\rVert_{2,\infty}/\delta\rceil^{d_{kq}}$ number of balls with radius $\delta/2$ (see e.g. \cite{verger2005covering,dumer2007covering,shalev2014understanding}), we can create a partition of the index set $[n]$ such that the size of the partition is $\lceil 4\lVert Q\rVert_{2,\infty}/\delta\rceil^{d_{kq}}$, and $\lVert q_i-q_{i'}\rVert_2\leq\delta$ for every $i,i'$ in the same partition. \footnote{We can create such a partition via the packing-covering duality in the following way (see e.g. Theorem 14.1, Theorem 14.2, and Example 14.1 of \cite{wu2017lec14}). First we create a maximal packing of the ball $B(\vec{0},\lVert Q\rVert_{2,\infty})\subset \mathbb{R}^{d_{kq}}$ using balls with radius $\delta/4$ greedily, where we greedily choose points $x_1,x_2,\dots$ such that the balls $B(x_i,\delta/4)$ are disjoint. We stop when no more such points can be added in $B(\vec{0},\lVert Q\rVert_{2,\infty})\subset \mathbb{R}^{d_{kq}}$. Then, by the packing-covering duality, the balls $B(x_i,\delta/2)$ cover $B(\vec{0},\lVert Q\rVert_{2,\infty})$. This is because any uncovered point can be added to the packing we created before, which contradicts the maximality of the packing.} Similarly, we can create a partition of the index set $[n]$ such that the size of the partition is $\lceil4\lVert K\rVert_{2,\infty}/\delta\rceil ^{d_{kq}}$, and $\lVert k_j-k_{j'}\rVert_2\leq\delta$ for every $j,j'$ in the same partition. Let $\mathcal{I}=\{I_1,\dots,I_\ell\}$ be the product partition of the above two partitions. Then $\ell=\lceil4\max\{\lVert Q\rVert_{2,\infty},\lVert K\rVert_{2,\infty}\}/\delta\rceil^{2d_{kq}}$, and for every $\ell'\in[\ell]$ and $i,i'\in I_{\ell'}$, we have $\lVert q_i-q_{i'}\rVert_2\leq\delta$ and $\lVert k_i-k_{i'}\rVert_2\leq\delta$. Therefore, for any $i,i'\in I_{\ell'}$ and $j,j'\in I_{\ell''}$, we have \begin{eqnarray*}
    \lvert q_i^\top k_j-q_{i'}^\top k_{j'}\rvert&=&\lvert q_i^\top(k_j-k_{j'})+k_{j'}^\top (q_i-q_{i'})\rvert\\
    &\leq& \lvert q_i^\top(k_j-k_{j'})\rvert+\lvert k_{j'}^\top (q_i-q_{i'})\rvert\\
    &\leq& \lVert q_i\rVert_2\lVert k_j-k_{j'}\rVert_2+\lVert k_{j'}\rVert_2\lVert q_i-q_{i'}\rVert_2 \\
    &\leq& 2\delta \max\{\lVert Q\rVert_{2,\infty},\lVert K\rVert_{2,\infty}\}.
\end{eqnarray*}
Then, because $2\delta \max\{\lVert Q\rVert_{2,\infty},\lVert K\rVert_{2,\infty}\}\leq1$, we have \begin{eqnarray*}
    \left\lvert\frac{\exp(q_i^\top k_j)}{\exp(q_{i'}^\top k_{j'})}-1\right\rvert&\leq&\left\lvert\exp(\lvert q_i^\top k_j-q_{i'}^\top k_{j'}\rvert)-1\right\rvert\\
    &\leq&\lvert\exp(2\delta \max\{\lVert Q\rVert_{2,\infty},\lVert K\rVert_{2,\infty}\})-1\rvert\\
    &\leq& 4\delta \max\{\lVert Q\rVert_{2,\infty},\lVert K\rVert_{2,\infty}\},
\end{eqnarray*}
where the last inequality follows by $\exp(x)\leq1+2x$ for $0\leq x\leq1$. Therefore, because $4\delta \max\{\lVert Q\rVert_{2,\infty},\lVert K\rVert_{2,\infty}\}\leq 1/35$, we have \begin{eqnarray*}
    \left\lvert\frac{w_{ij}}{w_{i'j'}}-1\right\rvert&=&\left\lvert\frac{\exp(q_i^\top k_j)/\sum_{m=1}^n\exp(q_i^\top k_m)}{\exp(q_{i'}^\top k_{j'})/\sum_{m'=1}^n\exp(q_{i'}^\top k_{m'})}-1\right\rvert\\
    &\leq& \left\lvert(1+4\delta \max\{\lVert Q\rVert_{2,\infty},\lVert K\rVert_{2,\infty}\})^2/(1-4\delta \max\{\lVert Q\rVert_{2,\infty},\lVert K\rVert_{2,\infty}\})^2-1\right\rvert\\
    &\leq&\frac{4\cdot35^2}{34^2}\cdot 4\delta \max\{\lVert Q\rVert_{2,\infty},\lVert K\rVert_{2,\infty}\}\\
    &\leq&17\delta \max\{\lVert Q\rVert_{2,\infty},\lVert K\rVert_{2,\infty}\},
\end{eqnarray*} 
where the first inequality follows since $1-4\delta \max\{\lVert Q\rVert_{2,\infty},\lVert K\rVert_{2,\infty}\}\leq\exp(q_i^\top k_j)/\exp(q_{i'}^\top k_{j'})\leq1+4\delta \max\{\lVert Q\rVert_{2,\infty},\lVert K\rVert_{2,\infty}\}$, and the second inequality follows since $(1+x)^2/(1-x)^2\leq 1+(4a^2/(a-1)^2)x$ for every $a>1$ and $0\leq x\leq 1/a$.

Now we construct our desired index set $I$. Let $\{S_1,\dots, S_\tau\}$ be a partition of the index set $[n]$ such that $f_i=f_j$ for every $\tau'\in[\tau]$ and $i,j\in S_{\tau'}$. We first construct an index set $J_{\tau'}\subset S_{\tau'}$ for each $\tau'\in[\tau]$, and then combine them to obtain $I$. 

Fix any $\tau'\in\tau$. For each index set $I_{\ell'}$, we only choose $k$ indices out of it to include in $J_{\tau'}$, namely the $k$ indices that are approximately the $k$ highest indices in $\{(Vu)_{i}\}_{i\in S_{\tau'}\cap I_{\ell'}}$.\footnote{For simplicity we assume $\lvert  S_{\tau'}\cap I_{\ell'}\rvert\geq k$. Otherwise we simply choose all indices in $S_{\tau'}\cap I_{\ell'}$.} Specifically, for each $\ell'\in[\ell]$, we run the given $\delta$-Approximate $k$-Nearest Neighbor oracle with given set of points $\bigcup_{i\in S_{\tau'}\cap I_{\ell'}}\{V_i\}\subset\mathbb{R}^{d_v}$, and query $u$, numbers $k$ and $\delta$ as inputs. We let $J_{\tau'}$ be the collection of all output indices for each $\ell'\in[\ell]$. Then because $\ell=\lceil4\max\{\lVert Q\rVert_{2,\infty},\lVert K\rVert_{2,\infty}\}/\delta\rceil^{2d_{kq}}$, we have \begin{equation*}
    |J_{\tau'}|=k\lceil4\max\{\lVert Q\rVert_{2,\infty},\lVert K\rVert_{2,\infty}\}/\delta\rceil^{2d_{kq}},
\end{equation*} and the expect amortized runtime of constructing each $J_{\tau'}$ is \begin{equation*}
    \lceil4\max\{\lVert Q\rVert_{2,\infty},\lVert K\rVert_{2,\infty}\}/\delta\rceil^{2d_{kq}}\cdot k\text{-ANN}(|S_{\tau'}|,d_v,k,\delta).
\end{equation*}
Let $I=\cup_{\tau'\in[\tau]}J_{\tau'}$. Then we have \begin{equation*}
    |I|=\tau k\lceil4\max\{\lVert Q\rVert_{2,\infty},\lVert K\rVert_{2,\infty}\}/\delta\rceil^{2d_{kq}},
\end{equation*} and the expect amortized runtime of constructing $I$ is \begin{align*}
    &\quad\text{ } \sum_{\tau'=1}^\tau\lceil4\max\{\lVert Q\rVert_{2,\infty},\lVert K\rVert_{2,\infty}\}/\delta\rceil^{2d_{kq}}\cdot k\text{-ANN}(|S_{\tau'}|,d_v,k,\delta)\\&\leq \lceil4\max\{\lVert Q\rVert_{2,\infty},\lVert K\rVert_{2,\infty}\}/\delta\rceil^{2d_{kq}}\tau\cdot k\text{-ANN}\left(\frac{n}{\tau},d_v,k,\delta\right).
\end{align*} The inequality follows since $\sum_{\tau'=1}^\tau |S_{\tau'}|=n$ and $k\textup{-ANN}(n,d,k,\epsilon)$ is concave in $n$.

Finally, we show that $\textup{OPT}_{\mathsf{P}(I)}\geq (1-g(\epsilon))\textup{OPT}_{\mathsf{P}}-kh(\epsilon)$ with our choice of $\delta$. Let $i^*_1,\dots,i_k^*$ be the non-zero coordinates of an optimal solution $x^*$ to the original problem (for simplicity we assume $x^*$ has $k$ non-zero entries, and other cases can be handled similarly). For each $m=1,\dots,k$, let $i_m$ be the index such that $i_m$ and $i^*_m$ are in the same $S_{\tau'}\cap I_{\ell'}$ and $(Vu)_{i_m}\geq (Vu)_{i^*_m}-\delta$. Then $f_{i_m}=f_{i^*_m}$. Let $x\in\{0,1\}^n$ such that $x_{i_m}=1$, then $x$ is feasible to $\mathsf{P}(I)$. We prove that the objective value of $\mathsf{P}(I)$ at $x$ is at least $(1-g(\epsilon))\textup{OPT}_{\mathsf{P}}-kh(\epsilon)$. Indeed, because $\epsilon\geq 34\delta\max\{\lVert Q\rVert_{2,\infty},\lVert K\rVert_{2,\infty}\}(Vu)_{\max}+\delta$, then for each $m=1,\dots,k$ we have

\begin{eqnarray*}
    f_{i_m}\left(\frac{(w_{i_m}\odot Vu)^\top x}{w_{i_m}^\top x}\right)&=&f_{i_m}\left(\frac{\sum_{j=1}^k (w_{i_m})_{i_j}(Vu)_{i_j}}{\sum_{j=1}^k (w_{i_m})_{i_j}}\right)\\
    &\geq&f_{i_m}\left(\frac{\sum_{j=1}^k (w_{i_m})_{i_j}(Vu)_{i^*_j}}{\sum_{j=1}^k (w_{i_m})_{i_j}}-\delta\right)\\
    &\geq& f_{i^*_m}\left(\frac{(1-\epsilon)\sum_{j=1}^k (w_{i^*_m})_{i^*_j}(Vu)_{i^*_j}}{(1+\epsilon)\sum_{j=1}^k (w_{i^*_m})_{i^*_j}}-\delta\right)\\
    &\geq& f_{i^*_m}\left(\frac{\sum_{j=1}^k (w_{i^*_m})_{i^*_j}(Vu)_{i^*_j}}{\sum_{j=1}^k (w_{i^*_m})_{i^*_j}}-34\delta \max\{\lVert Q\rVert_{2,\infty},\lVert K\rVert_{2,\infty}\}(Vu)_{\max}-\delta\right)\\
    &\geq& f_{i^*_m}\left(\frac{\sum_{j=1}^k (w_{i^*_m})_{i^*_j}(Vu)_{i^*_j}}{\sum_{j=1}^k (w_{i^*_m})_{i^*_j}}-\epsilon\right)\\
    &\geq &(1-g(\epsilon))f_{i^*_m}\left(\frac{\sum_{j=1}^k (w_{i^*_m})_{i^*_j}(Vu)_{i^*_j}}{\sum_{j=1}^k (w_{i^*_m})_{i^*_j}}\right)-h(\epsilon),
\end{eqnarray*}
where the first inequality follows since $(Vu)_{i_j}\geq (Vu)_{i^*_j}+\delta$, the second inequality follows since $1-\epsilon\leq(w_{i^*_m})_{i^*_j}/(w_{i_m})_{i_j}\leq1+\epsilon$ for every $i_m,i^*_m$ that are in the same partition in $\mathcal{I}$ and $i_j,i^*_j$ that are in the same partition in $\mathcal{I}$, and the third inequality follows since $(1-\epsilon)/(1+\epsilon)\geq 1-2\epsilon$ and $\sum_{j=1}^k (w_{i^*_m})_{i^*_j}(Vu)_{i^*_j}/\sum_{j=1}^k (w_{i^*_m})_{i^*_j}\leq (Vu)_{\max}$. Then \begin{eqnarray*}
   \textup{OPT}_{\mathsf{P}(I)}&\geq& \sum_{m=1}^k x_{i_m}f_{i_m}\left(\frac{(w_{i_m}\odot Vu)^\top x}{w_{i_m}^\top x}\right)\\
    &\geq & \sum_{m=1}^k x^*_{i^*_m}\left(g(\epsilon)f_{i^*_m}\left(\frac{\sum_{j=1}^k (w_{i^*_m})_{i^*_j}(Vu)_{i^*_j}}{\sum_{j=1}^k (w_{i^*_m})_{i^*_j}}\right)-h(\epsilon)\right)\\
    &\geq & (1-g(\epsilon))\textup{OPT}_{\mathsf{P}}-kh(\epsilon)
\end{eqnarray*} as desired. By our choice of $\delta$, we have \begin{equation*}
    |I|=k\left\lceil\frac{140(\max\{\lVert Q\rVert_{2,\infty},\lVert K\rVert_{2,\infty}\})^2}{(Vu)_{\max}\cdot\epsilon}\right\rceil^{2d_{kq}},
\end{equation*} and the expected amortized runtime of finding $I$ is \begin{equation*}
  \left\lceil\frac{140(\max\{\lVert Q\rVert_{2,\infty},\lVert K\rVert_{2,\infty}\})^2(Vu)_{\max}}{\epsilon}\right\rceil^{2d_{kq}} \tau \cdot k\textup{-ANN}\left(\frac{n
  }{\tau},d_v, k,\frac{\epsilon}{35\max\{\lVert Q\rVert_{2,\infty},\lVert K\rVert_{2,\infty}\}(Vu)_{\max}}\right).
\end{equation*}

\end{myproof}

Below we give the pseudo-code of our phase one algorithm in Proposition~\ref{prop: phase one}.

\begin{algorithm}[H]
\caption{Phase One: Preprocess}
\label{alg: phase one preprocess}
\DontPrintSemicolon
\SetAlgoLined

\BlankLine
\tcp{\textbf{Preprocess}}
\KwIn{Number of items $n$, maximum number of recommended items $k$, key matrix $K\in\mathbb{R}^{n\times d_{kq}}$, query matrix $Q\in\mathbb{R}^{n\times d_{kq}}$, value matrix $V\in\mathbb{R}^{n\times d_{v}}$, reward functions $\{f_i\}_{i=1}^n$, parameter $\epsilon > 0$, and $\epsilon$-Approximate $k$-Nearest Neighbor oracle}
\KwOut{Parameter $\delta>0$, partitions $\mathcal{I} = \{I_1, \ldots, I_\ell\}$ and $\mathcal{S} = \{S_1, \ldots, S_\tau\}$, and preprocessed $\epsilon$-Approximate $k$-Nearest Neighbor oracle with set of points $\{V_i : i \in I_{\ell'}\cap S_{\tau'}\}$ for each $\tau'\in[n]$ and $\ell'\in[\ell]$}
\BlankLine

Set $\delta > 0$ according to Eq.~\eqref{eqn:phase one delta}\;

Let $R \leftarrow \max\{\|Q\|_{2,\infty}, \|K\|_{2,\infty}\}$\;

Form a maximal $(\delta/2)$-separated set $C\subset B(0,R)$\; 

Set balls $\{B_i\}_{i=1}^m \leftarrow \{B(c,\delta/2): c\in C\}$ with $m=|C|$\;

Create partition $\mathcal{I} = \{I_1, \ldots, I_\ell\}$ of $[n]$: each $I_j$ contains indices whose $(K_i, Q_i)$ fall in the same pair of balls\;

Let $\tau
$ be the number of distinct functions among $f_1,\dots,f_n$\;

Create partition $\mathcal{S} = \{S_1, \ldots, S_\tau\}$ of $[n]$ where $f_i = f_j$ for all $i,j \in S_{\tau'}$\;

\BlankLine
\For{$\tau' = 1$ \KwTo $\tau$}{
    \For{$\ell' = 1$ \KwTo $\ell$}{
        \If{$I_{\ell'}\cap S_{\tau'} \neq \emptyset$}{
            Preprocess the $\epsilon$-Approximate $k$-Nearest Neighbor oracle with set of points $\{V_i : i \in I_{\ell'}\cap S_{\tau'}\}$\;
        }
    }
}
\Return{Parameter $\delta>0$, partitions $\mathcal{I} = \{I_1, \ldots, I_\ell\}$ and $\mathcal{S} = \{S_1, \ldots, S_\tau\}$, and preprocessed $\epsilon$-Approximate $k$-Nearest Neighbor oracle with set of points $\{V_i : i \in I_{\ell'}\cap S_{\tau'}\}$ for each $\tau'\in[n]$ and $\ell'\in[\ell]$}

\end{algorithm}

\begin{algorithm}[H]
\caption{Phase One: Query}
\label{alg: phase one query}
\DontPrintSemicolon
\SetAlgoLined

\BlankLine
\tcp{\textbf{Query}}

\KwIn{User vector $u \in \mathbb{R}^{d_v}$, maximum number of recommended items $k$, outputs of Algorithm \ref{alg: phase one preprocess}: parameter $\delta>0$, partitions $\mathcal{I} = \{I_1, \ldots, I_\ell\}$ and $\mathcal{S} = \{S_1, \ldots, S_\tau\}$, and preprocessed $\epsilon$-Approximate $k$-Nearest Neighbor oracle with set of points $\{V_i : i \in I_{\ell'}\cap S_{\tau'}\}$ for each $\tau'\in[n]$ and $\ell'\in[\ell]$}
\KwOut{Index set $I$ of candidate items}

\BlankLine
Initialize $I \leftarrow \emptyset$\;

\For{$\tau' = 1$ \KwTo $\tau$}{
    \For{$\ell' = 1$ \KwTo $\ell$}{
        \If{$I_{\ell'}\cap S_{\tau'} \neq \emptyset$}{
            $J \leftarrow$ Query the $\epsilon$-Approximate $k$-Nearest Neighbor oracle with set of points $\{V_i : i \in  I_{\ell'}\cap S_{\tau'}\}$, query $u$, and numbers $k$ and $\delta$ as inputs\;
            $I \leftarrow I \cup J$\;
        }
    }
}

\Return{$I$}

\end{algorithm}

\section{Proof of Proposition \ref{prop: phase two}.} \label{Appendix F}

\subsection{Step 1: Low Non-negative Rank Approximation}
First we bound the loss incurred by replacing $W$ with $W'$:

\begin{lemma}\label{lemma: replace W'}
    Let Problem $\mathsf{P}'(I)$ be defined as
    \begin{align*}\label{eqn:P'(I)}
    \tag{$\mathsf{P}'(I)$}
    \max \quad & f_{\mathsf{P}'(I)} = \sum_{i=1}^m x_i f_i\left(\frac{(w'_i \odot Vu)^\top x}{w'_i{^\top} x}\right) \\
    \textup{s.t.} \quad & x \in \{0,1\}^m,\ 1 \leq e^\top x \leq k.
\end{align*} 
 Let $x$ be a feasible solution to $\mathsf{P}'(I)$ (and hence also a feasible solution to $\mathsf{P}(I)$), and suppose $x$ satisfies $$f_{\mathsf{P}'(I)}(x)\geq (1-\alpha)\textup{OPT}_{\mathsf{P}'(I)}-\beta.$$ 
 Then we have
    \begin{align*}
f_{\mathsf{P}(I)}(x)
    &\geq (1-\alpha)(1-g(2\gamma(Vu)_{\max}))^2\textup{OPT}_{\mathsf{P}(I)}\\&\quad\text{ }-kh(2\gamma(Vu)_{\max})(1+(1-\alpha)(1-g(2\gamma(Vu)_{\max}))-\beta(1-g(2\gamma(Vu)_{\max})).
\end{align*}
\end{lemma}
\begin{myproof}[Proof of Lemma \ref{lemma: replace W'}]
    Let $x$ be a feasible solution to $\mathsf{P}'(I)$. Because $1-\gamma\leq W_{ij}/W'_{ij}\leq1+\gamma$, we have \begin{equation*}
       \frac{({w'_i}\odot Vu)^\top x}{{w'_i}^\top x}\geq\frac{(1-\gamma)}{(1+\gamma)}\frac{({w_i}\odot Vu)^\top x}{{w_i}^\top x}\geq(1-2\gamma)\frac{({w_i}\odot Vu)^\top x}{{w_i}^\top x}.
    \end{equation*}
    Therefore, for any feasible solution $x$, we have
    \begin{eqnarray}
        f_{\mathsf{P}'(I)}(x)&=& \sum_{i=1}^m x_i f_i\left(\frac{({w'_i}\odot Vu)^\top x}{{w'_i}^\top x}\right) \notag\\
        &\geq&\sum_{i=1}^m x_i f_i\left((1-2\gamma)\frac{({w_i}\odot Vu)^\top x}{w_i^\top x}\right) \notag\\
        &\geq&\sum_{i=1}^m x_i f_i\left(\frac{({w_i}\odot Vu)^\top x}{w_i^\top x}-2\gamma(Vu)_{\max}\right)\notag\\
        &\geq& (1-g(2\gamma(Vu)_{\max}))f_{\mathsf{P}(I)}(x)-kh(2\gamma(Vu)_{\max}), \label{eqn:lemma: replace W' 1}
    \end{eqnarray} where the third inequality follows since $({w_i}\odot Vu)^\top x^*/w_i^\top x^*\leq(Vu)_{\max}$. 
    
    Similarly, we also have 
    \begin{equation*}
       \frac{({w_i}\odot Vu)^\top x}{{w_i}^\top x}\geq(1-2\gamma)\frac{({w'_i}\odot Vu)^\top x}{{w'_i}^\top x},
    \end{equation*} 
    which gives 
    \begin{equation} \label{eqn:lemma: replace W' 2}
    f_{\mathsf{P}(I)}(x)\geq (1-g(2\gamma(Vu)_{\max}))f_{\mathsf{P}'(I)}(x)-kh(2\gamma(Vu)_{\max}).
    \end{equation}

       Now Let $x^*_{\mathsf{P}(I)}$ be an optimal solution to $\mathsf{P}(I)$. Then by Eq.~\eqref{eqn:lemma: replace W' 1}, we have \begin{align*}
    \mathrm{OPT}_{\mathsf{P}'(I)}&\geq f_{\mathsf{P}'(I)}(x^*_{\mathsf{P}'(I)})\\
    &\geq (1-g(2\gamma(Vu)_{\max}))f_{\mathsf{P}(I)}(x^*_{\mathsf{P}(I)})-kh(2\gamma(Vu)_{\max})\\
    &= (1-g(2\gamma(Vu)_{\max}))\mathrm{OPT}_{\mathsf{P}(I)}-kh(2\gamma(Vu)_{\max}).
\end{align*}
Finally, applying Eq.~\eqref{eqn:lemma: replace W' 2}, we conclude that
\begin{align*}
f_{\mathsf{P}(I)}(x)&\geq (1-g(2\gamma(Vu)_{\max}))f_{\mathsf{P}'(I)}(x)-kh(2\gamma(Vu)_{\max})\\
    &\geq (1-\alpha)(1-g(2\gamma(Vu)_{\max}))\textup{OPT}_{\mathsf{P}'(I)}-kh(2\gamma(Vu)_{\max})-\beta(1-g(2\gamma(Vu)_{\max}))\\
    &\geq (1-\alpha)(1-g(2\gamma(Vu)_{\max}))^2\textup{OPT}_{\mathsf{P}(I)}\\&\quad\text{ }-kh(2\gamma(Vu)_{\max})(1+(1-\alpha)(1-g(2\gamma(Vu)_{\max}))-\beta(1-g(2\gamma(Vu)_{\max})).
\end{align*}
\end{myproof}
Lemma \ref{lemma: replace W'} shows that, in order to approximately solve $\mathsf{P}(I)$, it is enough to approximately solve $\mathsf{P}'(I)$. 

Let $W'=AB^\top$ be the known non-negative factorization, where $A,B\in \mathbb{R}_{\geq 0}^{m\times r_+}$. Let $a_i^\top\in\mathbb{R}_{\geq 0}^{r_+}$ be the $i$-th row of $A$ and $b_j\in\mathbb{R}_{\geq 0}^{m}$ be the $j$-th column of $B$. Then $w'_i=\sum_{j=1}^{r_+} a_{ij}b_j$. Let \begin{equation}\label{eqn:b_j}
 d_j=b_j\odot (Vu).
\end{equation} Then $\mathsf{P}'(I)$ can be rewritten as:
\begin{align*}
    \tag{$\mathsf{P}'(I)$}
    \max \quad & f_{\mathsf{P}'(I)} = \sum_{i=1}^m x_i f_i\left(\frac{\sum_{j=1}^{r_+} a_{ij}d_j^\top x}{{w'_i}^\top x}\right) \\
    \textup{s.t.} \quad & x \in \{0,1\}^m,\ 1 \leq e^\top x \leq k.
\end{align*}
        
\subsection{Step 2: Enumeration of Partial Solutions}

Our algorithm enumerates a set of partial solutions (where a ``partial solution'' fixes the values of a subset of variables), and then for each partial solution, solves the remaining problem near-optimally. In this step we bound the total number of partial solutions, and in the next steps we show that for each partial solution, the remaining problem can be solved sufficiently fast. 

Let \begin{equation} \label{eqn: lambda}
    \lambda=\lceil(2r_++2)(Vu)_{\max}/\epsilon\rceil.
\end{equation}
Each partial solution that our algorithm considers is specified by a tuple $(X_1,\dots,X_{r_+})$, where each $X_j\subset [m]$ is an index set such that $1\leq |X_1|=\cdots=|X_{r_+}|\leq \lambda$. For each $j \in [r_+]$, let 
\begin{equation} \label{eqn:Xhatj}
\hat{X}_j=\{i\in[m]\setminus X_j\mid d_{ji}> \min_{i'\in X_j}\{d_{ji'}\}\}.
\end{equation}
In words, $\hat{X}_j$ consists of the indices outside of $X_j$ whose coefficients in $d_j$ are strictly greater than the minimum coefficient in $d_j$ across indices in $X_j$. 

We say a tuple $(X_1,\dots,X_{r_+})$, with corresponding index sets $\hat{X}_1,\ldots,\hat{X}_{r_+}$ defined according to \eqref{eqn:Xhatj}, is {\bf valid} if $\lvert\cup_{j}X_j\rvert\leq k$ and  $(\cup_{j}X_j)\bigcap(\cup_{j}\hat{X}_j)=\emptyset$. 
Then every feasible solution $z$ to $\mathsf{P}'(I)$ \textbf{corresponds} to a valid tuple $(X_1,\dots,X_{r_+})$ in the following way: 
Let $Z =\{i\in[m]\mid z_i=1\}$. 
For each $j\in [r_+]$, we define $X_j\subset Z$ to be the set of indices $i \in Z$ such that $d_{ji}$ is among the $\min\{\lambda,|Z|\}$ highest values in $Z$. That is, 
let $\pi_j: [|Z|] \to Z$ be a sorting of $Z$ according to $d_j$ such that $d_{j,\pi_j(1)}\geq \cdots\geq d_{j,\pi_j(|Z|)}$. Then $X_j=\{\pi_j(1),\dots,\pi_j(\min\{\lambda,|Z|\})\}$. Notice that $|Z|\leq k$, so $\lvert\cup_{j}X_j\rvert\leq k$. Also, we claim that $Z\cap \hat{X}_j=\emptyset$ for each $j\in[r_+]$. Supposing otherwise that $i\in Z\cap \hat{X}_j$, then by construction $d_{ji}>d_{j,\pi_j(\min\{\lambda,|Z|\})}$. Because $i\in Z$, we have that $i$ is in the image of $\pi_j$, so there exists $i'$ such that $\pi_j(i')=i$. Therefore $\pi_j(i')\geq \pi_j(\min\{\lambda,|Z|\})$, which shows $i\in X_j$. This contradicts $X_j\cap \hat{X}_j=\emptyset$. Therefore $Z\cap \hat{X}_j=\emptyset$ for each $j\in[r_+]$. Hence we have $(\cup_{j}X_j)\bigcap(\cup_{j}\hat{X}_j)=\emptyset$. Therefore $(X_1,\dots,X_{r_+})$ is indeed a valid tuple.

The notion of correspondence to valid tuples forms a partition of the set of feasible solutions to $\mathsf{P}'(I)$, so it suffices to solve $\mathsf{P}'(I)$ separately for each subset of this partition. This is formally stated in the following result:
\begin{lemma} \label{lemma: enumerate valid tuples}
    Suppose we are given an oracle $\textup{ALG}'$ that takes $\mathsf{P}'(I)$, any valid tuple $(X_1,\dots,X_{r_+})$, and any $\delta>0$ as inputs, and outputs a solution $x'_{(X_1,\dots,X_{r+})}$ of $\mathsf{P}'(I)$ that satisfies\begin{enumerate}
        \item $x'_{(X_1,\dots,X_{r+})}$ corresponds to $(X_1,\dots,X_{r_+})$, and
        \item for any $x_{(X_1,\dots,X_{r+})}$ that corresponds to $(X_1,\dots,X_{r_+})$, we have $$f_{\mathsf{P}'(I)}(x'_{(X_1,\dots,X_{r+})})\geq (1-g'(\delta))f_{\mathsf{P}'(I)}(x_{(X_1,\dots,X_{r+})})-h'(\delta),$$ where $0\leq g'(\delta)\leq 1$ and $h'(\delta)\geq 0$,
    \end{enumerate} with runtime $T(\delta)$. Then there exists an algorithm $\textup{ALG}$ that takes $\mathsf{P}'(I)$ and any $\delta>0$ as inputs, and outputs a solution of $\mathsf{P}'(I)$ that satisfies 
    $$\textup{ALG}_{\mathsf{P}'(I)}\geq (1-g'(\delta))\textup{OPT}_{\mathsf{P}'(I)}-h'(\delta)$$ 
    with runtime 
    $$r_+ m \log_2 m +\lambda m^\lambda+ \lambda r_+^2m^{\lambda r_+}+m^{\lambda r_+}T(\delta).$$
\end{lemma}

\begin{myproof}[Proof of Lemma \ref{lemma: enumerate valid tuples}.]
We will construct $\mathrm{ALG}$ explicitly. Now because every feasible solution $z$ to $\mathsf{P}'(I)$ corresponds to exactly one valid tuple, we can solve $\mathsf{P}'(I)$ by enumerating all valid tuples, and applying $\mathrm{ALG}'$ to each valid tuple. It turns out that pre-sorting the vectors $d_j$ allows for more-efficient enumeration. Let $\textup{ALG}$ take the following steps: 
\begin{enumerate}
\item {\bf Sort $d_j$ for each $j\in[r_+]$.} This takes runtime $r_+ m \log_2 m $.
    \item {\bf Enumerate all valid tuples with $|X_1|<\lambda$, and record the unique corresponding feasible solutions.} 
    We will show momentarily that when $|X_1|<\lambda$, there is a unique corresponding feasible solution, and as a result this step takes runtime $\lambda m^\lambda$.
    \item {\bf Enumerate all valid tuples with $|X_1|=\lambda$, and record the solution output by $\textup{ALG}'_{\mathsf{P}'(I)}$ for each such valid tuple.} We will show that it takes runtime $\lambda r_+^2m^{\lambda r_+}$ to enumerate all such valid tuples, and then because there are at most $m^{\lambda r_+}$ such valid tuples, the total runtime of this step is $\lambda r_+^2m^{\lambda r_+}+m^{\lambda r_+}T(\delta)$.
    \item {\bf Output a solution that is recorded with the highest objective value in} $\textsf{P}'(I)$. 
\end{enumerate}  Therefore the total runtime of $\textup{ALG}$ is $$r_+ m \log_2 m +\lambda m^\lambda+ \lambda r_+^2m^{\lambda r_+}+m^{\lambda r_+}T(\delta).$$
Finally, let $x^*_{(X^*_1,\dots,X^*_{r+})}$ be an optimal solution of $\mathsf{P}'(I)$ where $(X^*_1,\dots,X^*_{r+})$ is the valid tuple that it corresponds to. Let $x'_{(X^*_1,\dots,X^*_{r+})}$ be the solution that $\textup{ALG}'_{\mathsf{P}'(I)}$ outputs with input $\mathsf{P}'(I)$, valid tuple $(X^*_1,\dots,X^*_{r_+})$, and  $\delta>0$. Then \begin{align*}
\textup{ALG}_{\mathsf{P}'(I)}&\geq f_{\mathsf{P}'(I)}(x'_{(X^*_1,\dots,X^*_{r+})})\\&\geq (1-g'(\delta))f_{\mathsf{P}'(I)}(x^*_{(X^*_1,\dots,X^*_{r+})})-h'(\delta)
\\&= (1-g'(\delta))\textup{OPT}_{\mathsf{P}'(I)}-h'(\delta).
\end{align*}

It remains to analyze Steps 2 and 3 of $\mathrm{ALG}$.

\paragraph{Step 2:} 
Fix any valid tuple $(X_1,\dots,X_{r_+})$ such that $|X_1|<\lambda$. Assume there exists a feasible solution $z$ to $\mathsf{P}'(I)$ that corresponds to the valid tuple, and let $Z =\{i\in[m]\mid z_i=1\}$. Then because $|X_1|=\min\{\lambda,|Z|\}=|Z|$ and $X_1\subset Z$, we have $X_1=Z$. Similarly, $X_j=Z$ for all $j\in[r_+]$. Therefore, there exists a feasible solution to $\mathsf{P}'(I)$ that corresponds to $(X_1,\dots,X_{r_+})$ only if $X_1=\cdots=X_{r+}$. There are at most $\sum_{i=1}^{\lambda-1} {m\choose i}\leq\lambda m^\lambda$ such tuples. Moreover, there is a unique feasible solution $z$ that corresponds to $(X_1,\dots,X_{r_+})$, namely $z_i=1$ for every $i\in X_1$ and $z_i=0$ for every $i\notin X_1$. Thus, the runtime of enumerating all corresponding feasible solutions is bounded by $\lambda m^\lambda$.

\paragraph{Step 3:}
Fix any valid tuple $(X_1,\dots,X_{r_+})$ such that $|X_1|= \lambda$. Then by construction we must have  $z_i=1$ for every $i\in \cup_{j}X_j$, and $z_i=0$ for every $i\in \cup_{j}\hat{X}_j$. Therefore every feasible solution $z$ that corresponds to $(X_1,\dots,X_{r_+})$ must lie in the following set:
    \begin{align*}
    \{  z\in\{0,1\}^m\mid & z_i=1\text{ } \forall i\in \cup_{j}X_j, \\
    & z_i=0 \text{ }\forall i\in \cup_{j}\hat{X}_j,\\
    & 1\leq e^\top z\leq k\}.
    \end{align*}


There are ${m\choose\lambda}^{r_+}\leq m^{\lambda r_+}$ tuples such that $|X_1| = \lambda$. We enumerate all such tuples, and check each for validity according to the following procedure:

\begin{enumerate}\setcounter{enumi}{-1}
\item From Step 1 of $\mathrm{ALG}$, let $\pi_j: [m] \to [m]$ be a sorting of $[m]$ according to $d_j$ such that $d_{j,\pi_j(1)}\geq \cdots\geq d_{j,\pi_j(m)}$. 
    \item Fix any given tuple $(X_1,\dots,X_{r_+})$. For each $j\in[r_+]$, let $i(j)\in[m]$ be an index such that $d_{j,i(j)}\in X_j$ and $d_{j,i(j)}\leq d_{ji}$ for all $d_{ji}\in X_j$. Without loss of generality, if $d_{j,i(j)}=d_{ji'}$ and $d_{ji'}\notin X_j$ for some index $i'$, we let $\pi_j(i(j))>\pi_j(i')$ for tie-breaking in the sorting. Also, if $d_{j,i(j)}=d_{ji'}$ and $d_{ji'}\in X_j$ for some index $i'$, we let $\pi_j(i(j))<\pi_j(i')$ for tie-breaking in the sorting. Then by our construction \begin{eqnarray*} 
\hat{X}_j&=&\{i\in[m]\setminus X_j\mid d_{ji}> \min_{i'\in X_j}\{d_{ji'}\}\}\\&=&\{i\in[m]\setminus X_j\mid d_{ji}> d_{j,i(j)}\}\}\\&=&\{i\in[m]\setminus X_j\mid \pi_j(i)>\pi_j(i(j))\}.
\end{eqnarray*} Therefore $X_j\cup \hat{X}_j=\{i\in[m]\mid \pi_j(i)\geq\pi_j(i(j))\}$.
    \item Note that $\sum_{j=1}^{r_+}\lvert X_j\rvert=\lambda r_+$. Therefore, in order to check whether $|\cup_{j}X_j|\leq k$, we just need to count the number of overlaps among all $X_j$'s. Set a counter $c=0$ to count the overlaps. For each $j\in [r_+]$ and for each $i\in X_j$, we check if $\pi_{j'}(i)\geq \pi_{j'}(i(j'))$, and do the following:
    \begin{itemize}
        \item If $\pi_{j'}(i)< \pi_{j'}(i(j'))$, then we have $i\notin X_{j'}\cup \hat{X}_{j'}$. We do nothing in this case.
        \item If $\pi_{j'}(i)\geq \pi_{j'}(i(j'))$ and $\pi_{j'}(i)\in X_{j'}$, then $i$ appears in both $X_j$ and $X_{j'}$. Therefore we increase $c$ by 1.
        \item  If $\pi_{j'}(i)\geq \pi_{j'}(i(j'))$ and $\pi_{j'}(i)\notin X_{j'}$, then we must have $\pi_{j'}(i)\in \hat{X}_{j'}$. Therefore $(\cup_{j}X_j)\bigcap(\cup_{j}\hat{X}_j)\neq\emptyset$, so we can terminate the process and declare that $(X_1,\dots,X_{r_+})$ is not a valid tuple.    \end{itemize}
    We iterate all $j\in [r_+]$ and $i\in X_j$. Notice that $c$ counts the number of overlaps (with multiplicity) of elements in $X_j$, we have $|\cup_{j}X_j|=\sum_{j=1}^{r_+}\lvert X_j\rvert-c=\lambda r_+-c$. Therefore we can check if $|\cup_{j}X_j|\leq k$. Also, if the above procedure never encounters $\pi_{j'}(i)\in \hat{X}_{j'}$, then we have $(\cup_{j}X_j)\bigcap(\cup_{j}\hat{X}_j)=\emptyset$. Therefore this procedure allows us to check the validity of $(X_1,\dots,X_{r_+})$.

\end{enumerate} 

 Because each $d_{j'}$ is sorted, the above procedure takes a constant runtime for each $j'\in[r_+]$, so the runtime for each fixed $j\in [r_+]$ and $i\in X_j$ is $r_+$. Because $\sum_{j=1}^{r_+}\lvert X_j\rvert=\lambda r_+$, there are $\lambda r_+$ combinations of $j\in [r_+]$ and $i\in X_j$. Therefore the runtime to check the validity is $\lambda r_+^2$ for any given tuple $(X_1,\dots,X_{r_+})$. Because there are at most $m^{\lambda r_+}$ such tuples, the runtime of enumerating all such tuples is $\lambda r_+^2m^{\lambda r_+}$.

\end{myproof}

Below we give the pseudo-code of the algorithm $\textup{ALG}$ in Lemma~\ref{lemma: enumerate valid tuples}.

\begin{algorithm}[H]
\caption{Phase Two: Enumeration of Partial Solutions}
\label{alg: phase two enumeration}
\DontPrintSemicolon
\SetAlgoLined

\KwIn{Instance of Problem \ref{eqn:P'(I)} from Lemma \ref{lemma: replace W'}, oracle $\textup{ALG}'$ from Lemma \ref{lemma: enumerate valid tuples}, parameter $\delta > 0$}
\KwOut{Solution $x$ to $\mathsf{P}(I)$}
\BlankLine
Set $\lambda > 0$ according to Eq.~\eqref{eqn: lambda}\;

\BlankLine
\tcp{Enumerate valid tuples with $|X_j|<\lambda$}
\For{$\lambda' = 1$ \KwTo $\lambda - 1$}{
    \For{each $X_1 \subset [m]$ with $1 \leq |X_1| \leq k$}{
        Set $x \in \{0,1\}^m$ where $x_i = 1 \iff i \in X_1$; record $(x, f_{\mathsf{P'}(I)}(x))$\;
    }
}

\BlankLine
\tcp{Enumerate valid tuples with $|X_j|=\lambda$}

\For{$j \in [m]$}{
    $d_j \leftarrow b_j \odot (Vu)$\;
    $\pi_j \leftarrow$ permutation sorting $[m]$ by $d_j$ descending\;
}

\For{each tuple $(X_1, \ldots, X_{r_+})$ where $|X_j| = \lambda$ for all $j$}{
    \For{$j \in [r_+]$}{
        $i(j) \leftarrow \arg\min_{i \in X_j} d_{j,i}$\;
        $\hat{X}_j \leftarrow \{i \in [m] \setminus X_j : \pi_j(i) > \pi_j(i(j))\}$\;
    }
    
    $c \leftarrow 0$; $\textit{valid} \leftarrow \textbf{true}$\;
    
    \For{$j \in [r_+]$ and $i \in X_j$}{
        \For{$j' \in [r_+]$}{
            \lIf{$\pi_{j'}(i) \leq \pi_{j'}(i(j'))$}{$c \leftarrow c + 1$}
            \lIf{$\pi_{j'}(i) \leq \pi_{j'}(i(j'))$ and $i \notin X_{j'}$}{$\textit{valid} \leftarrow \textbf{false}$; \textbf{break}}
        }
        \lIf{not $\textit{valid}$}{\textbf{break}}
    }
    
    \If{$\textit{valid}$ and $c \leq \lambda r_+ - k$}{
        Let $x$ be the output of $\textup{ALG}'$ with inputs \ref{eqn:P(x)} where valid tuple $X=(X_1, \ldots, X_{r_+})$, and parameter $\delta$\; record $(x, f_{\mathsf{P'}(I)}(x))$\;
    }
}

\Return{recorded $x$ with maximum $f_{\mathsf{P'}(I)}(x)$}

\end{algorithm}

The consequence of this result is that we have reduced to the task of, for each valid tuple $X=(X_1,\dots,X_{r_+})$ with $|X_1|=\lambda$, solving $\mathsf{P}'(I)$ with the additional constraint that the solution must correspond to the valid tuple, as stated in Problem $\mathsf{P}(X)$ below:

\begin{align*} \label{eqn:P(x)}
    \tag{$\mathsf{P}(X)$}
    \max \quad & f_{\mathsf{P}(X)}(x)=\sum_{i=1}^m x_i f_i\left(\frac{\sum_{j=1}^{r_+} a_{ij} d_j^\top x}{w'_i{}^\top x}\right) \\
    \textup{s.t.} \quad & x \in \{0,1\}^m, \\
    & 1 \leq e^\top x \leq k, \\
    & x_i = 1 \quad \forall i \in \cup_j X_j, \\
    & x_i = 0 \quad \forall i \in \cup_j \hat{X}_j.
\end{align*}

\subsection{Step 3: Linearization}

In order to solve $\mathsf{P}(X)$, we define the following auxiliary problem $\mathsf{P}(X,t)$ for each $t\in\mathbb{R}_+^m$: 
\begin{align*}\label{eqn:P(X,t)}
    \tag{$\mathsf{P}(X,t)$}
    \max \quad & f_{\mathsf{P}(X,t)}(x)=\sum_{i=1}^m x_i f_i\left( \frac{\sum_{j=1}^{r_+} a_{ij} d_j^\top x}{t_i} \right) \\
    \textup{s.t.} \quad & x \in \{0,1\}^m,\\\quad &e^\top x \leq k, \\
    & w'_i{}^\top x \leq t_i \quad \forall i \in [m], \\
    & x_i = 1 \quad \forall i \in \cup_j X_j, \\
    & x_i = 0 \quad \forall i \in \cup_j \hat{X}_j.
\end{align*}

\noindent Note that we have dropped the constraint $1\leq e^\top x$ in $\mathsf{P}(X,t)$. This is inconsequential: because we assume $\textup{OPT}_{\mathsf{P}}$ is positive, the solution $x$ with all entries equal to zero is not an optimal solution to $\mathsf{P}$. Indeed, the only reason we have maintained the $1\leq e^\top x$ constraint until now has been to rule out notational edge cases (such as dividing by zero).

We prove an important property regarding the relationship between optimal solutions of $\mathsf{P}(X)$ and those of $\mathsf{P}(X,t)$.

\begin{lemma}\label{lemma: auxiliary problem}
    Fix any valid tuple $X$. Let $t^*\in\arg\max_{t\in\mathbb{R}_+^m}\textup{OPT}_{\mathsf{P}(X, t)}$, and let $x^*$ be an optimal solution to $\mathsf{P}(X, t^*)$. Then $W'x^*=t^*$, and $x^*$ is an optimal solution to $\mathsf{P}(X)$.
\end{lemma}

\begin{myproof}[Proof of Lemma \ref{lemma: auxiliary problem}.]
First, we show that $W'x^* =t^*$. Suppose otherwise, and let $t'=W'x^*$. Then $t'_i\leq t^*_i$ for every $i\in[m]$ and $t'_i< t^*_i$ for some $i$. Therefore $$\textup{OPT}_{\mathsf{P}(X, t^*)}=\sum_{i=1}^m x^*_i f_i\left(\frac{\sum_{j=1}^{r_+} a_{ij}d_j^\top x^*}{t^*_i}\right)<\sum_{i=1}^m x^*_i f_i\left(\frac{\sum_{j=1}^{r_+} a_{ij}d_j^\top x^*}{t'_i}\right)\leq \textup{OPT}_{\mathsf{P}(X, t')},$$ contradicting the definition of $t^*$.  

Now we show that $x^*$ is an optimal solution to the $\mathsf{P}(X)$. For the sake of contradiction, suppose that $x'$ gives a higher objective value than $x^*$ to $\mathsf{P}(X)$, that is, $$\sum_{i=1}^m x'_i f_i\left(\frac{\sum_{j=1}^{r_+} a_{ij}d_j^\top x'}{{w_i'}^\top x'}\right)> \sum_{i=1}^m x^*_i f_i\left(\frac{\sum_{j=1}^{r_+} a_{ij}d_j^\top x^*}{{w_i'}^\top x^*}\right)=\textup{OPT}_{\mathsf{P}(X, t^*)}.$$ Let $t'=Wx'$. Then $x'$ is a feasible solution to $P(X,t')$. Therefore, we have $$\textup{OPT}_{\mathsf{P}(X, t')}\geq \sum_{i=1}^m x'_i f_i\left(\frac{\sum_{j=1}^{r_+} a_{ij}d_j^\top x'}{{w_i'}^\top x'}\right)>\textup{OPT}_{\mathsf{P}(X, t^*)},$$ contradicting the definition of $t^*$. 
\end{myproof}

Because $t^*=W'x^*$ and $\sum_{j=1}^m w_{ij}\leq 1+\gamma$ for each $i\in [m]$, we have $t^*\in[W'_{\min},1+\gamma]^m$. Thus, from here on we will only consider the problems $\mathsf{P}(X,t)$ with $t\in[W'_{\min},1+\gamma]^m$. 

\subsection{Step 4: Dimensionality Reduction and Discretization of Auxiliary Problems:}

Lemma \ref{lemma: auxiliary problem} shows that, in order to solve $\mathsf{P}(X)$, it is enough to solve for $\arg\max_{t\in\mathbb{R}_+^m}\textup{OPT}_{\mathsf{P}(X, t)}$. 
Below we show that it suffices to solve the auxiliary problem  $\mathsf{P}(X, t)$ for a smaller, discretized set of $t$'s.

\begin{lemma}\label{lemma: auxiliary oracle is enough}
    Suppose we are given an oracle $\textup{ALG}'$ that takes  $\mathsf{P}(X,t)$ with any $t\in[W'_{\min},1+\gamma]^m$ and any $\delta>0$ as inputs, and outputs a solution of $\mathsf{P}(X,t)$ that satisfies $$\textup{ALG}'_{\mathsf{P}(X,t)}\geq (1-g'(\delta))\textup{OPT}_{\mathsf{P}(X,t)}-h'(\delta)$$ with runtime $T(\delta)$, where $0\leq g'(\delta)\leq 1$ and $h'(\delta)\geq 0$. Then there exists an algorithm $\textup{ALG}$ that takes $\mathsf{P}(X)$ and any $\delta>0$ as inputs, and outputs a solution of $\mathsf{P}(X)$ that satisfies 
    $$\textup{ALG}_{\mathsf{P}(X)}\geq (1-g'(\delta))\left(\left(1-g\left(\frac{(1+\gamma)\delta}{W'_{\min}}\right)\right)\textup{OPT}_{\mathsf{P}(X)}+kh\left(\frac{(1+\gamma)\delta}{W'_{\min}}\right)\right)-h'(\delta)$$ 
    with runtime 
 $$\left\lceil\left(\frac{4(1+\gamma)k}{\epsilon W'_{\min}}\right)^{r_+}\right\rceil T(\delta).$$
\end{lemma}


\begin{myproof}[Proof of Lemma \ref{lemma: auxiliary oracle is enough}.]

Recall that we have known non-negative factorization $W'=AB^\top$, where $A,B \in \mathbb{R}_{\geq 0}^{m \times r_+}$. First, we make the following observation on the scale of $A$ and $B$:
\begin{observation}\label{obs: scale}
    There exists $A',B'\in \mathbb{R}_{\geq 0}^{m\times r_+}$ where $W'=A'{B'}^\top$, such that $\lVert a'_i\rVert_1\leq 1+\gamma$ for every $i\in[m]$ and $\lVert b'_j\rVert_1= 1$ for every $j\in[r_+]$.
\end{observation}
\begin{myproof}[Proof of Observation \ref{obs: scale}.]
    We construct $A'$ and $B'$ explicitly. Let the rows of $B'$ be the rows of $B$ that are rescaled so that $\lVert b_j\rVert_1= 1$. That is, let $b'_{jk}=b_{jk}/\sum_{k'=1}^m b_{jk'}$ for every $j\in[r_+]$ and $k\in[m]$. Let the columns of $A'$ be the columns of $A$ that are rescaled accordingly. That is, let $a'_{ij}=a_{ij}\sum_{k'=1}^m b_{jk'}$ for every $i\in[m]$ and $j\in[r_+]$. Then we have $a'_{ij}b'_{jk}=a_{ij}b_{jk}$. Therefore $W'=A'{B'}^\top$. Finally, since $1-\gamma\leq W_{ij}/W'_{ij}\leq1+\gamma$ and $\sum_{j=1}^m w_{ij}=1$ for every $i\in[m]$, we have $\sum_{j=1}^m w'_{ij}\leq1+\gamma$ for every $i\in[m]$. Therefore for every $i\in[m]$, we have 
    $$1+\gamma\geq \sum_{j=1}^m w'_{ij}=\sum_{j=1}^{r_+} a'_{ij}\sum _{k=1}^{m}b'_{jk}=\sum_{j=1}^{r_+} a'_{ij}=\lVert a'_i\rVert_1.$$
\end{myproof}

By Observation \ref{obs: scale}, we may assume $\lVert a_i\rVert_1\leq 1+\gamma$ for every $i\in[m]$ and $\lVert b_j\rVert_1= 1$ for every $j\in[r_+]$ from now on. Let \begin{equation}\label{eqn:Y}
    Y=\{B^\top x\mid x\textup{ is a feasible solution to }\mathsf{P}'(I)\}\subset \mathbb{R}^{r_+}.
\end{equation}
We will partition the $t$-space $[W'_{\min},1+\gamma]^m$ by partitioning $Y$. Let $\delta'$ be the quantity
\begin{equation}\label{eqn:delta'}
    \delta'=\frac{W'_{\min}\delta}{(1+\gamma)\sqrt{k}},
\end{equation}
where the reason for this choice will be specified momentarily. Notice that $\lVert y\rVert_2\leq\sqrt{k}$ for every $y\in Y$. As seen in the proof of Proposition \ref{prop: phase one}, we can create a cover of a ball with radius $\sqrt{k}$ in $\mathbb{R}^{r_+}$ using $\lceil(4\sqrt{k}/\delta')^{r_+}\rceil$ number of balls with radius $\delta'/2$. Therefore we can create a partition $\mathcal{Y}=\{Y_1,\dots,Y_{\ell}\}$ of $Y$ such that $\ell=\lceil(4\sqrt{k}/\delta')^{r_+}\rceil$, and $\lVert y-y'\rVert_2\leq\delta'$ for every $\ell'\in[\ell]$ and $y,y'\in Y_{\ell'}$. 

Fix any row $a_i^\top$ of $A$. For every $\ell'\in[\ell]$ and $i,i'\in I_{\ell'}$, we have
\begin{align*}
    \left\lvert\frac{a_i^\top y}{a_{i}^\top y'}-1\right\rvert&=\left\lvert\frac{a_i^\top (y-y')}{a_{i}^\top y}\right\rvert\\
    &\leq \frac{\lVert a_i\rVert_2\lVert y-y'\rVert_2}{|a_i^\top y|}\\
    &\leq \frac{\delta' (1+\gamma)\sqrt{k}}{W'_{\min}}.
\end{align*}
Thus, for our particular choice of $\delta'$, we have that $\lvert a_i^\top y/a_{i}^\top y'-1\rvert\leq \delta$ for every $\ell'\in[\ell]$ and $y,y'\in Y_{\ell'}$.

For every $\ell'\in[\ell]$, let $T_{\ell'}=\{Ay\mid y\in Y_{\ell'}\}\subset\mathbb{R}^m$. Then for every $\ell'\in[\ell]$ and $t,t'\in T_{\ell'}$, we have that $\lvert (t_i/t_i')-1\rvert\leq \delta$ for every $i\in[m]$. Fix any $t_1\in T_1,\dots,t_{\ell}\in T_{\ell}$. The algorithm $\textup{ALG}$ will use the given oracle $\textup{ALG}'$ to obtain a solution for each $\mathsf{P}(X,t_{\ell'})$, and then output the solution that has the highest objective value of $\mathsf{P}(X,t_{\ell'})$. That is, $$\textup{ALG}_{\mathsf{P}(X)}=\max_{\ell'\in[\ell]} \textup{ALG}_{\mathsf{P}(X,t_{\ell '})}'.$$ 
Because $\ell=\lceil(4\sqrt{k}/\delta)^{r_+}\rceil=\lceil(4(1+\gamma)k/\epsilon W'_{\min})^{r_+}\rceil$, the runtime of $\textup{ALG}$ is $$\lceil(4(1+\gamma)k/\epsilon W'_{\min})^{r_+}\rceil T(\delta).$$

Finally, we prove the performance guarantee for $\textup{ALG}$. Let $t^*\in\arg\max_{t\in\mathbb{R}_+^m}\textup{OPT}_{\mathsf{P}(X, t)}$. Assume $t^*\in T_{\ell'}$. Then $\lvert (t^*_i/(t_{\ell'})_i)-1\rvert\leq \delta$. Let $x^*$ be the corresponding optimal solution to $\mathsf{P}(X, t^*)$. Then by Lemma \ref{lemma: auxiliary problem}, $x^*$ is an optimal solution to $\mathsf{P}(X)$. Let $x'$ be an optimal solution of $\mathsf{P}(X, t_{\ell'})$. Because $t_{\ell'} \in[W'_{\min},1+\gamma]^m$,  we have $\sum_{j=1}^{r_+} a_{ij}d_j^\top x^*/(t_{\ell'})_i={w'_i}^\top x^*/(t_{\ell'})_i\leq (1+\gamma)/W'_{\min}$. Also, since $f_i(x-\epsilon)\geq (1-g(\epsilon))f_i(x)-h(\epsilon)$ for all $x$, we have $f_i(x+\epsilon)\leq(f_i(x)+h(\epsilon))/(1-g(\epsilon))$ for all $x$. As a consequence of Lemma \ref{lemma: auxiliary problem}, we have $\textup{OPT}_{\mathsf{P}(X)}= \textup{OPT}_{\mathsf{P}(X,t^*)}$. Then 
\begin{align*}
    \textup{OPT}_{\mathsf{P}(X)}&= \textup{OPT}_{\mathsf{P}(X,t^*)}\\
    &=\sum_{i=1}^m x^*_i f_i\left(\frac{\sum_{j=1}^{r_+} a_{ij}d_j^\top x^*}{t^*_i}\right)\\
    &\leq \sum_{i=1}^m x^*_i f_i\left((1+\delta)\frac{\sum_{j=1}^{r_+} a_{ij}d_j^\top x^*}{(t_{\ell'})_i}\right)\\
    &\leq \frac{\sum_{i=1}^m x^*_i f_i\left(\frac{\sum_{j=1}^{r_+} a_{ij}d_j^\top x^*}{(t_{\ell'})_i}\right)+kh(\delta(1+\gamma)/W'_{\min})}{1-g(\delta (1+\gamma)/W'_{\min})}\\
    &\leq \frac{\textup{OPT}_{\mathsf{P}(X,t_{\ell'})}+kh(\delta(1+\gamma)/W'_{\min})}{1-g(\delta (1+\gamma)/W'_{\min})}.
\end{align*}
Rearranging the above, we have $$\textup{OPT}_{\mathsf{P}(X,t_{\ell'})}\geq (1-g(\delta (1+\gamma)/W'_{\min}))\textup{OPT}_{\mathsf{P}(X)}-kh(\delta(1+\gamma)/W'_{\min}).$$
Therefore,
\begin{eqnarray*}
    \textup{ALG}_{\mathsf{P}(X)}&\geq&  \textup{ALG}'_{\mathsf{P}(X,t_{\ell'})}\\
     &\geq& (1-g'(\delta))\textup{OPT}_{\mathsf{P}(X,t_{\ell'})}-h'(\delta)\\
     &\geq& (1-g'(\delta))\left(\left(1-g\left(\frac{(1+\gamma)\delta}{W'_{\min}}\right)\right)\textup{OPT}_{\mathsf{P}(X)}-kh\left(\frac{(1+\gamma)\delta}{W'_{\min}}\right)\right)-h'(\delta).
\end{eqnarray*} 
\end{myproof}

Below we give the pseudo-code of the algorithm $\textup{ALG}$ in Lemma~\ref{lemma: auxiliary oracle is enough}.

\begin{algorithm}[H]
\caption{Phase Two: Discretization of Auxiliary Problems}
\label{alg: phase two discretization}
\DontPrintSemicolon
\SetAlgoLined

\KwIn{Instance of Problem \ref{eqn:P(x)}, oracle $\textup{ALG}'$ from Lemma \ref{lemma: auxiliary oracle is enough}, parameter $\delta > 0$}
\KwOut{Solution $x$ to $\mathsf{P}(X)$}

\BlankLine

\tcp{Discretize the space of auxiliary variable}
Set $\delta' > 0$ according to Eq.~\eqref{eqn:delta'}\;

Construct maximal $(\delta'/2)$-separated set $C \subset B(0, \sqrt{k})$ with $|C| = \lceil(4\sqrt{k}/\delta')^{r_+}\rceil$\;

Let $\mathcal{Y} = \{Y_1, \ldots, Y_\ell\} \leftarrow \{B(c, \delta'/2) : c \in C\}$\;

\BlankLine
\tcp{Solve each discretized auxiliary problem}

Set $\delta_1\gets \delta$ and $\delta_2\gets k\delta$\;
\For{$\ell' = 1$ \KwTo $\ell$}{
    $T_{\ell'} \leftarrow \{Ay : y \in Y_{\ell'}\}$\;
    Choose arbitrary $t_{\ell'} \in T_{\ell'}$\;
    Let $x$ be the output of $\textup{ALG}'$ with inputs $\mathsf{P}(X,t_{\ell'})$ and parameters $\delta_1,\delta_2$\; 
    Record $(x_{\ell'}, f_{\mathsf{P}(X, t_{\ell'})}(x_{\ell'}))$\;
}

\BlankLine
\Return{$\arg\max_{x_{\ell'}} f_{\mathsf{P}(X, t_{\ell'})}(x_{\ell'})$ among recorded solutions}

\end{algorithm}

\subsection{Step 5: Complete Linearization of Auxiliary Problems}
Lemma \ref{lemma: auxiliary oracle is enough} shows that, to approximately solve $\mathsf{P}(X)$, it suffices to construct an oracle that approximately solves $\mathsf{P}(X,t)$ for any given $t\in[W'_{\min},1+\gamma]^m$. Below we give such an oracle. Let $c_i=a_i/t_i\in\mathbb{R}^{r_+}_{\geq 0}$ for each $i\in[m]$. 
Then $\mathsf{P}(X,t)$ can be equivalently formulated as
\begin{align*}
    \tag{$\mathsf{P}(X,t)$}
    \max \quad & f_{\mathsf{P}(X,t)}(x)=\sum_{i=1}^m x_i f_i\left( \sum_{j=1}^{r_+} c_{ij} d_j^\top x \right) \\
    \textup{s.t.} \quad & x \in \{0,1\}^m,\\\quad &e^\top x \leq k, \\
    & w'_i{}^\top x \leq t_i \quad \forall i \in [m], \\
    & x_i = 1 \quad \forall i \in \cup_j X_j, \\
    & x_i = 0 \quad \forall i \in \cup_j \hat{X}_j.
\end{align*}

    To solve $\mathsf{P}(X,t)$, we partition the space of possible values $(d_1^\top x,\dots,d_{r_+}^\top x)\in\mathbb{R}^{r_+}$, as well as the space of the objective value $f_{\mathsf{P}(X,t)}(x)$. 

\begin{lemma} \label{lemma: auxiliary oracle on each tuple is enough}
    Fix any $t\in[W'_{\min},1+\gamma]^m$. Suppose we are given an oracle with runtime $T(\delta_1,\delta_2)$ that takes $\mathsf{P}(X,t)$, any  $\theta=(\theta_1,\dots,\theta_{r_+})\in\mathbb{R}^{r_+}$, any $\zeta \ge 0$, and any $\delta_1,\delta_2>0$ as inputs, and either \begin{enumerate}
        \item Scenario one: correctly declares that there is no feasible $x$ to $\mathsf{P}(X,t)$ such that $d_j^\top x\geq \theta_j$ for every $j\in[r_+]$ and $\sum_{i=1}^m x_i f_i(c_i^\top \theta)\geq \zeta$, or
        \item Scenario two: outputs a feasible $x$ to $\mathsf{P}(X,t)$ such that  $d_j^\top x+\delta_1\geq \theta_j$ for every $j\in[r_+]$ and $\sum_{i=1}^m x_i f_i(c_i^\top \theta)+\delta_2\geq \zeta$.
    \end{enumerate}
     Then there exists an algorithm $\textup{ALG}$ that satisfies $$\textup{ALG}_{\mathsf{P}(X,t)}\geq\left(1-g\left(\frac{\delta_1(1+\gamma)}{W'_{\min}}\right)\right)(\textup{OPT}_{\mathsf{P}(X,t)}-2\delta_2)-kh\left(\frac{\delta_1(1+\gamma)}{W'_{\min}}\right)$$ with runtime $$\lceil ((Vu)_{\max}-\min\{0,(Vu)_{\min})\}/\delta_1\rceil^{r_+}\cdot\left\lceil k\max_{i\in[m]}\{f_i((Vu)_{\max})\}/\delta_2\right\rceil T(\delta_1,\delta_2).$$ 
\end{lemma}

\begin{myproof}[Proof of Lemma \ref{lemma: auxiliary oracle on each tuple is enough}.]
Let $(Vu)_{\min}$ be the minimum entry of $Vu$ (possibly negative). Because $\lVert b_j\rVert_1=1$ and $d_j=b_j\odot(Vu)$ for every $j\in[r_+]$, 
we have $d_j^\top x^*\in[\min\{0,(Vu)_{\min}\},(Vu)_{\max}]$ for every $j \in [r_+]$. Also, because each $f_i$ is non-decreasing, $\textup{OPT}_{\mathsf{P}(X,t)}\in[0,k\max_{i\in[m]}\{f_i((Vu)_{\max})\}]$. Similar to the proof of Lemma \ref{lemma: auxiliary oracle is enough}, we will create a partition of the space of possible values $(d_1^\top x,\dots,d_{r_+}^\top x)\in\mathbb{R}^{r_+}$, as well as the space of the objective value $f_{\mathsf{P}(X,t)}(x)$. We then show that it is sufficient to solve $\mathsf{P}(X,t)$ in each subset of the partition. Let $\Delta_{\ell}=\min\{0,(Vu)_{\min}\}+\delta_1(\ell-1)$ for $\ell=1,\dots, \lceil ((Vu)_{\max}-\min\{0,(Vu)_{\min}\})/\delta_1\rceil$. Let $\Delta'_{s}=\delta_2(s-1)$ for $s=1,\dots, \lceil k\max_{i\in[m]}\{f_i((Vu)_{\max})\}/\delta_2\rceil$. Consider all tuples $(\ell_1,\dots,\ell_{r_+},s)$. There are in total $$\lceil ((Vu)_{\max}-\min\{0,(Vu)_{\min}\})/\delta_1\rceil^{r_+}\cdot\left\lceil k\max_{i\in[m]}\{f_i((Vu)_{\max})\}/\delta_2\right\rceil$$ such tuples. Moreover, there exists a tuple $(\ell^*_1,\dots,\ell^*_{r_+},s^*)$ such that $\Delta_{\ell^*_i}\leq x^*_i\leq \Delta_{\ell^*_i+1}$ for each $i\in [m]$ and $\Delta'_{s^*}\leq\textup{OPT}_{\mathsf{P}(X,t)}\leq \Delta'_{s^*+1}$.

For each tuple $(\ell_1,\dots,\ell_{r_+},s)$, our desired algorithm uses the given oracle to determine whether there exists a feasible $x$ to $\mathsf{P}(X,t)$ that satisfies the conditions in scenario two with $\theta_i=\Delta_{\ell_i}$ for each $i\in[r_+]$ and $\zeta=\Delta'_{s}$. Then our desired algorithm returns the $x$ with the highest objective value of $\mathsf{P}(X,t)$ among all tuples. Note that $x^*$ satisfies the conditions in scenario two on the tuple $(\ell^*_1,\dots,\ell^*_{r_+},s^*)$. Therefore the given oracle would return some feasible $x'$ to $\mathsf{P}(X,t)$ that satisfies the conditions in scenario two with $\theta_i=\Delta_{\ell^*_i}$ for each $i\in[r_+]$ and $\zeta=\Delta'_{s^*}$. Notice that $c_{ij}=a_{ij}/t_i\leq (1+\gamma)/W'_{\min}$. Then by the conditions in scenario two we have \begin{eqnarray*}
     \textup{ALG}_{\mathsf{P}(X,t)}&\geq& \sum_{i=1}^m x'_if_i\left( \sum_{j=1}^{r_+}c_{ij}d_j^\top x'\right)\\
     &\geq& \sum_{i=1}^m x'_if_i\left( \sum_{j=1}^{r_+}c_{ij}(\theta_j-\delta_1)\right)\\
      &\geq& \sum_{i=1}^m x'_if_i\left( \sum_{j=1}^{r_+}c_{ij}\theta_j-\delta_1 (1+\gamma)/W'_{\min}\right)\\
        &\geq& (1-g(\delta_1 (1+\gamma)/W'_{\min}))\sum_{i=1}^m x'_if_i\left( \sum_{j=1}^{r_+}c_{ij}\theta_j\right)-kh(\delta_1 (1+\gamma)/W'_{\min})\\  &\geq& (1-g(\delta_1 (1+\gamma)/W'_{\min}))(\Delta'_{s^*}-\delta_2)-kh(\delta_1 (1+\gamma)/W'_{\min})\\
 &\geq& (1-g(\delta_1 (1+\gamma)/W'_{\min}))(\textup{OPT}_{\mathsf{P}(X,t)}-2\delta_2)-kh(\delta_1 (1+\gamma)/W'_{\min}),
\end{eqnarray*} where the second and the fifth inequalities follow from the conditions in scenario two, and the last inequality follows since $\Delta'_{s^*}\leq\textup{OPT}_{\mathsf{P}(X,t)}\leq \Delta'_{s^*+1}$.

Because there are in total $$\lceil ((Vu)_{\max}-\min\{0,(Vu)_{\min}\})/\delta_1\rceil^{r_+}\cdot\left\lceil k\max_{i\in[m]}\{f_i((Vu)_{\max})\}/\delta_2\right\rceil$$ number of tuples $(\ell_1,\dots,\ell_{r_+},s)$, the runtime of our algorithm is $$\left\lceil ((Vu)_{\max}-\min\{0,(Vu)_{\min}\})/\delta_1\right\rceil^{r_+}\cdot\left\lceil k\max_{i\in[m]}\{f_i((Vu)_{\max})\}/\delta_2\right\rceil T(\delta_1,\delta_2).$$ 
\end{myproof}

Below we give the pseudo-code of the algorithm $\textup{ALG}$ in Lemma~\ref{lemma: auxiliary oracle on each tuple is enough}.

\begin{algorithm}[H]
\caption{Phase Two: Complete Linearization of Auxiliary Problems}
\label{alg: phase two linearization}
\DontPrintSemicolon
\SetAlgoLined

\KwIn{Instance of Problem \ref{eqn:P(X,t)}, oracle from Lemma \ref{lemma: auxiliary oracle on each tuple is enough}, parameters $\delta_1, \delta_2 > 0$}
\KwOut{Solution $x$ to $\mathsf{P}(X,t)$}

\BlankLine
\tcp{Discretize value and objective spaces}
Let $\nu_{\min} \leftarrow \min\{0, (Vu)_{\min}\}$ and $\nu_{\max} \leftarrow (Vu)_{\max}$\;
Let $L \leftarrow \lceil (\nu_{\max} - \mu_{\min})/\delta_1 \rceil$ and $S \leftarrow \lceil k \cdot \max_{i \in [m]} f_i(\nu_{\max})/\delta_2 \rceil$\;

\For{$\ell = 1$ \KwTo $L$}{
    $\Delta_\ell \leftarrow \nu_{\min} + \delta_1(\ell - 1)$\;
}

\For{$s = 1$ \KwTo $S$}{
    $\Delta'_s \leftarrow \delta_2(s - 1)$\;
}

\BlankLine
\tcp{Enumerate and solve linearized problems}
\For{each tuple $(\ell_1, \ldots, \ell_{r_+}, s)$ where $\ell_i \in [L]$ and $s \in [S]$}{
    Set $\theta_i \leftarrow \Delta_{\ell_i}$ for each $i \in [r_+]$ and $\zeta \leftarrow \Delta'_s$\;
    
    Run oracle from Lemma \ref{lemma: auxiliary oracle on each tuple is enough} with inputs 
         $\mathsf{P}(X,t)$, value tuple $(\theta_1, \ldots, \theta_{r_+}, \zeta)$, and parameters $\delta_1, \delta_2$\;
    
    \If{oracle returns scenario two with feasible solution $x$}{
        Record $(x, f_{\mathsf{P}(X,t)}(x))$\;
    }
}

\BlankLine
\Return{$\arg\max_x f_{\mathsf{P}(X,t)}(x)$ among recorded solutions}

\end{algorithm}

\subsection{Step 6: Approximation of Linearized Auxiliary Problems via LP Rounding}

Lemma \ref{lemma: auxiliary oracle on each tuple is enough} shows that, to solve $P(X, t)$ for any given $t\in[W'_{\min},1+\gamma]^m$, it is enough to give an oracle described in Lemma \ref{lemma: auxiliary oracle on each tuple is enough}. Similar to the idea of enumerating partial solutions by constructing valid tuples based on the values of each $d_j$, we further construct index sets based on the values of $f_i(c_i^\top \theta)$ and enumerate all possible index sets. Recall that via the valid tuple $(X_1,\dots,X_{r_+})$, we have already fixed at least $\lambda$ indices of any feasible solution to $\mathsf{P}(X,t)$ to be equal to $1$, namely the indices in $\cup_{j} X_j$. Fix \begin{equation}\label{eqn: lambda'}
    \lambda'=\lceil(2r_++2)\max_{i\in[m]}\{f_i((Vu)_{\max})\}/\epsilon\rceil.
\end{equation} Let $X'\subset [m]\setminus (\cup_{j} X_j)\cup(\cup_j \hat{X}_j)$ be an index set such that $0\leq |X'|\leq \lambda'$. Let 
\begin{equation} \label{eqn: X'}
\hat{X'}=\{i\in[m]\setminus X'\cup(\cup_{j} X_j)\cup(\cup_j \hat{X}_j)\mid f_i(c_i^\top \theta)> \min_{i'\in X'}\{f_{i'}(c_{i'}^\top \theta)\}\}.
\end{equation}
In words, $\hat{X'}$ consists of the indices outside of $X'\cup(\cup_{j} X_j)\cup(\cup_j \hat{X}_j)$ whose corresponding values of $f_i(c_{i'}^\top \theta)$ are strictly greater than the minimum value of $f_i(c_{i}^\top \theta)$ across indices in $X'$. Then every feasible solution $z$ to $\mathsf{P}(X,t)$ \textbf{corresponds} to an index set $X'$ in the following way: 
let $Z =\{i\in[m]\setminus(\cup_{j} X_j)\cup(\cup_j \hat{X}_j)\mid z_i=1\}$. We define $X'\subset Z$ to be the set of indices $i \in Z$ such that $f_i(c_i^\top \theta)$ is among the $\min\{\lambda',|Z|\}$ highest values in $Z$. That is, 
let $\pi': [|Z|] \to Z$ be a sorting of $Z$ according to $f_i(c_{i}^\top \theta)$ such that $f_{\pi'(1)}(c_{\pi'(1)}^\top \theta)\geq \cdots\geq f_{\pi'(|Z|)}(c_{\pi'(|Z|)}^\top \theta)$. Then $X'=\{\pi'(1),\dots,\pi'(\min\{\lambda',|Z|\})\}$. Similar to before, if $z$ corresponds to $X'$, we must have $z_i=1$ for $i\in X'$ and $z_i=0$ for $i\in \hat{X'}$.

As in Lemma \ref{lemma: enumerate valid tuples}, the notion of correspondence to index sets forms a partition of the set of feasible solutions to $\mathsf{P}(X,t)$. Therefore, in order to give an oracle described in Lemma \ref{lemma: auxiliary oracle on each tuple is enough}, it suffices to give an oracle described in Lemma \ref{lemma: auxiliary oracle on each tuple is enough} separately for each subset of this partition. This is formally stated in the following result:

\begin{observation} \label{obs: partition value}
    Suppose we are given an oracle with runtime $T(\delta_1,\delta_2)$ that takes $\mathsf{P}(X,t)$, any  $\theta=(\theta_1,\dots,\theta_{r_+})\in\mathbb{R}^{r_+}$, any $\zeta \ge 0$, any $\delta_1,\delta_2>0$, and any index set $X'\subset [m]\setminus (\cup_{j} X_j)\cup(\cup_j \hat{X}_j)$ such that $0\leq |X'|\leq \lambda'$ as inputs, and either \begin{enumerate}
        \item Scenario one: correctly declares that there is no feasible $x$ to $\mathsf{P}(X,t)$ that corresponds to $X'$ such that $d_j^\top x\geq \theta_j$ for every $j\in[r_+]$ and $\sum_{i=1}^m x_i f_i(c_i^\top \theta)\geq \zeta$, or
        \item Scenario two: outputs a feasible $x$ to $\mathsf{P}(X,t)$ that corresponds to $X'$ such that  $d_j^\top x+\delta_1\geq \theta_j$ for every $j\in[r_+]$ and $\sum_{i=1}^m x_i f_i(c_i^\top \theta)+\delta_2\geq \zeta$.
    \end{enumerate}
    Then there exists an oracle described in Lemma \ref{lemma: auxiliary oracle on each tuple is enough} with runtime $\lambda'm^{\lambda'}T(\delta_1,\delta_2)$.
\end{observation}

\begin{myproof}[Proof of Observation \ref{obs: partition value}.]
    Because every feasible solution $x$ to $\mathsf{P}(X,t)$ corresponds to exactly one index set $X'\subset [m]\setminus (\cup_{j} X_j)\cup(\cup_j \hat{X}_j)$ such that $0\leq |X'|\leq \lambda'$, we can give an oracle described in Lemma \ref{lemma: auxiliary oracle on each tuple is enough} by enumerating all such index sets $X'$ and applying the oracle  in Observation \ref{obs: partition value}. 
    
    More specifically, for the oracle described in Lemma \ref{lemma: auxiliary oracle on each tuple is enough} with inputs $\theta,\zeta,\delta_1,\delta_2$:
    \begin{itemize}
        \item If the oracle in Observation \ref{obs: partition value} outputs an $x$ in scenario two with inputs $\theta,\zeta,\delta_1,\delta_2,X'$ for some $X'$, then  the oracle described in Lemma \ref{lemma: auxiliary oracle on each tuple is enough} also outputs this $x$ in scenario two.
        \item If the oracle in Observation \ref{obs: partition value} declares scenario one with inputs $\theta,\zeta,\delta_1,\delta_2,X'$ for all $X'$, then the oracle described in Lemma \ref{lemma: auxiliary oracle on each tuple is enough} also declares scenario one.
    \end{itemize}
We can enumerate all index sets $X'\subset [m]\setminus (\cup_{j} X_j)\cup(\cup_j \hat{X}_j)$ such that $0\leq |X'|\leq \lambda'$ by simply enumerating all combinations of $|X'|$ indices from $[m]\setminus (\cup_{j} X_j)\cup(\cup_j \hat{X}_j)$. Because there are at most $\sum_{k=0}^{\lambda'} {m\choose\lambda'}\leq \lambda'm^{\lambda'}$ such index sets $X'$, the runtime of the oracle described in Lemma \ref{lemma: auxiliary oracle on each tuple is enough} is $\lambda'm^{\lambda'}T(\delta_1,\delta_2)$.
\end{myproof}

By Observation \ref{obs: partition value}, it suffices to give an oracle as described. Below we give such an oracle.

First, suppose $|X'|<\lambda'$. Assume there exists a feasible solution $z$ to $\mathsf{P}(X,t)$ that corresponds to $X'$. Let $Z =\{i\in[m]\setminus (\cup_{j} X_j)\cup(\cup_j \hat{X}_j)\mid z_i=1\}$. Then because $|X'|=\min\{\lambda',|Z|\}=|Z|$ and $X'\subset Z$, we have $X'=Z$. Therefore, there is a unique feasible solution $z$ of $\mathsf{P}(X,t)$ that corresponds to $X'$, namely $z_i=1$ for every $i\in X'\cup (\cup_{j} X_j)$ and $z_i=0$ for every $i\notin X'\cup (\cup_{j} X_j)$. Therefore, if $|X'|<\lambda'$, the oracle in Observation \ref{obs: partition value} can directly check this unique feasible solution of $\mathsf{P}(X,t)$ that corresponds to $X'$, and outputs the correct scenario accordingly.

From now on, we assume $|X'|=\lambda'$. Fix any $\theta\in\mathbb{R}^{r_+}$, any $\zeta\in \mathbb{R}_+$,  any $\delta_1,\delta_2>0$, and any $X'$. The oracle essentially needs to determine the existence of a feasible binary solution to a system of linear inequalities, which is NP-hard in general. However, the linear constraints of $\mathsf{P}(X,t)$ lie in a lower-dimensional subspace, a structure we can exploit by solving a relaxation of the system, obtained by replacing binary variables with continuous ones, and rounding its solution back to a binary solution. Because of the rounding, it is possible that the values $(d_1^\top x,\dots,d_{r_+}^\top x)$ and $\sum_{i=1}^m x_i f_i(c_i^\top \theta)$ of the rounded solution are out of the desired ranges. However, the valid tuple $X$ and the index set $X'$ we fixed before ensures that the gaps between the values and the desired ranges are within small constants. 

We define a polyhedron $PH\subset \mathbb{R}^{m}$ as follows:
\begin{equation}\tag{$PH$}\label{eqn: PH}
\begin{aligned}
\sum_{j=1}^{r_+} a_{ij}\, b_j^\top x &\le t_i &&\quad\text{for } i=1,\dots,m,\\
e^\top x &\le k,\\
d_j^\top x &\ge \theta_j &&\quad\text{for } j=1,\dots,r_+,\\
\sum_{i=1}^m x_i\, f_i(c_i^\top \theta) &\ge \zeta,\\
x_i &= 1 &&\quad\text{for } i\in \big(\bigcup_j X_j\big)\cup X',\\
x_i &= 0 &&\quad\text{for } i\in \big(\bigcup_j \hat X_j\big)\cup \hat X',\\
x_i &\in [0,1] &&\quad\text{for } i\notin \big(\bigcup_j X_j\big)\cup\big(\bigcup_j \hat X_j\big)\cup X'\cup \hat X'.
\end{aligned}
\end{equation}
Then $PH$ is a polyhedron in $\mathbb{R}^{m}$ defined by at most $3m+r_++2$ inequalities. Let $\textup{LP}(m,n)$ be the runtime of solving a linear program with $m$ variables and $n$ constraints. Then checking whether $PH$ is non-empty and return a point in $PH$ if $PH$ is non-empty can be done in runtime $\text{LP}(m,3m+r_++2)$. For more on the runtime of solving a linear program, we refer the readers to e.g. \cite{khachiyan1980polynomial,grotschel1981ellipsoid} (ellipsoid methods) and \cite{nesterov1994interior,vaidya1996new} (interior point methods). In what follows we assume that $PH \neq\emptyset$, otherwise the oracle outputs scenario one. 

\begin{lemma}\label{lemma: constant fractional components}
    If $PH\neq\emptyset$, then we can find a point $z\in PH$ with at most $2r_++2$ fractional components in runtime $\textup{LP}(m,3m+r_++2)+\textup{LP}(m,2m+2r_++2)$.
\end{lemma}

\begin{myproof}[Proof of Lemma \ref{lemma: constant fractional components}.]
Let $z^*\in PH$ be an (arbitrary) point found in runtime $\textup{LP}(m,3m+r_++2)$. Let $PH(z^*)\subset \mathbb{R}^{m-|(\cup_{j} X_j)\cup(\cup_j \hat{X}_j)\cup X'\cup \hat{X'}|}$ be the polyhedron on variable $y$, where the index set of $y$ is taken to be $I_y=[m]\setminus (\cup_{j} X_j)\cup(\cup_j \hat{X}_j)\cup X'\cup \hat{X'}$, with the following constraints:
\begin{equation}\tag{$PH(z^*)$}\label{eqn: PH(z)}
\begin{aligned}
\sum_{i\in I_y} b_{ji}\, y_i &\le \sum_{i\in I_y} b_{ji}\, z_i^* &&\text{for } j=1,\dots,r_+,\\
\sum_{i\in I_y} y_i &\le \sum_{i\in I_y} z_i^*,\\
\sum_{i\in I_y} d_{ji}\, y_i &\ge \theta_j - \sum_{i\in (\bigcup_j X_j)\cup X'} d_{ji} &&\text{for } j=1,\dots,r_+,\\
\sum_{i\in I_y} y_i\, f_i(c_i^\top\theta) &\ge \zeta - \sum_{i\in (\bigcup_j X_j)\cup X'} f_i(c_i^\top\theta),\\
y_i &\in [0,1] &&\text{for } i\in I_y.
\end{aligned}
\end{equation}

Note that $PH(z^*)\neq\emptyset$ since the projection of $z^*$ on $\mathbb{R}^{|I_y|}$ is in $PH(z^*)$. Because $PH(z^*)$ has $2r_++2$ linear inequalities other than the inequalities $y_i \in[0,1] \text{ for }i\in I_y$, we can compute a vertex $y^*$ of $PH(z^*)$ with at most $2r_++2$ fractional components with runtime $\textup{LP}(m,2m+2r_++2)$ (see a standard textbook on linear programming, e.g., \cite{schrijver1998theory}).

Let $z\in[0,1]^{m}$ where \begin{equation*}
    z_i = \begin{cases}
            1 & \text{if } i \in (\cup_j X_j) \cup X'\\
            0 & \text{if } i \in (\cup_j \hat{X}_j) \cup \hat{X'}\\
            y^*_i & \text{if } i \in I_y
        \end{cases}.
\end{equation*} Then $z$ has at most $2r_++2$ fractional components. We show that $z\in PH$. Because $z_i=z^*_i=1$ for $i\in (\cup_{j} X_j)\cup X'$ and $z_i=z^*_i=0$ for $i\in (\cup_{j} \hat{X}_j)\cup\hat{X'}$, the last three sets of constraints of $PH$ is satisfied. By the first set of constraints of $PH(z^*)$ we have $\sum_{i\in I_y} b_{ji}z_i \leq \sum_{i\in I_y} b_{ji}z^*_i$. Because $a_{ij}\geq0$ for every $i,j$, 
the first set of constraints of $PH$ is satisfied. Similarly the second constraint of $PH$ is also satisfied. By the third set of constraints of $PH(z^*)$ we have $$d_j^\top z=\sum_{i\in (\cup_{j} X_j)\cup X'} d_{ji}+\sum_{i\in I_y} d_{ji} z_i\geq \sum_{i\in (\cup_{j} X_j)\cup X'} d_{ji}+\left(\theta_j-\sum_{i\in (\cup_{j} X_j)\cup X'} d_{ji}\right)=\theta_j,$$ so the third set of constraints of $PH$ is satisfied. Similarly the fourth constraint of $PH$ is also satisfied. Therefore $z\in PH$ is the desired point.
\end{myproof}

Let $z$ be the point obtained in Lemma \ref{lemma: constant fractional components}. Then $z$ satisfies all the constraints of $\mathsf{P}(X,t)$ except the integrality constraints. 
We round $z$ down to obtain a feasible solution: let $\bar{z}\in\{0,1\}^m$ where $\bar{z}_i=\lfloor z_i\rfloor$ for each $i$. Notice that since $w'_{ij}>0$ for every $i,j$, we have $e^\top \bar{z}\leq e^\top z\leq k$ and ${w'_i}^\top\bar{z}\leq {w'_i}^\top z\leq t_i$ for every $i\in[m]$. Therefore $\bar{z}$ is feasible to $\mathsf{P}(X,t)$. Moreover, since $\bar{z}=1$ for $i\in X'$ and $\bar{z}=0$ for $i\in\hat{X'}$, we have that $\bar{z}$ corresponds to $X'$. 

In the final step, we show that by setting $\lambda,\lambda'$ appropriately, $\bar{z}$ satisfies the conditions in scenario two of Observation \ref{obs: partition value}, hence completing the oracle in Observation \ref{obs: partition value}.

\begin{lemma}\label{lemma: z satisfies auxiliary oracle}
    Set $\lambda=\lceil(2r_++2)(Vu)_{\max}/\delta_1\rceil$ and $\lambda'=\lceil(2r_++2)k\max_{i\in[m]}\{f_i((Vu)_{\max})\}/\delta_2\rceil$. Then $d_j^\top \bar{z}+\delta_1\geq \theta_j$ for every $j\in[r_+]$, and $\sum_{i=1}^m \bar{z}_i f_i(c_i^\top \theta)+\delta_2\geq \zeta$.
\end{lemma}

\begin{myproof}[Proof of Lemma \ref{lemma: z satisfies auxiliary oracle}.]
Because $z\in PH$, we have $d_j^\top z\geq \theta_j$ for every $j\in[r_+]$ and $\sum_{i=1}^m z_i f_i(c_i^\top \theta)\geq \zeta$. Fix $j\in[r_+]$ and let $\ell\in X_j$ be an index where $d_{j\ell}=\min_{\ell'\in X_j}\{d_{j\ell'}\}$. Then since $|X_j|=\lambda$, we have $$d_j^\top z\geq \sum_{\ell'\in X_j}d_{j\ell'}z_{\ell'}\geq \lambda d_{j\ell}.$$ On the other hand,  $\bar{z}$ is obtained by rounding $z$ down. Notice that $z_{\ell'}\in\{0,1\}$ for all $\ell'\in X_j\cup\hat{X}_j$, that is, for all $\ell'$ such that $d_{j\ell'}>d_{j\ell}$. Therefore for all $\ell'$ such that $d_{j\ell'}>d_{j\ell}$, we have $z_{\ell'}=\bar{z}_{\ell'}$. By Lemma \ref{lemma: constant fractional components}, $z$ has at most $2r_++2$ fractional components. 
Therefore, because $\theta_j\leq d_j^\top z\leq(Vu)_{\max}$, we have \begin{eqnarray*}
    d_j^\top \bar{z}&\geq& d_j^\top z - (2r_++2)d_{j\ell}\\&\geq& d_j^\top z-(2r_++2)d_j^\top z/\lambda\\
    &\geq& \theta_j-(2r_++2)(Vu)_{\max}/\lambda\\
    &\geq& \theta_j-\delta_1.
\end{eqnarray*}

Similarly, let $p\in X'$ be an index where $f_p(c_p^\top \theta)=\min_{p'\in X'}\{f_{p'}(c_{p'}^\top \theta)\}$. Then since $|X'|= \lambda'$, we have $$ \sum_{i=1}^m z_i f_i(c_i^\top \theta)\geq \sum_{p'\in X'}z_{p'} f_{p'}(c_{p'}^\top \theta)\geq \lambda'f_p(c_p^\top \theta).$$ On the other hand,  $\bar{z}$ is obtained by rounding $z$ down. Notice that $z_{p'}\in\{0,1\}$ for all $p'\in X'\cup\hat{X'}$, that is, for all $p'$ such that $f_{p'}(c_{p'}^\top \theta)>f_p(c_p^\top \theta)$. Therefore for all $p'$ such that $f_{p'}(c_{p'}^\top \theta)>f_p(c_p^\top \theta)$, we have $z_{\ell'}=\bar{z}_{\ell'}$. By Lemma \ref{lemma: constant fractional components} $z$ has at most $2r_++2$ fractional components. Therefore, because $\zeta\leq \sum_{i=1}^m z_i f_i(c_i^\top \theta)\leq k\max_{i\in[m]}\{f_i((Vu)_{\max})\}$, \begin{eqnarray*}
    \sum_{i=1}^m \bar{z}_i f_i(c_i^\top \theta)&\geq& \sum_{i=1}^m z_i f_i(c_i^\top \theta) - (2r_++2)f_p(c_p^\top \theta)\\&\geq& \sum_{i=1}^m z_i f_i(c_i^\top \theta)-(2r_++2) \sum_{i=1}^m z_i f_i(c_i^\top \theta)/\lambda'\\
    &\geq& \zeta-(2r_++2)k\max_{i\in[m]}\{f_i((Vu)_{\max})\}/\lambda'\\
    &\geq&\zeta-\delta_2.
\end{eqnarray*}

\end{myproof}

Below we give the pseudo-code of the oracle described in Lemma~\ref{lemma: auxiliary oracle on each tuple is enough}.

\begin{algorithm}[H]
\caption{Phase Two: Approximation of Linearized Auxiliary Problems}
\label{alg: phase two LP rounding}
\DontPrintSemicolon
\SetAlgoLined

\KwIn{Instance of Problem \ref{eqn:P(X,t)}, thresholds $\theta \in \mathbb{R}^{r_+}$, $\zeta \in \mathbb{R}_+$, parameters $\delta_1, \delta_2 > 0$}
\KwOut{One of two scenarios:
\begin{enumerate}
    \item \textbf{Scenario 1:} Certify no feasible $x$ to $\mathsf{P}(X,t)$ satisfies $d_j^\top x \geq \theta_j$ $\forall j \in [r_+]$ and $\sum_{i=1}^m x_i f_i(c_i^\top \theta) \geq \zeta$
    \item \textbf{Scenario 2:} Return feasible $x$ to $\mathsf{P}(X,t)$ with $d_j^\top x + \delta_1 \geq \theta_j$ $\forall j \in [r_+]$ and $\sum_{i=1}^m x_i f_i(c_i^\top \theta) + \delta_2 \geq \zeta$
\end{enumerate}}

\BlankLine
Set $\lambda' > 0$ according to Eq.~\eqref{eqn: lambda'}\;
Let $\mathcal{M} \leftarrow [m] \setminus (\cup_j X_j) \cup (\cup_j \hat{X}_j)$\;

\BlankLine
\tcp{Check small solutions}
\For{each $X' \subset \mathcal{M}$ with $0 \leq |X'| \leq \lambda' - 1$}{
    Set $\hat{X'}$ according to Eq.~\eqref{eqn: X'}\;
    Define $z \in \{0,1\}^m$: $z_i = 1$ if $i \in X' \cup (\cup_j X_j)$, else $z_i = 0$\;
    
    \If{$z$ feasible and satisfies thresholds with slack $(\delta_1, \delta_2)$}{
        \Return{Scenario 2 with solution $z$}
    }
}

\BlankLine
\tcp{Solve via LP rounding for larger solutions}
\For{each $X' \subset \mathcal{M}$ with $|X'| = \lambda'$}{
    Construct polyhedron $PH$ according to Eq.~\eqref{eqn: PH}\;
    
    \If{$PH \neq \emptyset$}{
        Choose arbitrary $z^* \in PH$\;
        Construct $PH(z^*)$ according to Eq.~\eqref{eqn: PH(z)}\;
        Compute vertex $y^*$ of $PH(z^*)$ with $\leq 2r_+ + 2$ fractional components\;
        
        Set $I_y \leftarrow \mathcal{M} \setminus (X' \cup \hat{X'})$\;
        Define $z \in [0,1]^m$: $z_i = \begin{cases}
            1 & \text{if } i \in (\cup_j X_j) \cup X'\\
            0 & \text{if } i \in (\cup_j \hat{X}_j) \cup \hat{X'}\\
            y^*_i & \text{if } i \in I_y
        \end{cases}$\;
        
        \Return{Scenario 2 with solution $z$}
    }
}

\Return{Scenario 1}

\end{algorithm}

\subsection{Completing the Proof}

Finally, we analyze our algorithm's overall performance and runtime. Let $\delta_1=\epsilon$ and $\delta_2=k\epsilon$, and in order to apply Lemma \ref{lemma: z satisfies auxiliary oracle}, we set $\lambda=(2r_++2)(Vu)_{\max}/\epsilon$ and $\lambda'=(2r_++2)\max_{i\in[m]}\{f_i((Vu)_{\max})\}/\epsilon$. 
We will treat the performance guarantee and runtime analysis separately.

\paragraph{Performance Guarantee:} Let $$c_{\epsilon,\gamma,W'_{\min}}=\frac{(1+\gamma)\epsilon}{W'_{\min}}.$$
The algorithm $\textup{ALG}$ (for solving $\mathsf{P}(X,t)$) in Lemma \ref{lemma: auxiliary oracle on each tuple is enough} satisfies \begin{align*}
    \textup{ALG}_{\mathsf{P}(X,t)}&\geq\left(1-g\left(\frac{\delta_1(1+\gamma)}{W'_{\min}}\right)\right)(\textup{OPT}_{\mathsf{P}(X,t)}-2\delta_2)-kh\left(\frac{\delta_1(1+\gamma)}{W'_{\min}}\right)\\&=(1-g(c_{\epsilon,\gamma,W'_{\min}}))\textup{OPT}_{\mathsf{P}(X,t)}-k(2\epsilon(1-g(c_{\epsilon,\gamma,W'_{\min}}))+h(c_{\epsilon,\gamma,W'_{\min}}))
    .
\end{align*} This gives the $\textup{ALG}'$ (for solving $\mathsf{P}(X,t)$) in Lemma \ref{lemma: auxiliary oracle is enough} with $$g'(\epsilon)=g(c_{\epsilon,\gamma,W'_{\min}})$$ and $$h'(\epsilon)=k(2\epsilon(1-g(c_{\epsilon,\gamma,W'_{\min}}))+h(c_{\epsilon,\gamma,W'_{\min}})).$$ 
Therefore, the algorithm $\textup{ALG}$ (for solving $\mathsf{P}(X)$) in Lemma \ref{lemma: auxiliary oracle is enough} satisfies \begin{align*}
    \textup{ALG}_{\mathsf{P}(X)}&\geq (1-g'(\epsilon))((1-g(c_{\epsilon,\gamma,W'_{\min}}))\textup{OPT}_{\mathsf{P}(X)}+kh(c_{\epsilon,\gamma,W'_{\min}}))-h'(\epsilon)\\
    &= (1-g(c_{\epsilon,\gamma,W'_{\min}}))((1-g(c_{\epsilon,\gamma,W'_{\min}}))\textup{OPT}_{\mathsf{P}(X)}+kh(c_{\epsilon,\gamma,W'_{\min}}))\\&\quad\text{ }-k(2\epsilon(1-g(c_{\epsilon,\gamma,W'_{\min}}))+h(c_{\epsilon,\gamma,W'_{\min}}))\\
    &=(1-g(c_{\epsilon,\gamma,W'_{\min}}))^2\textup{OPT}_{\mathsf{P}(X)}-k(2\epsilon(1-g(c_{\epsilon,\gamma,W'_{\min}}))+g(c_{\epsilon,\gamma,W'_{\min}})h(c_{\epsilon,\gamma,W'_{\min}})).
\end{align*}
Therefore, the algorithm $\textup{ALG}$ (for solving $\mathsf{P}'(I)$) in Lemma \ref{lemma: enumerate valid tuples}  satisfies $$\textup{ALG}_{\mathsf{P}'(I)}\geq (1-g(c_{\epsilon,\gamma,W'_{\min}}))^2\textup{OPT}_{\mathsf{P}'(I)}-k(2\epsilon(1-g(c_{\epsilon,\gamma,W'_{\min}}))+g(c_{\epsilon,\gamma,W'_{\min}})h(c_{\epsilon,\gamma,W'_{\min}})).$$
Finally, we apply Lemma \ref{lemma: replace W'} by plugging in $1-\alpha= (1-g(c_{\epsilon,\gamma,W'_{\min}}))^2$ and $\beta=-k(2\epsilon(1-g(c_{\epsilon,\gamma,W'_{\min}}))+g(c_{\epsilon,\gamma,W'_{\min}})h(c_{\epsilon,\gamma,W'_{\min}}))$. This gives \begin{align*}
    \mathrm{ALG}_{\mathsf{P}(I)}
    &\geq (1-g(2\gamma(Vu)_{\max}))^2(1-g(c_{\epsilon,\gamma,W'_{\min}}))^2\textup{OPT}_{\mathsf{P}(I)}\\&\quad\text{ }-k(1-g(2\gamma(Vu)_{\max}))(2\epsilon(1-g(c_{\epsilon,\gamma,W'_{\min}}))+g(c_{\epsilon,\gamma,W'_{\min}})h(c_{\epsilon,\gamma,W'_{\min}})\\&\quad\text{ }+(1-g(c_{\epsilon,\gamma,W'_{\min}}))^2h(2\gamma(Vu)_{\max})+h(2\gamma(Vu)_{\max})).
\end{align*}

\paragraph{Runtime Analysis:} By Lemma \ref{lemma: constant fractional components}, we give an oracle described in Observation \ref{obs: partition value} with runtime $$T_{\textup{LP}}:=\textup{LP}(m,3m+r_++2)+\textup{LP}(m,2m+2r_++2).$$ 
Therefore, by Observation \ref{obs: partition value}, we give an oracle 
described in Lemma \ref{lemma: auxiliary oracle on each tuple is enough} with runtime $$\lambda'm^{\lambda'}T_{\textup{LP}}.$$
Therefore, the algorithm $\textup{ALG}'$ (for solving $\mathsf{P}(X,t)$) in Lemma \ref{lemma: auxiliary oracle is enough} has runtime
\begin{align*}
    &\quad\text{ }\lceil ((Vu)_{\max}-\min\{0,(Vu)_{\min}\})/\delta_1\rceil^{r_+}\cdot\left\lceil k\max_{i\in[m]}\{f_i((Vu)_{\max})\}/\delta_2\right\rceil\cdot\lambda'm^{\lambda'}T_{\textup{LP}}\\&=\lceil ((Vu)_{\max}-\min\{0,(Vu)_{\min}\})/\epsilon\rceil^{r_+}\cdot\left\lceil \max_{i\in[m]}\{f_i((Vu)_{\max})\}/\epsilon\right\rceil\cdot\lambda'm^{\lambda'}T_{\textup{LP}}.
\end{align*} 
Therefore, the algorithm $\textup{ALG}'$ (for solving $\mathsf{P}(X)$) in Lemma \ref{lemma: enumerate valid tuples} has runtime \begin{align*}
    \left\lceil\left(\frac{4(1+\gamma)k}{\epsilon W'_{\min}}\right)^{r_+}\right\rceil \cdot \lceil ((Vu)_{\max}-\min\{0,(Vu)_{\min}\})/\epsilon\rceil^{r_+}\cdot\left\lceil \max_{i\in[m]}\{f_i((Vu)_{\max})\}/\epsilon\right\rceil\cdot\lambda'm^{\lambda'}T_{\textup{LP}}.
\end{align*} 
Finally, by Lemma \ref{lemma: enumerate valid tuples}, our algorithm's runtime is
\begin{align*}
    &\quad \text{ } r_+ m \log_2 m +\lambda m^\lambda+ \lambda r_+^2m^{\lambda r_+}\\&+\left\lceil\left(\frac{4(1+\gamma)k}{\epsilon W'_{\min}}\right)^{r_+}\right\rceil \cdot \left\lceil \frac{(Vu)_{\max}-\min\{0,(Vu)_{\min}\}}{\epsilon}\right\rceil^{r_+}\cdot\left\lceil \frac{\max_{i\in[m]}\{f_i((Vu)_{\max})\}}{\epsilon}\right\rceil\cdot\lambda'm^{\lambda r_++\lambda'}T_{\textup{LP}}.
\end{align*} 

Below we give the pseudo-code of our phase two algorithm in Proposition~\ref{prop: phase two}.

\begin{algorithm}[H]
\caption{Phase Two}
\label{alg: phase two}
\DontPrintSemicolon
\SetAlgoLined

\KwIn{Attention matrix $W=\textup{softmax}(QK^\top)\in\mathbb{R}_+^{n\times n}$, low-rank approximation $W'\in\mathbb{R}_+^{n\times n}$ with factorization $W'=AB^\top$ and element wise guarantee $1-\gamma \le W_{ij}/W'_{ij} \le 1+\gamma$, value matrix $V\in\mathbb{R}^{n\times d_v}$, user vector $u\in\mathbb{R}^{d_v}$, maximum number of recommended items $k$, candidate index set $I$, parameter $\epsilon>0$.}
\KwOut{Solution $x$ to Problem \eqref{eqn:problem}.}

\BlankLine

Construct the instance of Problem \eqref{eqn:problem2} from the inputs\;

\BlankLine
\tcp{Low Non-negative Rank Approximation}
Form the instance of Problem \eqref{eqn:P'(I)} by replacing $W$ with $W'$ in \eqref{eqn:problem2}\;


\BlankLine
\tcp{Enumeration of Partial Solutions}
$x \leftarrow$ \textbf{Run} Algorithm \ref{alg: phase two enumeration} on instance \eqref{eqn:P'(I)} with parameter $\epsilon$, which internally invokes:\;
\Indp

\BlankLine
\tcp{Discretization of Auxiliary Problems}
\(\triangleright\) \textbf{Algorithm \ref{alg: phase two discretization}} with \textit{Input:} Instance \eqref{eqn:P(x)}, parameter $\epsilon$; \textit{Calls:} Algorithm \ref{alg: phase two linearization}\; 

\BlankLine
\tcp{Complete Linearization of Auxiliary Problems}
\(\triangleright\) \textbf{Algorithm \ref{alg: phase two linearization}} with \textit{Input:} Instance $\mathsf{P}(X,t_{\ell'})$, parameters $\delta_1\leftarrow \epsilon$, $\delta_2\leftarrow k\epsilon$; \textit{Calls:} Algorithm \ref{alg: phase two LP rounding}\;

\BlankLine
\tcp{Approximation of Linearized Auxiliary Problems via LP Rounding}
\(\triangleright\) \textbf{Algorithm \ref{alg: phase two LP rounding}} with \textit{Input:} Instance $\mathsf{P}(X,t)$, thresholds $(\theta_1,\ldots,\theta_{r_+},\zeta)$, parameters $\delta_1\leftarrow \epsilon$, $\delta_2\leftarrow k\epsilon$\;
\Indm

\BlankLine
\Return{$x$}
\end{algorithm}

\end{document}